\documentclass[preprint,review]{elsarticle}

\usepackage{lineno}

\usepackage{amsmath,amssymb}
\usepackage{stmaryrd}
\usepackage{enumerate}
\usepackage{algorithm}
\usepackage[noend]{algpseudocode}
\usepackage{calrsfs}
\usepackage{physics}
\usepackage{bibentry}
\usepackage{graphicx}
\usepackage{amssymb}
\usepackage{array}
\usepackage{amsmath}
\usepackage{epstopdf}
\usepackage{epsfig}
\usepackage{subfigure}
\usepackage{pstricks}
\usepackage{fancyheadings}
\usepackage[english]{babel}
\usepackage{amsthm}

\usepackage{stackengine}

\usepackage[margin=1.0in]{geometry}
\DeclareMathAlphabet{\pazocal}{OMS}{zplm}{m}{n}

\usepackage{graphicx}

\usepackage{stackengine}

\makeatletter
\def\BState{\State\hskip-\ALG@thistlm}
\makeatother

\newcommand{\lJump}{[\![}
\newcommand{\rJump}{]\!]}

\newtheorem{theorem}{Theorem}
\newtheorem{lemma}{Lemma}
\newtheorem{remark}{Remark}

\newtheorem{Definition}{Definition}

\begin{document}
\begin{frontmatter}
\title{A stable discontinuous Galerkin method for linear elastodynamics in 3D geometrically complex elastic solids using physics based numerical fluxes}
\author[label1,label5]{Kenneth Duru}
\author[label4]{Leonhard Rannabauer }
\author[label2]{Alice-Agnes Gabriel}
\author[label2,label3]{On Ki Angel Ling}
\author[label2]{Heiner Igel}
\author[label4]{Michael Bader}
 \address[label1]{Mathematical Sciences Institute, The Australian National University, Canberra, Australia}
 \address[label2]{Department of Earth and Environmental Sciences, Ludwig-Maximilians-Universit{\"a}t M{\"u}nchen, Germany}
   \address[label3]{ETH Zurich, Switzerland}
 \address[label4]{Technical University of Munich, Germany}
 \address[label5]{Corresponding author: kenneth.duru@anu.edu.au}

\pagenumbering{arabic}

\begin{abstract}
High order accurate and explicit time-stable solvers are well suited for hyperbolic wave propagation problems.
As a result of the complexities of real geometries, internal interfaces and nonlinear boundary and interface conditions, discontinuities and sharp wave fronts may become fundamental features of the solution.
Thus, geometrically flexible and adaptive numerical algorithms are critical for high fidelity and efficient simulations of wave phenomena in many applications. 
%
Adaptive curvilinear meshes hold promise to minimise the effort to represent complicated geometries or heterogeneous material data avoiding the bottleneck of feature-preserving meshing. 
To enable the design of stable DG methods on three space dimensional (3D) curvilinear elements we construct a structure preserving anti-symmetric coordinate transformation motivated by the underlying physics.
Using a physics-based numerical penalty-flux, we develop a 3D provably energy-stable discontinuous Galerkin finite element approximation of the elastic wave equation in geometrically complex and heterogenous media. By construction, our numerical flux is upwind and yields a discrete energy estimate analogous to the continuous energy estimate.
The ability to treat conforming and non-conforming curvilinear elements  allows for flexible adaptive mesh refinement strategies.
The numerical scheme has been implemented in ExaHyPE, a simulation engine for parallel dynamically adaptive simulations of wave problems on adaptive Cartesian meshes.
 %
We present 3D numerical experiments of wave propagation in heterogeneous isotropic and anisotropic elastic solids demonstrating stability and high order accuracy. 
We demonstrate the potential of our approach for computational seismology in a regional wave propagation scenario in a geologically constrained 3D model including the geometrically complex free-surface topography of Mount Zugspitze, Germany.
\end{abstract}

\end{frontmatter}
\section{Introduction}\label{sec:s1}
Numerical algorithms based on the discontinuous Galerkin (DG) method \cite{ReedHill1973, CockburnShu1989, CockburnHouShu1990, HesthavenWarburton2008, HesthavenWarburton2002} have shown to be flexible, high order accurate, provably stable, and well suited for complex large scale wave propagation problems  \cite{HesthavenWarburton2002, DumbserKaeser2006}. 
In computational seismology, they have been successfully applied to extreme-scale simulations of grand-challenge scenarios, often exploiting the largest-available supercomputers \cite{Burstedde:2010,Breuer2014,Heineckeetal2014,Uphoff:2017}. 

In this study, we will pay special attention to scattering of linear elastic waves in heterogeneous, isotropic and anisotropic, solid Earth models with complex free surface topography. 
Accurate and efficient numerical simulation of seismic surface and interface waves, and scattering of high frequency waves by complex nonplanar topography are critical for assessing and quantifying seismic risks and hazards \cite{bielak2010shakeout,chaljub2010grenoble,Graves_etal2011}. 
Surface and interface waves \cite{Rayleigh1885} are often the largest amplitude waves modes. 
On a regional scale, physics-based 3D ground motion simulations which include realistic three-dimensional Earth structure and topography
are now able to resolve frequencies that are relevant for building response 
(from static displacements at zero frequency up to 5-10~Hz or higher, \cite[e.g.]{cui2013sc, rodgers2018hayward}), however, many computational studies are forced to use highly smoothed representations of surface topography due to limitations either in terms of algorithmic restrictions 
(e.g., spatial discretisations limited to structured and or hexahedral approaches, 
limited flexibility of free-surface boundary conditions, limited applicability of graded meshes) and/or computational efficiency.
For efficient representation of complex geometries, in this study, we use adaptive boundary conforming curvilinear elements. 
The elastic wave equation is transformed from Cartesian coordinates to curvilinear coordinates using a structure preserving coordinate transformation. The transformation is local within the element, and inside the element the curvilinear elements have logical Cartesian coordinates.  Essentially, we solve the equations on adaptive Cartesian meshes, and the complex geometries are moved into variable metric terms and are used to define transformed variable material parameters.
Our numerical scheme allows adaptive non-conforming refinement of the hexahedral elements, which enables its implementation on tree- or block-structured Cartesian meshes. 
We implemented the scheme in the ExaHyPE engine for solving hyperbolic PDE systems \cite{ExaHyPE2019}, which realises DG with ADER time stepping \cite{Dumbser:2013}.  The Peano adaptive mesh refinement framework \cite{Weinzierl2019, WeinzierlMehl2011} provides dynamically adaptive tree-structured Cartesian meshes and parallelisation in shared and distributed memory systems.

\paragraph{DG and physics-based numerical fluxes}
The DG method is an increasingly attractive method for approximating partial differential equations (PDEs), and has been successfully used in computational seismology applications including geometrically and rheologically complicated wave propagation and dynamic rupture simulations  \cite[e.g.,]{Puente2007,Puente2008,delaPuenteAmpueroKaser2009,PeltiesdelaPuenteAmpueroBrietzkeKaser2012,Wolf2020}.
An important component of the DG method is the numerical flux  \cite{DeGraziaMengaldoMoxeyVincentSherwin2013,Huynh2007}.
For hyperbolic PDEs this flux is based on approximate or exact solutions of the Riemann problem \cite{Godunov1959, Rusanov1961}. 
The choice of a numerical flux and (approximate) solutions of the Riemann problem are critical for accuracy and stability of the DG method \cite{Qiu2008, KirbyKarniadakis2005, KoprivaNordstromGassner2016}. 
The Rusanov flux \cite{Rusanov1961} (also called local Lax-Friedrichs flux) is widely used, because of its simplicity and robustness. However, the Rusanov flux might be inappropriate for simulating seismic surface waves \cite{DuruGabrielIgel2017}. 
The Godunov flux with the exact solution to the Riemann problem has been demonstrated as a reasonable choice for seismological applications \cite{KaeserDumbser2006,Pelties2014}. However, the Godunov flux requires a complete eigenvector and eigenvalue decomposition of the coefficient matrices of the spatial operator. In general anisotropic media, the eigen-decomposition can be nontrivial. The eigen-decomposition will  become even more cumbersome if elements are curved. In \cite{ChangWarburton2017, Wilcox2010} numerical penalty fluxes are introduced, penalizing the normal components of the primitive variables across inter-element faces, and avoiding the eigen-decomposition of coefficient matrices. However, many dynamic boundary conditions such as empirical constitutive friction laws \cite{Scholz1998, Rice1983, JRiceetal_83} are formulated using derived quantities, such as tractions in local coordinates. We argue that it will be more natural to anchor numerical fluxes in elastic solids on derived quantities such as local velocity and traction vectors, rotated into local orthogonal coordinates. This will enable the development of a unified provably stable and robust adaptive DG framework in complex geometries for the numerical treatment of 1) nonlinear frictional sliding in elastic solids, 2) for coupling classical DG inter-element interfaces in elastic solids where slip is not permitted, and 3) numerical enforcement of external well-posed boundary conditions modeling various solid mechanics and geophysical phenomena. 

In \cite{DuruGabrielIgel2017}, we introduce (in 1D and 2D) an alternative and accurate approach to couple locally adjacent DG elements, using physical conditions such as friction. The nascent physics-based numerical flux obeys the eigen-structure of the PDE and the underlying physics at the internal and external DG element boundaries in a provably stable manner. Our formulation does not require a complete eigenvector-eigenvalue decomposition of the spatial coefficient matrices and can be easily adapted to model linear and nonlinear boundary and interface wave phenomena. In spirit, our approach is analogous to the method used in a finite difference framework \cite{DuruandDunham2016} to model frictional sliding during dynamic earthquake rupture. However, static and/or dynamic adaptive mesh refinement in a finite difference setting is arduous. 

\paragraph{Extension to 3D curvilinear elements}
In this study, we will extend the physics-based numerical flux \cite{DuruGabrielIgel2017} to 3D geometrically complex elastic media. We extract tractions and particle velocities on an element face in Cartesian coordinates and rotate them into the curvilinear coordinates. The rotated tractions and particle velocities, in conjunction with the material impedance, can be used to extract characteristics, such as plane shear waves and compressional waves,  propagating along the boundary surface. Then, we construct boundary and interface data by solving a Riemann-like problem and constrain the solutions against the physical conditions acting at the element faces. The physics based fluctuations are constructed by penalisation of boundary and interface data against the incoming characteristics. These fluctuations are then appended to the discrete equations with physically motivated penalty weights chosen such that the semi-discrete approximation satisfies an energy estimate analogous to the continuous energy estimate. The energy estimate proves the asymptotic stability of the semi-discrete approximation.

The semi-discrete DG approximation is integrated in time using the Arbitrary DERivative (ADER) time integration \cite{Toro1999}. The ADER time discretisation is summarised in \ref{sec:ADER}. The cell-wise local character of DG  can be readily combined with an ADER scheme leading to high-order accuracy in time within a single-step. The combination of the DG approximation in space and the ADER time integration is often referred to as the ADER-DG scheme \cite{KaeserDumbser2006}. For the ADER-DG scheme, the numerical flux fluctuation is evaluated only once for any order of accuracy. The implication is that most of the computations are performed within the element to evaluate spatial derivatives. 

 We will here present 3D ADER-DG numerical simulations verifying accuracy and stability of the method, using community developed benchmark problems \cite{Kristekova_etal2006, Kristekova_etal2009, Favretto-Cristini_etal2011, Komatitsch_etal2011} and complex geologically constrained geometries. All of our simulations employ provable stable PML boundary conditions \cite{DuruRannabauerGabrielKreissBader2019, Duru2019}  to effectively prevent artificial boundary conditions from contaminating the simulations.  A stable numerical implementation of the PML for  3D linear elastodynamics is nontrivial and allows the generation of high quality seismograms. We will perform error analysis and compute error parameters relevant to computational seismology  \cite{Kristekova_etal2006, Kristekova_etal2009}.



The remaining part of the paper will proceed as follows. In the next section, we introduce a general model for linear elastodynamics. Curvilinear coordinates and structure preserving curvilinear transformations are presented in section 3. In section 4, we present physical boundary and interface conditions and derive energy estimates. In section 5, the physics based numerical flux and the algorithms for solving the Riemann problem in heterogeneous media and arbitrary curvilinear coordinates are presented. Numerical experiments based on the implementation in ExaHyPE are presented in section 8 verifying accuracy and stability, and demonstrating the potentials of the method in simulating complex wave phenomena. In the last section we draw conclusions and suggest future work.
\section{First order linear hyperbolic PDE}
Consider the   3D  first order linear  hyperbolic system in a source--free heterogeneous medium
{
\small
\begin{equation}\label{eq:linear_wave}
\begin{split}
\mathbf{P}^{-1}\frac{\partial{\mathbf{Q}}}{\partial t} = \sum_{\xi = x,y,z}\mathbf{A}_{\xi}\frac{\partial{\mathbf{Q}}}{\partial \xi},
  \end{split}
  \end{equation}
}
where $\mathbf{P} = \mathbf{P}^T$ with $ \mathbf{Q}^T\mathbf{P}\mathbf{Q} > 0$ and $\mathbf{A}_\xi = \mathbf{A}_\xi^T$, $\xi = x, y, z$.
Here, $t \ge 0$ denotes time and $(x,y,z)\in \Omega $ are the Cartesian coordinates of the spatial domain $ \Omega \subset \mathbb{R}^{3}$.
At the initial time, $t = 0$, we set the initial condition
\begin{align}\label{eq:initial_data}
\mathbf{Q}(x,y,z,0) = \mathbf{Q}^0(x,y,z) \in L^2(\Omega),
\end{align}
belonging to the  space of square integrable functions. 
In general, the symmetric positive definite matrix $\mathbf{P}$ depends on the spatial coordinates $x, y, z$, and encodes the material parameters of the underlying medium. 
The  non-dimensional constant  matrices $\mathbf{A}_\xi$  encapsulate the underlying linear conservation law and the corresponding linear constitutive relation. 

We introduce the energy density and the physical energy,  defined by
{
\small
\begin{equation}\label{eq:physical_energy}
\frac{dE(x,y,z,t)}{dxdydz} = \frac{1}{2}\left[\mathbf{Q}^T\mathbf{P}^{-1}\mathbf{Q}\right]  > 0, \quad  E(t) = \int_{\Omega}  dE(x,y,z,t)  > 0.
\end{equation}
}
Let $\Gamma$ denote the boundary of the domain $\Omega$,  and $\mathbf{n} = \left(n_x, n_y, n_z\right)^T$ the outward unit normal on the boundary.  We multiply   \eqref{eq:linear_wave}  with $\boldsymbol{\phi}^T(x,y,z) $ from the left, where $\boldsymbol{\phi}(x,y,z) \in L^2\left(\Omega\right)$ is an arbitrary test function, and integrate over the whole spatial domain, $\Omega$.  We have
{
\small
\begin{align}\label{eq:product_1}
\int_{\Omega}\boldsymbol{\phi}^T\mathbf{P}^{-1}\frac{\partial{\mathbf{Q}}}{\partial t}dxdydz &= \int_{\Omega}\boldsymbol{\phi}^T\left(\sum_{\xi = x,y,z}\mathbf{A}_{\xi}\frac{\partial{\mathbf{Q}}}{\partial \xi}\right)dxdydz.
\end{align}
}
In the right hand side of \eqref{eq:product_1}, integrating-by-parts, and using the fact that the coefficient matrices are constant and symmetric, $\mathbf{A}_\xi  = \mathbf{A}_\xi^T$, gives
{
\small
\begin{align}\label{eq:product_2}
\int_{\Omega}\boldsymbol{\phi}^T\mathbf{P}^{-1}\frac{\partial{\mathbf{Q}}}{\partial t}dxdydz 
& = \frac{1}{2}\int_{\Omega}\left(\sum_{\xi = x,y,z}\left[\boldsymbol{\phi}^T\mathbf{A}_{\xi}\frac{\partial{\mathbf{Q}}}{\partial {\xi}}-\mathbf{Q}^T\mathbf{A}_{\xi}\frac{\partial{\boldsymbol{\phi}}}{\partial {\xi}}\right]  \right)dxdydz + \frac{1}{2}\oint_{\Gamma}\boldsymbol{\phi}^T \left(\sum_{\xi = x, y, z}n_{\xi}\mathbf{A}_{\xi}\right)\mathbf{Q} dS.
\end{align}
}
 Replacing $\boldsymbol{\phi}$ with $\mathbf{Q}$ in \eqref{eq:product_2}, in the right hand side,  the volume terms vanish, having
{
\small
\begin{align}\label{eq:product_3}
\int_{\Omega}\mathbf{Q}^T\mathbf{P}^{-1}\frac{\partial{\mathbf{Q}}}{\partial t}dxdydz = \frac{1}{2}\oint_{\Gamma}\mathbf{Q}^T \left(\sum_{\xi = x, y, z}n_{\xi}\mathbf{A}_{\xi}\right)\mathbf{Q} dS.
\end{align}
}
 The decay condition, $|\mathbf{Q}| \to 0$ at  the boundary $(x,y,z) \in \Gamma $, yields the energy equation
{
\small
\begin{align}\label{eq:energy_equation}
\frac{d}{dt}E(t) = 0.
\end{align}
}
The energy is conserved, $E(t) = E(0)$ for all $t\ge 0$. 



Depending on the coefficient matrices $\mathbf{P}, \mathbf{A}_\xi$ the  system \eqref{eq:linear_wave} can describe acoustic waves, electromagnetic waves, linear MHD waves, linear elastic waves, or the interaction of acousto--elastic waves propagating in a heterogeneous medium. Our general  interest is the development of  a provably stable DG  approximations for first order  linear hyperbolic PDEs of the form \eqref{eq:linear_wave},  in heterogeneous and geometrically complex elastic solids, subject to well-posed boundary conditions. 

In this study, we will pay special attention to scattering of linear elastic waves in heterogeneous, isotropic and anisotropic,  elastic solids   with complex free surface topography.
We note, that with limited modifications the method we develop can be extended to other problems modeled by \eqref{eq:linear_wave} in future work.

\subsection{Waves in elastic solids}
To describe wave propagation in elastic solids, we introduce the unknown wave fields 
\begin{align}\label{eq:velocity_stress}
\mathbf{Q}\left(x,y,z,t\right) = \begin{bmatrix}
\mathbf{v}(x,y, z,t) \\
\boldsymbol{\sigma}(x,y,z,t)
 \end{bmatrix},
\end{align}
with the particle velocity vector, $\mathbf{v}(x,y, z,t) = \left[ v_x, v_y, v_z \right]^T$, and the stress vector, \\
$\boldsymbol{\sigma}(x,y,z,t) = \left[ \sigma_{xx}, \sigma_{yy}, \sigma_{zz}, \sigma_{xy},  \sigma_{xz},  \sigma_{yz}\right]^T$.  
The symmetric constant coefficient matrices $\mathbf{A}_\xi $ describing the conservation of momentum and the constitutive relation, defined by Hooke's law, are given by
%
{
\footnotesize
\begin{align}\label{eq:elastic_coeff}
\mathbf{A}_{\xi} = 
\begin{pmatrix}
\mathbf{0}_3 & \mathbf{a}_{\xi}\\
\mathbf{a}_{\xi}^T & \mathbf{0}_6
\end{pmatrix},
\quad
\mathbf{a}_x = 
\begin{pmatrix}
1& 0& 0& 0&0& 0\\
0& 0& 0& 1&0& 0\\
0& 0& 0& 0&1& 0
\end{pmatrix},
\quad
\mathbf{a}_y = 
\begin{pmatrix}
0& 0& 0& 1&0& 0\\
0& 1& 0& 0&0& 0\\
0& 0& 0& 0&0& 1
\end{pmatrix},
\quad
\mathbf{a}_z = 
\begin{pmatrix}
0& 0& 0& 0&1& 0\\
0& 0& 0& 0&0& 1\\
0& 0& 1& 0&0& 0\\
\end{pmatrix},
\end{align}
}
where $\mathbf{0}_3$ and $\mathbf{0}_6$ are  the $3$-by-$3$ and $6$-by-$6$ zero  matrices.

The symmetric positive definite material parameter matrix $\mathbf{P}$  is defined by
{
\footnotesize
\begin{align}\label{eq:material_coeff}
\mathbf{P} = 
\begin{pmatrix}
\rho^{-1} \mathbf{1}  & \mathbf{0}\\
 \mathbf{0}^T  & \mathbf{C}
\end{pmatrix}
 ,
  \quad
  \mathbf{1} =   \begin{pmatrix}
  1 & 0 & 0 \\
  0 & 1 & 0 \\
  0 & 0 & 1 
  \end{pmatrix}
,
  \quad
  \mathbf{0} =   \begin{pmatrix}
  0 & 0 & 0& 0 & 0 & 0 \\
  0 & 0 & 0& 0 & 0 & 0 \\
  0 & 0 & 0& 0 & 0 & 0 
  \end{pmatrix}
  ,
\end{align}
}
where  $\rho(x,y,z) > 0$ is the  density of the medium, and  $\mathbf{C} = \mathbf{C}^{T} > 0$ is the symmetric positive definite  matrix of elastic constants.
With the unknown wave fields  prescribed by \eqref{eq:velocity_stress}, and  the coefficient matrices defined in  \eqref{eq:elastic_coeff} and \eqref{eq:material_coeff}, the first three equations in \eqref{eq:linear_wave} are the conservation of momentum and the last six equations are the time derivative of the constitutive relation, defined by Hooke's law, relating  stress fields to strains where the constant of proportionality is the stiffness matrix of elastic coefficients $\mathbf{C}$.

The mechanical energy density is the sum of the kinetic energy density and the strain energy density
\begin{align}\label{eq:elastic_energy_density}
\frac{dE}{dxdydz} := \frac{1}{2}[\mathbf{Q}^T\mathbf{P}^{-1}\mathbf{Q}] = \frac{\rho}{2}\left(v_x^2 + v_y^2 + v_z^2\right) + \frac{1}{2}\boldsymbol{\sigma}^T\mathbf{S}\boldsymbol{\sigma} > 0,
\end{align}
where $\mathbf{S} = \mathbf{C}^{-1}$ is the compliance matrix.

In a general anisotropic medium the stiffness matrix $\mathbf{C}$ is described by 21 independent elastic coefficients. 
In an orthotropic anisotropic  medium  the stiffness matrix has 9 independent elements,
\begin{align}\label{eq:orthotropic_stiffness_tensor}
\mathbf{C} =   \begin{pmatrix}
  c_{11} & c_{12} & c_{13} & 0 & 0 & 0\\
  c_{12} & c_{22} & c_{23}& 0  & 0 & 0\\
  c_{13} & c_{23} & c_{33}& 0 &  0 & 0\\
  0 & 0 & 0& c_{44} & 0 & 0 \\
  0 & 0 & 0 & 0 & c_{55} & 0\\
  0 & 0 & 0 & 0 & 0 & c_{66}
  \end{pmatrix}.
\end{align}
 In the isotropic case, the medium is described by two independent elastic coefficients, the Lam\'e parameters $\mu > 0$, $\lambda > -\mu $, thus we have
\begin{align}\label{eq:isotropic_stiffness_tensor}
\mathbf{C} =   \begin{pmatrix}
  2\mu + \lambda &  \lambda &  \lambda & 0 & 0 & 0\\
   \lambda & 2\mu + \lambda &  \lambda & 0  & 0 & 0\\
   \lambda &  \lambda & 2\mu + \lambda & 0 &  0 & 0\\
  0 & 0 & 0& \mu  & 0 & 0 \\
  0 & 0 & 0 & 0 & \mu & 0\\
  0 & 0 & 0 & 0 & 0 & \mu 
  \end{pmatrix}.
\end{align}
However, as shown in \cite{MarsdenHughes1994} (pages 241--243), for strong ellipticity  we must  have $\mu > 0$, $\lambda > -1/2\mu $. 

We introduce the 3D canonical  basis vectors
\begin{align}\label{eq:canonical_basis}
\mathbf{e}_x = \left(1, 0, 0\right)^T,  \quad \mathbf{e}_y = \left(0, 1, 0\right)^T,  \quad \mathbf{e}_z = \left(0, 0, 1\right)^T,
\end{align}
and define the  Cartesian components of the velocity and traction vectors
\begin{align}\label{eq:velocity_tractions}
\mathbf{v} = \begin{pmatrix}
v_{x}\\
v_{y}\\
v_{z}
\end{pmatrix},
 \quad
\mathbf{T}^{(\xi)} = \begin{pmatrix}
T_{x}^{(\xi)} \\
T_{y}^{(\xi)} \\
T_{z}^{(\xi)} 
\end{pmatrix} = \bar{\bar{\sigma}}\mathbf{e}_\xi =  \boldsymbol{a}_{\xi}\bar{\bar{\sigma}}, \quad \xi = x, y, z, 
\quad
 \bar{\bar{\sigma}} =
\begin{pmatrix}
\sigma_{xx}&\sigma_{xy}&\sigma_{xz}\\
\sigma_{xy}&\sigma_{yy}&\sigma_{yz}\\
\sigma_{xz}&\sigma_{yz}&\sigma_{zz}
\end{pmatrix}.
\end{align}
Next, we define the coefficient matrices 
\begin{align}\label{eq:coeff}
\widetilde{\mathbf{A}}_\xi = {\mathbf{P}}{\mathbf{A}}_\xi.
\end{align}
In the absence of boundaries and discontinuous interfaces, the elastic wave equation supports two families of solutions, primary (p-waves) and secondary (s-waves) waves.
The p-wave and s-wave modes are related to the nontrivial eigenfunctions of $\widetilde{\mathbf{A}}_\xi$, and the 
 nontrivial eigenvalues of $\widetilde{\mathbf{A}}_\xi$ are
$\pm c_{p\xi},  \quad \pm c_{sh\xi}, \quad \pm c_{sv\xi},$
and correspond to the  p-wave and  s-wave speeds.
Note that $c_{p\xi} > 0$ are the p-wave speeds, $c_{sh\xi} > 0$ are the wave speeds of the horizontally polarized s-wave and $c_{sv\xi} > 0$ are the wave speeds of the vertically polarized s-wave.
The negative and positive going  p-wave and s-wave modes are given by
\begin{align}\label{eq:eigenfunction_elatistic}
 Z_{p\xi} v_\xi \mp T_\xi^{(\xi)}, \quad  Z_{sh\xi} v_{\eta} \mp T_{\eta}^{(\xi)} ,\quad  Z_{sv\xi} v_{\theta} \mp T_{\theta}^{(\xi)}, 
\end{align}
where $ Z_{p\xi} = \rho c_{p\xi}$, $Z_{sh\xi} = \rho  c_{sh\xi}$, $Z_{sv\xi} = \rho  c_{sv\xi}$ are the impedances.
Here, $\xi, \eta, \theta = x, y, z$,  $\eta \ne \xi$,  and $\theta \ne \xi, \eta$.
The wave modes defined in \eqref{eq:eigenfunction_elatistic} are the plane p-waves and plane s-waves propagating along the ${\xi}$-axis, and are related to the 1D Riemann invariant.
The eigenvalues of the matrices $\widetilde{\mathbf{A}}_\xi$ can be easily determined. For example in orthotropic anisotropic media, 
with $\mathbf{C}$ defined in \eqref{eq:orthotropic_stiffness_tensor}, the eigenvalues are given by
\begin{align}\label{eq:anisotropic_wavespeed}
c_{px} = \sqrt{\frac{c_{11}}{\rho}} ,\quad c_{shx} = \sqrt{\frac{c_{44}}{\rho}}, \quad c_{svx} = \sqrt{\frac{c_{55}}{\rho}},
\\
\nonumber
c_{py} = \sqrt{\frac{c_{22}}{\rho}} ,\quad c_{shy} = \sqrt{\frac{c_{66}}{\rho}}, \quad c_{svy} = \sqrt{\frac{c_{44}}{\rho}},
\\
\nonumber
c_{pz} = \sqrt{\frac{c_{33}}{\rho}} ,\quad c_{shz} = \sqrt{\frac{c_{55}}{\rho}}, \quad c_{svz} = \sqrt{\frac{c_{66}}{\rho}}.
\end{align}
Note that for the same wave mode the wave speed can vary in all directions. 
In isotropic media, the eigenvalues are uniform in all directions
\begin{align}\label{eq:isotropic_wavespeed}
c_{p\xi} = c_p = \sqrt{\frac{\lambda + 2\mu}{\rho}} ,\quad c_{sh\xi} =  c_{sv\xi} = c_s =\sqrt{\frac{\mu}{\rho}}.
\end{align}
Thus, in isotropic media a wave mode propagates with a uniform wave speed in all directions. In particular, the vertically polarized and the horizontally polarized  s-waves have identical wave-speed $c_{sh\xi} =  c_{sv\xi} = c_s$ in all directions. 

We recall our main interest being the development of an energy stable DG approximation of the equation of motion \eqref{eq:linear_wave} defined by \eqref{eq:velocity_stress},  \eqref{eq:elastic_coeff} and \eqref{eq:material_coeff}, in heterogeneous and geometrically complex elastic solid with complicated nonplanar free-surface topography.
%
To this end, we will next introduce the anti-symmetric split form. This is necessary for the development of provably stable approximations on curvilinear meshes for complex geometries and generally arbitrary heterogeneous media. 
\subsection{Anti-symmetric splitting}
To enable effective numerical treatments, we introduce the split form of  the equation of motion  \eqref{eq:linear_wave} defined by
  \begin{align}\label{eq:gen_hyp_antisymmetry}
{\mathbf{P}}^{-1} \frac{\partial \mathbf{Q} }{\partial t} =   \sum_{\xi = x,y,z}\mathbf{A}_{\xi}\frac{\partial{\mathbf{Q}}}{\partial \xi} = \div \mathbf{F} \left(\mathbf{Q} \right) + \sum_{\xi = x, y, z}\mathbf{B}_{\xi}\left(\grad\mathbf{Q}\right), \quad \mathbf{P} = \mathbf{P}^T > 0,
\end{align}
 for some flux function $\mathbf{F} \left(\mathbf{Q} \right) = [\mathbf{F}_{x}\left(\mathbf{Q} \right), \mathbf{F}_{y}\left(\mathbf{Q} \right), \mathbf{F}_{z} \left(\mathbf{Q} \right)]^T$ and coefficient matrices $\mathbf{B}_{\xi}$.

  \begin{Definition}
  Consider the equation of motion \eqref{eq:linear_wave}.  Suppose that the spatial operator can be split into:
  \begin{align}\label{eq:anti_split_operator}
   \sum_{\xi = x,y,z}\mathbf{A}_{\xi}\frac{\partial{\mathbf{Q}}}{\partial \xi} = \div \mathbf{F} \left(\mathbf{Q} \right) + \sum_{\xi = x, y, z}\mathbf{B}_{\xi}\left(\grad\mathbf{Q}\right),
   \end{align}
     where the first term, with $\div \mathbf{F} \left(\mathbf{Q} \right)$, is called the conservative flux term and the second term, with $\mathbf{B}_{\xi}\left(\grad\mathbf{Q}\right)$, is the non-conservative-products flux term.
   The split operators, given in \eqref{eq:anti_split_operator}, are anti-symmetric if 
\begin{align}\label{eq:anti_symmetry}
  \mathbf{Q}^T\mathbf{B}_\xi \left(\grad\mathbf{Q}\right) - \frac{\partial \mathbf{Q}^T}{\partial \xi} \mathbf{F}_\xi \left(\mathbf{Q} \right) = 0.
\end{align}
  \end{Definition}
The anti-symmetric split-form, \eqref{eq:gen_hyp_antisymmetry} with \eqref{eq:anti_symmetry}, can be useful when designing provably stable approximations on curvilinear meshes for complex geometries, and variable material properties.

There are several ways of casting a hyperbolic PDE, such as \eqref{eq:linear_wave},  in the anti-symmetric form, see the examples in \cite{Norstrom2006, Kopriva2006, ThomasLombad1979, KoprivaGassner2014, KoprivaGassner2015}. Here, our choice of the anti-symmetric form is motivated by the underlying physics.
For linear elasticity we will use
  {
\small
\begin{align}\label{eq:conervative_flux_and_nonconervative_product}
\mathbf{F}_{\xi} \left(\mathbf{Q} \right)  =   \begin{pmatrix}
{e}_{\xi x}\sigma_{xx} + {e}_{\xi y}\sigma_{xy} + {e}_{\xi z}\sigma_{xz}\\
{e}_{\xi x}\sigma_{xy} + {e}_{\xi y}\sigma_{yy} + {e}_{\xi z} \sigma_{yz}\\
{e}_{\xi x}\sigma_{xz} + {e}_{\xi y}\sigma_{yz} + {e}_{\xi z}\sigma_{zz}\\
0 \\
0\\
0 \\
0\\
0\\
  0\\
  \end{pmatrix},
  \quad
\mathbf{B}_{\xi} \left(\grad\mathbf{Q}\right)  =  \begin{pmatrix}
 0\\
 0\\
  0\\
  {e}_{\xi x}\frac{\partial v_x}{\partial {\xi}}  \\
{e}_{\xi y}\frac{\partial v_y}{\partial {\xi}}  \\
{e}_{\xi z}\frac{\partial v_z}{\partial {\xi}}  \\
{e}_{\xi y}\frac{\partial v_x}{\partial {\xi}}  + {e}_{\xi x}\frac{\partial v_y}{\partial {\xi}}  \\
{e}_{\xi z}\frac{\partial v_x}{\partial {\xi}}  + {e}_{\xi x}\frac{\partial v_z}{\partial {\xi}}  \\
{e}_{\xi z}\frac{\partial v_y}{\partial {\xi}}  + {e}_{\xi y}\frac{\partial v_z}{\partial {\xi}}   \\
  \end{pmatrix},
  \quad
  \mathbf{e}_{\xi} = \left({e}_{\xi x}, {e}_{\xi y}, {e}_{\xi z}\right)^T,
\end{align}
  }
  where the conservative flux comes from conservation of momentum, that is Newton's second law of motion. The non-conservative product term comes from Hooke's law, the constitutive relation relating stress and strain. In a Cartesian coordinate system, $\mathbf{e}_{\xi}$ are the canonical bases defined in \eqref{eq:canonical_basis}. In general curvilinear coordinates the bases $\mathbf{e}_{\xi}$ are arbitrary nonzero vectors, with $|\mathbf{e}_{\xi}| > 0$.
\begin{lemma}\label{lem:anti_symmetry}
Consider the split form of the equation \eqref{eq:gen_hyp_antisymmetry}, with the split operators defined in  \eqref{eq:conervative_flux_and_nonconervative_product}.
For  arbitrary  bases vectors $\mathbf{e}_{\xi}$, with $|\mathbf{e}_{\xi}| > 0$, the conservative flux term and the non-conservative-products flux term  satisfy the anti-symmetric property \eqref{eq:anti_symmetry}, that is
\[
\mathbf{Q}^T\mathbf{B}_\xi \left(\grad\mathbf{Q}\right) - \frac{\partial \mathbf{Q}^T}{\partial \xi} \mathbf{F}_\xi \left(\mathbf{Q} \right) = 0.
\]
\end{lemma}
\begin{proof}
The proof of the lemma follows from direct calculations, by using \eqref{eq:conervative_flux_and_nonconervative_product} and evaluating the products in \eqref{eq:anti_symmetry}.
\end{proof}

Another important consequence of the choice of the split operators  \eqref{eq:conervative_flux_and_nonconervative_product} is the fact
\begin{align}\label{state_times_flux}
 \mathbf{Q}^T\mathbf{F}_\xi \left(\mathbf{Q} \right) =  \mathbf{v}^T\left(\bar{\bar{\sigma}}\mathbf{e}_\xi\right) = \mathbf{v}^T\mathbf{T}^{(\xi)}. 
\end{align}
\begin{theorem}\label{Theo:energy_estimate_BC}
 Consider the  split form of   equation of motion \eqref{eq:gen_hyp_antisymmetry} for the linear elastic wave equation, defined by \eqref{eq:conervative_flux_and_nonconervative_product} in the spatial domain $(x,y,z)\in \Omega$ where $\Gamma$ denotes the boundary of the domain and $\mathbf{n}$  is the outward unit normal on the boundary.  Let $\mathbf{T}  = \bar{\bar{\sigma}}\mathbf{n}$ denote the traction on the boundary,  and the energy density $dE$ defined by \eqref{eq:elastic_energy_density} for an elastic medium. The solutions of  the anti-symmetric split form \eqref{eq:gen_hyp_antisymmetry} satisfy 
  \begin{align}\label{eq:energy_estimate_bc}
\frac{d}{dt}E\left(t\right) = \oint_{\Gamma} \mathbf{v}^T\mathbf{T} dS, \quad E\left(t\right) =   \int_{\Omega} dE > 0.
 \end{align}
 \end{theorem}
\begin{proof}
From the left, multiply the  split form of   equation of motion \eqref{eq:gen_hyp_antisymmetry} by the transpose of the solution $\mathbf{Q}^T$ and integrate over the whole domain $\Omega$, having
    \begin{align}\label{eq:energy_transformed_1_cart}
\int_{\Omega}\mathbf{Q}^T{\mathbf{P}}^{-1} \frac{\partial }{\partial t} \mathbf{Q} {dx dy dz} = \int_{\Omega}\mathbf{Q}^T \div \mathbf{F} \left(\mathbf{Q} \right) {dx dy dz} + \sum_{{\xi}= x, y, z}\int_{\Omega}\mathbf{Q}^T\mathbf{B}_{\xi}\left(\grad\mathbf{Q}\right) {dx dy dz} .
\end{align}
 On the left hand side  of \eqref{eq:energy_transformed_1} we recognize time derivative of the energy $ E\left(t\right)$.
  On the right hand side of \eqref{eq:energy_transformed_1_cart}, integrate--by--parts the first (conservative flux) term, we have
      \begin{align}\label{eq:energy_transformed_4}
\frac{d}{dt}E\left(t\right) = \sum_{{\xi}= x, y, z}\int_{\Omega}\left[ \mathbf{Q}^T\mathbf{B}_\xi \left(\grad\mathbf{Q}\right) - \frac{\partial \mathbf{Q}^T}{\partial \xi} \mathbf{F}_\xi \left(\mathbf{Q} \right)  \right]{dx dy dz} + \oint_{\Gamma} \mathbf{v}^T\mathbf{T} dS.
\end{align}
  By the anti-symmetric property \eqref{eq:anti_symmetry}, the volume terms vanish  having 
  \begin{align}\label{eq:energy_transformed_3_cart}
\frac{d}{dt}E\left(t\right)  = \oint_{\Gamma} \mathbf{v}^T\mathbf{T} dS.
 \end{align}
\end{proof}

In the coming section below, we will introduce curvilinear coordinates and transformations to model geometrically complex elastic solids.
\section{Curvilinear coordinates and structure preserving coordinate transformations}
To simplify the presentation, we consider two DG elements filled with heterogeneous elastic media separated by an  interface at $\left(\widetilde{x}\left(y, z\right), y, z\right)$, where $\widetilde{x}\left(y, z\right)$ is an arbitrarily smooth level surface describing the surface of the interface.  The interface can be a pre-existing fault in the medium, 
or the interface between two adjacent DG elements.  
\begin{figure}[h]
    \centering
    \includegraphics[width=0.6\linewidth]{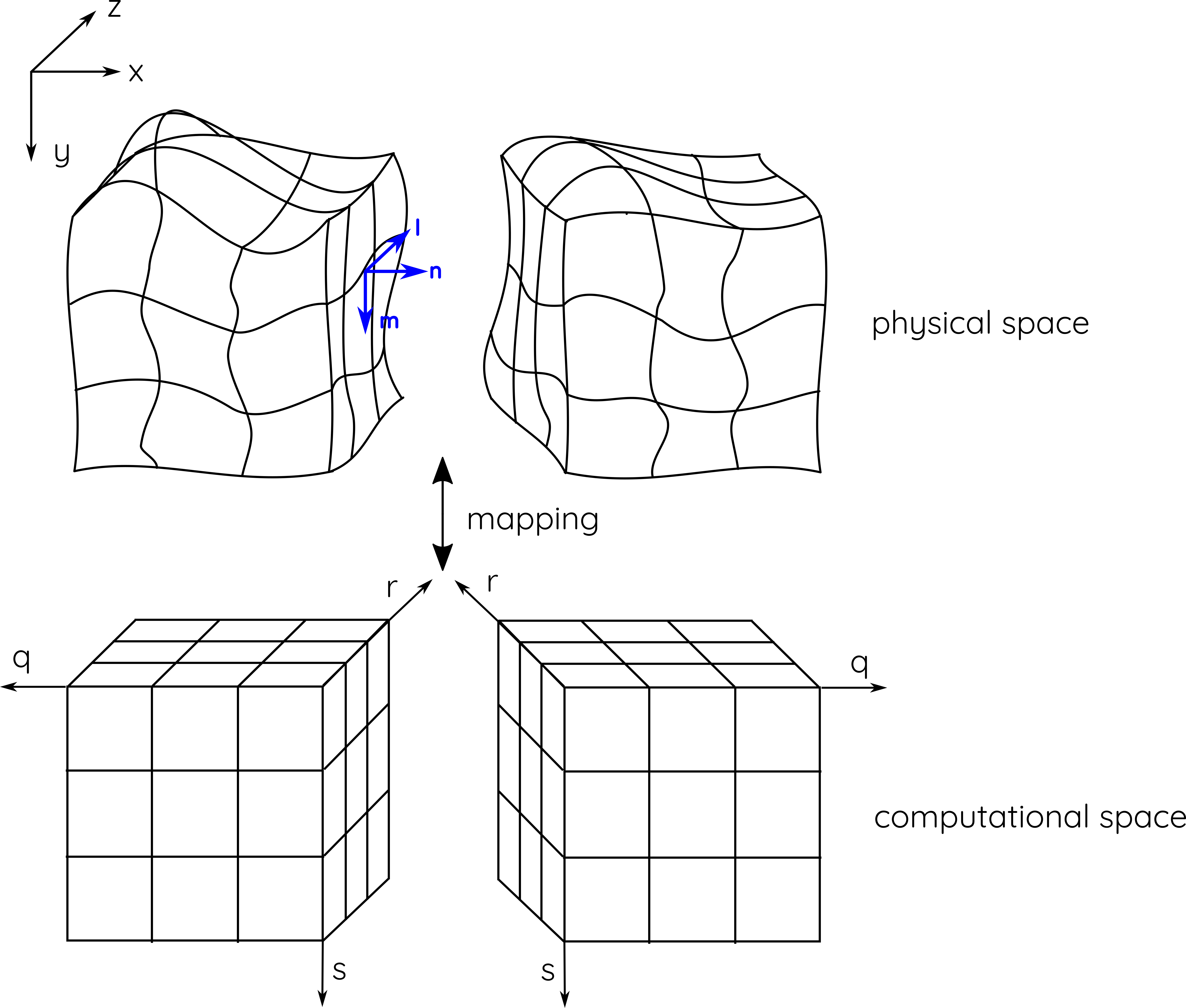}
    \caption{Internal and external boundary conforming curvilinear meshes and coordinate transformation}
    \label{fig:split_blocks}
\end{figure} 
For geometrically complex models, numerical treatments can be simplified if there is a smooth coordinate transformation from the physical space to a reference element, the unit cube $\left(q, r, s\right)\in \widetilde{\Omega} = [0,1]\times [0, 1]\times [0,1]$. We split the elements along the interface and map each element to the reference element, defined by
 \begin{equation}\label{eq:grid_transformation}
 \left(x\left(q, r, s\right), y\left(q, r, s\right), z\left(q, r, s\right)\right) \leftrightarrow \left(q\left(x, y, z\right), r\left(x, y, z\right), s\left(x, y, z\right)\right) .
 \end{equation}
Note that the mapping and coordinate transformations are element local. 


Let us denote the metric coefficients  and the Jacobian of the transformation \eqref{eq:grid_transformation} by 
 \[
 J = x_q\left(y_rz_s - z_ry_s\right) - y_q\left(x_rz_s - z_rx_s\right) + z_q\left(x_ry_s - y_rx_s\right) > 0,
 \]
 \begin{align*}
 q_x &= \frac{1}{J}\left(y_rz_s - z_ry_s\right), \quad r_x = \frac{1}{J}\left(z_qy_s - y_qz_s\right) , \quad s_x = \frac{1}{J}\left(y_qz_r - z_qy_r\right),\\
 q_y &= \frac{1}{J}\left(z_rx_s - x_rz_s\right), \quad r_y = \frac{1}{J}\left(x_qz_s - z_qx_s\right) , \quad s_y = \frac{1}{J}\left(z_qx_r - x_qz_r\right),\\
 q_z &= \frac{1}{J}\left(z_ry_s - y_rx_s\right), \quad r_z = \frac{1}{J}\left(y_qx_s - y_sx_q\right) , \quad s_z = \frac{1}{J}\left(x_qy_r - x_ry_q\right).
  \end{align*}
  Here, the subscripts denote partial derivatives, that is $x_q = \frac{\partial x}{\partial q}$,  $x_r = \frac{\partial x}{\partial r}$, etc. 
  
 To construct structure preserving coordinate transformation, we will use two different transformations of the spatial derivatives in the transformed coordinates \cite{DuruandDunham2016},  the
 conservative form:
 \begin{align}\label{eq:conservative}
\frac{\partial v }{\partial x} = \frac{1}{ J}\left(\frac{\partial }{\partial q}\left(Jq_x v\right) + \frac{\partial }{\partial r}\left(Jr_x v\right) + \frac{\partial }{\partial s}\left(Js_x v\right)\right),
 \end{align}
  and the non-conservative form:
 \begin{align}\label{eq:non-conservative}
 \frac{\partial v}{\partial x} = q_x\frac{\partial v}{\partial q} + r_x\frac{\partial v}{\partial r}  + s_x\frac{\partial v}{\partial s} .
 \end{align}
Note that in the continuous setting, the transformed derivatives \eqref{eq:conservative} and \eqref{eq:non-conservative} are mathematically equivalent.
 In nontrivial geometries, discrete approximations of the transformed derivatives \eqref{eq:conservative} and \eqref{eq:non-conservative}  will yield two different  discrete spatial operators.

In order to preserve the anti-symmetric structure, in  \eqref{eq:gen_hyp_antisymmetry}, it becomes natural to transform the derivatives in the conservative flux term using the conservative transformation \eqref{eq:conservative}, and the derivatives in the non-conservative-products term using the non-conservative transformation \eqref{eq:non-conservative}. 
 The elastic wave equation \eqref{eq:gen_hyp_antisymmetry} in the transformed curvilinear  coordinates $(q, r, s)$ is
  \begin{align}\label{eq:gen_hyp_transformed}
\widetilde{\mathbf{P}}^{-1} \frac{\partial }{\partial t} \mathbf{Q} = \div \mathbf{F} \left(\mathbf{Q} \right) + \sum_{\xi = q, r, s}\mathbf{B}_{\xi}\left(\grad\mathbf{Q}\right).
\end{align}
Here, $\grad = \left(\partial /\partial q, \partial /\partial r, \partial /\partial s\right)^T$ is the gradient operator,
$
\widetilde{\mathbf{P}} = J^{-1}\mathbf{P} $,
where $\mathbf{P}$ is the material matrix defined in \eqref{eq:elastic_coeff}, and
{
\small
\begin{align}\label{eq:flux_ncp}
 \mathbf{F}_{\xi} \left(\mathbf{Q} \right)  =   \begin{pmatrix}
J\left({\xi}_x\sigma_{xx} + {\xi}_y\sigma_{xy} + {\xi}_z\sigma_{xz}\right)\\
J\left({\xi}_x\sigma_{xy} + {\xi}_y\sigma_{yy} + {\xi}_z\sigma_{yz}\right)\\
J\left({\xi}_x\sigma_{xz} + {\xi}_y\sigma_{yz} + {\xi}_z\sigma_{zz}\right)\\
0 \\
0\\
0 \\
0\\
0\\
  0\\
  \end{pmatrix},
  \quad
 \mathbf{B}_{\xi} \left(\grad\mathbf{Q}\right)  =  \begin{pmatrix}
 0\\
 0\\
  0\\
J{\xi}_x\frac{\partial v_x}{\partial {\xi}}  \\
J{\xi}_y\frac{\partial v_y}{\partial {\xi}}  \\
J{\xi}_z\frac{\partial v_z}{\partial {\xi}} \\
J\left({\xi}_y\frac{\partial v_x}{\partial {\xi}}  + {\xi}_x\frac{\partial v_y}{\partial {\xi}} \right)  \\
J\left({\xi}_z\frac{\partial v_x}{\partial {\xi}}  + {\xi}_x\frac{\partial v_z}{\partial {\xi}} \right) \\
J\left({\xi}_z\frac{\partial v_y}{\partial {\xi}}  + {\xi}_y\frac{\partial v_z}{\partial {\xi}} \right)  \\
  \end{pmatrix},
  \quad
  \end{align}
  }
  where $\xi = q, r, s$, with $\xi_{\eta} = \partial \xi /\partial \eta$ and $\eta  = x, y, z$.
Note that in \eqref{eq:gen_hyp_transformed} all spatial derivatives for the stress fields are transformed using the conservative form \eqref{eq:conservative} and all spatial derivatives for the velocity fields are transformed using the non-conservative form  \eqref{eq:non-conservative}. 
When discrete approximations are introduced, this is crucial in order minimize the number of floating point operations, and also prove numerical stability \cite{DuruandDunham2016}. 

We will show that the transformed equation of motion, \eqref{eq:gen_hyp_transformed} with \eqref{eq:flux_ncp}, preserves the anti-symmetric property \eqref{eq:anti_symmetry}.
\begin{lemma}\label{Lem:Anti-symmetry}
Consider the transformed equation of motion \eqref{eq:gen_hyp_transformed}, in curvilinear coordinates, with the conservative flux terms and non-conservative-products flux  terms given by  \eqref{eq:flux_ncp}.
The corresponding spatial operators satisfy the anti-symmetric property \eqref{eq:anti_symmetry}, that is
\begin{align*}
  \mathbf{Q}^T\mathbf{B}_\xi \left(\grad\mathbf{Q}\right) - \frac{\partial \mathbf{Q}^T}{\partial \xi} \mathbf{F}_\xi \left(\mathbf{Q} \right) = 0.
\end{align*}
\end{lemma}
\begin{proof}
With $ \mathbf{e}_{\xi} = J\left({\xi}_x, {\xi}_y, {\xi}_z\right)^T$, the proof of Lemma \ref{Lem:Anti-symmetry} follows directly from Lemma \ref{lem:anti_symmetry}.
\end{proof}

Let $\xi = q, r, s$, and $\Gamma$ denote a boundary face at $\xi = 0$ or $\xi = 1$.
The the positively pointing unit normals on the boundary  are given by
\begin{align}\label{eq:unit_normal}
\mathbf{n} 
=
 \frac{1}{\sqrt{\xi_x^2 + \xi_y^2 + \xi_z^2}}
  \begin{pmatrix}
 \xi_x\\
 \xi_y\\
\xi_z
 \end{pmatrix}.
\end{align}
Note  again that 
\begin{align}\label{state_times_flux_00}
 \mathbf{Q}^T\mathbf{F}_\xi \left(\mathbf{Q} \right) = J\sqrt{\xi_x^2 + \xi_y^2 + \xi_z^2}\mathbf{v}^T\mathbf{T},
\end{align}
where $\mathbf{v} = [v_x, v_y, v_z]^T$ is the velocity vector, and $\mathbf{T} = [T_x, T_y, T_z]^T = \bar{\bar{\sigma}}\mathbf{n}$ is the traction vector.

Let us introduce  the reference boundary surface $\widetilde{\Gamma} = [0, 1]\times[0, 1]$, and define the boundary term
  \begin{align}\label{eq:boundaryterm_101}
& \mathrm{BTs}\left(v , T \right) := \oint_{\Gamma} \mathbf{v}^T\mathbf{T} dS \nonumber \\
 &=\sum_{\xi = q, r, s}\left(\int_{\widetilde{\Gamma}} \left(J\sqrt{\xi_x^2 + \xi_y^2 + \xi_z^2}\right)\mathbf{v}^T\mathbf{T} \Big|_{\xi = 1}\frac{dqdrds}{d\xi} - \int_{\widetilde{\Gamma}} \left(J\sqrt{\xi_x^2 + \xi_y^2 + \xi_z^2}\right)\mathbf{v}^T\mathbf{T} \Big|_{\xi = 0} \frac{dqdrds}{d\xi} \right).
  \end{align}
 Another important consequence of the the transformed anti-symmetric split form \eqref{eq:gen_hyp_transformed} with \eqref{eq:flux_ncp} is the following lemma
 \begin{lemma}\label{Lem:Anti-symmetry_00}
Consider the conservative flux term defined in \eqref{eq:flux_ncp}. We have
\begin{align}\label{eq:flux_integration_by_parts}
 \int_{\widetilde{\Omega}}\mathbf{Q}^T \div \mathbf{F} \left(\mathbf{Q} \right) dqdrds =  &-\sum_{{\xi}= q, r, s}\int_{\widetilde{\Omega}}\frac{\partial \mathbf{Q}^T}{\partial \xi} \mathbf{F}_\xi \left(\mathbf{Q} \right) dqdrds +  \mathrm{BTs}\left(v , T \right),
\end{align}
where the boundary term $\mathrm{BTs}\left(v , T \right)$  is defined in \eqref{eq:boundaryterm_101}
\end{lemma}
\begin{proof}
Consider 
  \begin{align}\label{eq:flux_integration_by_parts_0}
 \int_{\widetilde{\Omega}}\mathbf{Q}^T \div \mathbf{F} \left(\mathbf{Q} \right) dqdrds &=  \sum_{\xi = q, r, s}\int_{\widetilde{\Omega}} \mathbf{Q}^T  \frac{\partial \mathbf{F}_\xi \left(\mathbf{Q}  \right)}{\partial \xi}dqdrds, 
\end{align}
 and integrate-by-parts, we have
   \begin{equation}\label{eq:flux_integration_by_parts_1}
   \begin{split}
 \int_{\widetilde{\Omega}}\mathbf{Q}^T \div \mathbf{F} \left(\mathbf{Q} \right) dqdrds =  &\sum_{\xi = q, r, s}\left(-\int_{\widetilde{\Omega}}\frac{\partial \mathbf{Q}^T}{\partial \xi} \mathbf{F}_\xi \left(\mathbf{Q} \right) dqdrds 
 +  \int_{\widetilde{\Gamma}}\mathbf{Q}^T\mathbf{F}_\xi \left(\mathbf{Q} \right)\Big|_{\xi = 0}^{\xi = 1}\frac{dqdrds}{d\xi} \right)
 \end{split}
\end{equation}
Using the fact \eqref{state_times_flux_00}, and the boundary term $\mathrm{BTs}\left(v , T \right)$  defined in \eqref{eq:boundaryterm_101} completes the proof.
\end{proof}

Introduce the energy density in the transformed space
\begin{align}\label{eq:energy_density_transformed_space}
\frac{d\widetilde{E}}{dqdrds} = \frac{1}{2}{\left[\mathbf{Q}^T \widetilde{\mathbf{P}}^{-1} \mathbf{Q}\right]} > 0.
\end{align}
Analogous to Theorem \ref{Theo:energy_estimate_BC}, we also have
\begin{theorem}\label{Theo:energy_estimate_curvilinear}
 Consider the transformed equation of motion \eqref{eq:gen_hyp_transformed}, in curvilinear coordinates,  $(q,r,s)\in \widetilde{\Omega}$, with the flux terms and non-conservative products  terms given by  \eqref{eq:flux_ncp}. Let $\xi = q, r, s$, and $\Gamma$ denote a boundary face at $\xi = 0$ or $\xi = 1$, where $\mathbf{T}  = \bar{\bar{\sigma}}\mathbf{n}$ is the traction vector on the boundary. The solutions of the transformed equation \eqref{eq:gen_hyp_transformed} satisfy 
  \begin{align}\label{eq:energy_estimate_bc_0}
\frac{d}{dt}E\left(t\right) = \mathrm{BTs}\left(v , T \right), \quad E\left(t\right) =   \int_{\widetilde{\Omega}} d\widetilde{E} > 0,
 \end{align}
 where $\mathrm{BTs}\left(v , T \right)$ is the boundary term defined in \eqref{eq:boundaryterm_101}.
 \end{theorem}
\begin{proof}
From the left, multiply the  transformed  equation of motion \eqref{eq:gen_hyp_transformed} by $\mathbf{Q}^T$ and integrate over the whole domain $\widetilde{\Omega}$, having
    \begin{align}\label{eq:energy_transformed_1}
\int_{\widetilde{\Omega}}\mathbf{Q}^T\widetilde{\mathbf{P}}^{-1} \frac{\partial }{\partial t} \mathbf{Q} dqdrqs = \int_{\widetilde{\Omega}}\mathbf{Q}^T \div \mathbf{F} \left(\mathbf{Q} \right) dqdrqs + \sum_{{\xi}= q, r, s}\int_{\widetilde{\Omega}}\mathbf{Q}^T\mathbf{B}_{\xi}\left(\grad\mathbf{Q}\right) dqdrqs .
\end{align}
 On the left hand side  of \eqref{eq:energy_transformed_1} we recognize time derivative of the energy.
  On the right hand side of \eqref{eq:energy_transformed_1}, we use Lemma \ref{Lem:Anti-symmetry_00}, and replace the conservative flux term with \eqref{eq:flux_integration_by_parts}, we have
      \begin{align}\label{eq:energy_transformed_4}
\frac{d}{dt}E\left(t\right) = \sum_{{\xi}= q, r, s}\int_{\widetilde{\Omega}}\left[\mathbf{Q}^T\mathbf{B}_\xi \left(\grad\mathbf{Q}\right) - \frac{\partial \mathbf{Q}^T}{\partial \xi} \mathbf{F}_\xi \left(\mathbf{Q} \right)\right]dqdrqs + \mathrm{BTs} ,
\end{align}
 where $\mathrm{BTs}$ is the boundary term  defined in \eqref{eq:boundaryterm_101}.
 By Lemma \ref{Lem:Anti-symmetry}, the volume terms in \eqref{eq:energy_transformed_4} vanish,  having 
  \begin{align}\label{eq:energy_transformed_3}
\frac{d}{dt}E\left(t\right) = \mathrm{BTs}.
 \end{align}
\end{proof}

The boundary term, $\mathrm{BTs}$, is the rate of the  work done by the traction, $\mathbf{T}  = \bar{\bar{\sigma}}\mathbf{n}$, against the boundary. Note that the energy rate is controlled by the boundary term, $\mathrm{BTs}$. In a bounded domain, well-posed boundary conditions are designed such that the boundary term is never positive, $\mathrm{BTs} \le 0$.
\section{Boundary and interface conditions}
 We will now introduce physical boundary and interface conditions prescribed at element faces. 
  In Figure \ref{fig:split_blocks}, for each of the two blocks, there are five external boundaries, and one  internal boundary connecting the two blocks. 
 First, we will consider boundary conditions posed at external element boundaries, and proceed later to interface conditions acting at the internal element face, connecting the two locally adjacent elements.
  \subsection{Boundary conditions}
 Here, we present linear well-posed boundary conditions closing the external boundaries.
   Consider the 5 external boundaries of each elastic block, see Figure \ref{fig:split_blocks}. To define well-posed boundary conditions,  we introduce the basis vectors   denoted  by $\mathbf{n} = \left(n_x, n_y, n_z\right), \mathbf{m} =  \left(m_x, m_y, m_z\right), \textbf{{l}} = \left(l_x, l_y, l_z\right)$, where, $\mathbf{n}$ is the unit normal, defined by \eqref{eq:unit_normal}, on the boundary pointing in the positive $\xi$-direction, with $\xi = q, r, s$. 
\begin{figure}[h]
    \centering
    \includegraphics[width=0.45\linewidth]{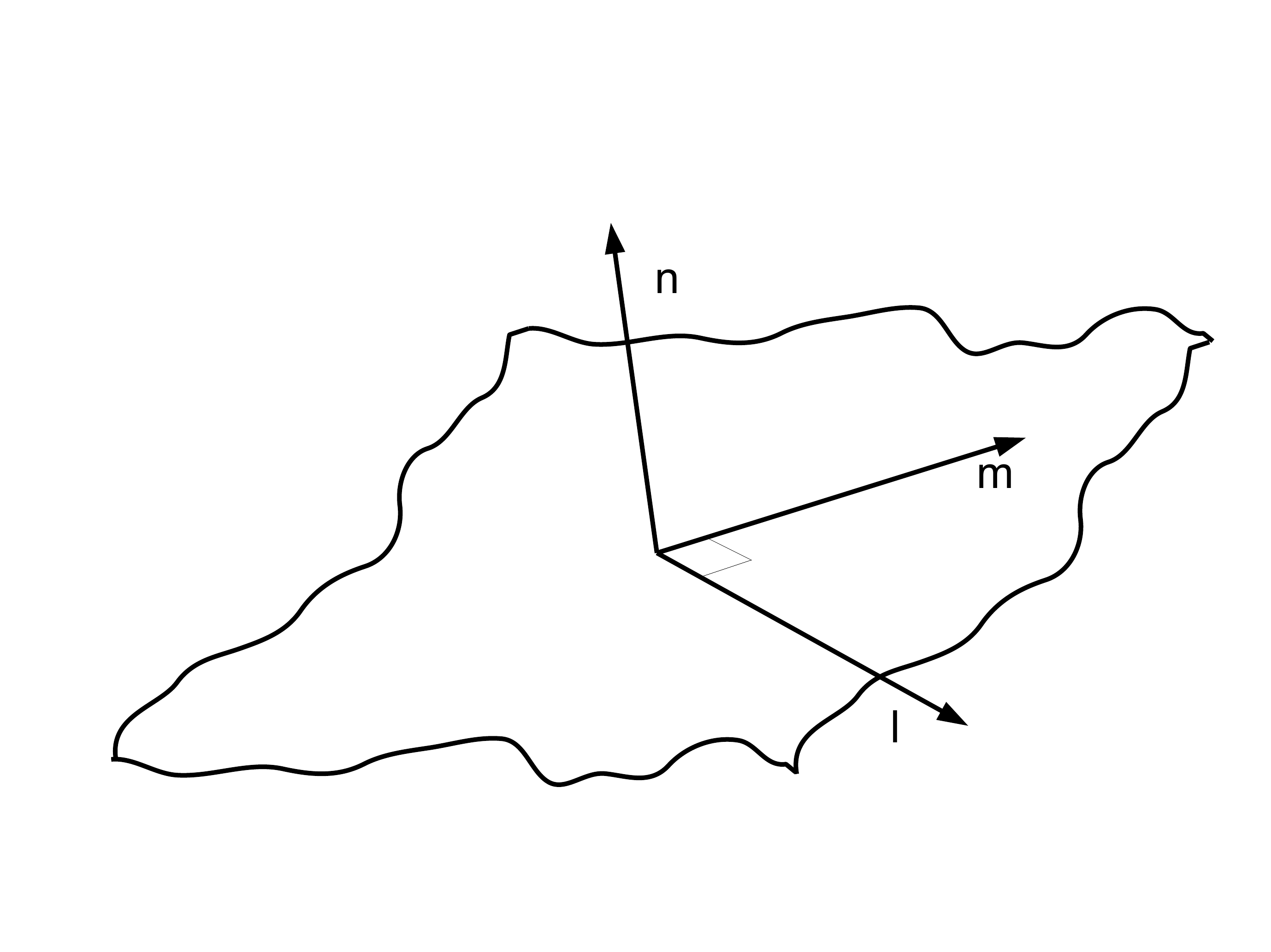}
    \caption{An element face and local basis vectors.}
    \label{fig:Local_basis_vectors}
\end{figure} 
Note that the unit vectors are locally defined on each point on the boundary, see also Figure \ref{fig:Local_basis_vectors}.
Given the unit normal $\mathbf{n} $,  defined by \eqref{eq:unit_normal}, we can construct the other two  basis vectors $\mathbf{l}, \mathbf{m}$, as follows:
\begin{align}\label{eq:basis_vectors}
\mathbf{m} = \frac{\mathbf{m}_0 - \left(\mathbf{n}^T\mathbf{m_0}\right) \mathbf{n}}{|\mathbf{m}_0 - \left(\mathbf{n}^T\mathbf{m_0}\right) \mathbf{n}|}, \quad \mathbf{l} = \mathbf{n}\times \mathbf{m}.
\end{align}
Here, $\mathbf{m}_0  \left(\ne \pm \mathbf{n}\right)$ is an arbitrary unit vector.
Let the  local impedances  at the boundary be denoted by $Z_\eta, \eta = l, m,n,$
where $Z_n = \rho c_n$ is the p--wave impedance and $Z_m =\rho c_m$, $Z_l= \rho c_l$ are   the s--wave impedances. Here, $c_n, c_m, c_l$ are the corresponding effective wavespeeds defined by
\begin{align}\label{eq:effectivewavespeed}
 c_{n} =  \sqrt{\sum_{\xi=x,y,z}\left(n_{\xi}c_{p\xi}\right)^2}, \quad
  c_{m} =  \sqrt{\sum_{\xi=x,y,z}\left(n_{\xi}c_{sh\xi}\right)^2}, \quad
   c_{l} =  \sqrt{\sum_{\xi=x,y,z}\left(n_{\xi}c_{sv\xi}\right)^2}.
\end{align}
In anisotropic media with a geometrically complex curvilinear coordinate system,  the effective wavespeeds  depend on the orientation of the normal vector $\mathbf{n}$.
Note that ${n_x^2 + n_y^2  + n_z^2} = 1$.
Thus, in an isotropic medium  the effective wavespeeds  are given by
$
c_n = c_p, \quad c_m =c_l = c_s.
$

On each point on the boundary, we denote the particle velocity vector, traction vector, and the local rotation matrix  on the boundary  by 
\begin{align}\label{eq:characteristics}
\mathbf{v} = \begin{pmatrix}
v_{x}\\
v_{y}\\
v_{z}
\end{pmatrix},
 \quad
\mathbf{T} = \begin{pmatrix}
T_{x}\\
T_{y}\\
T_{z}
\end{pmatrix} = \bar{\bar{\sigma}}\mathbf{n}, \quad \mathbf{R} = \begin{pmatrix}
n_{x}&n_{y}&n_{z}\\
m_{x}&m_{y}&m_{z}\\
l_{x}&l_{y}&l_{z}
\end{pmatrix},
\end{align}
where $\det{(\mathbf{R})} \ne 0$ and $ \mathbf{R}^{-1} = \mathbf{R}^T$. 

Next, rotate the particle velocity and traction vectors into the local orthonormal basis, $\mathbf{l}$ ,  $\mathbf{m}$  and  $\mathbf{n}$, having 
\begin{align}\label{eq:local_velocity_tractions}
v_\eta  = \left(\mathbf{R}\mathbf{v}\right)_\eta , \quad T_\eta  = \left(\mathbf{R}\mathbf{T}\right)_\eta , \quad \eta = l, m, n.
\end{align}
 The corresponding in and out of the domain characteristics at the boundary  are
 \begin{align}\label{eq:characteristics}
{q}_\eta  = \frac{1}{2}\left({Z}_\eta {v}_\eta  + {T}_\eta \right), \quad {p}_\eta  = \frac{1}{2}\left({Z}_\eta {v}_\eta  - {T}_\eta \right), \quad Z_\eta  > 0.
\end{align}
Here, the characteristics defined in \eqref{eq:characteristics} are plane p--waves and plane s--waves propagating along the normal vector $\mathbf{n}$ on the boundary. 

At the boundary $\xi = 1$ $(\xi = 0)$, if $Z_\eta  > 0$ then ${q}_\eta $ (${p}_\eta $)  are the characteristics  going into the domain  and  ${p_\eta }$ (${q_\eta }$) the  characteristics going out of the domain.
The number of boundary conditions must correspond to the number of characteristics going into the domain, see \cite{DuruandDunham2016, GustafssonKreissOliger1995}.
 We consider linear boundary conditions, 
{
\begin{equation}\label{eq:BC_General2}
\begin{split}
&{q}_\eta  - {\gamma_\eta }{p}_\eta  = 0 \iff \frac{Z_{\eta}}{2}\left({1-\gamma_\eta }\right){v}_\eta  -\frac{1+\gamma_\eta }{2} {T}_\eta  = 0,  \quad \xi = 0, \\
&{p}_\eta  - {\gamma_\eta }{q}_\eta  = 0 \iff \frac{Z_{\eta}}{2} \left({1-\gamma_\eta }\right){v}_\eta  + \frac{1+\gamma_\eta }{2}{T}_\eta  = 0,  \quad  \xi = 1,
 \end{split}
\end{equation}
}
where the reflection coefficients $\gamma_\eta $ are real numbers with $ 0 \le |\gamma_\eta |\le 1$.
The boundary conditions \eqref{eq:BC_General2} specify the ingoing characteristics on the boundary in terms of the outgoing characteristics. 
In an elastic medium, we have $Z_\eta  > 0$ for all $\eta = l, m, n$, and there are three boundary conditions at each boundaries $\xi = 1$, ($\xi = 0$).
 The boundary condition \eqref{eq:BC_General2}, can describe several physical situations.
 We have a free-surface boundary condition  if $\gamma_\eta  = 1$, an absorbing boundary condition  if $\gamma_\eta  = 0$ and a clamped boundary condition  if $\gamma_\eta  = -1$.
For later use in deriving energy estimates, we note that
  \begin{align}\label{eq:simplify_3}
& \text{at} \quad \xi = 0, \quad v_\eta T_\eta  =  \frac{Z_\eta \left(1-\gamma_\eta \right)}{\left(1+\gamma_\eta \right)}v_{\eta}^2 = \frac{\left(1+\gamma_\eta \right)}{Z_\eta \left(1-\gamma_\eta \right)}T_{\eta}^2  >0, \quad \forall |\gamma_\eta | < 1, \quad \text{and} \quad  v_\eta T_\eta  = 0, \quad \forall |\gamma_\eta | = 1, \nonumber
 \\
 & \text{at} \quad \xi = 1, \quad  v_\eta T_\eta  = - \frac{Z_\eta \left(1-\gamma_\eta \right)}{\left(1+\gamma_\eta \right)}v_{\eta}^2 = - \frac{\left(1+\gamma_\eta \right)}{Z_\eta \left(1-\gamma_\eta \right)}T_{\eta}^2  < 0, \quad \forall |\gamma_\eta | < 1, \quad \text{and} \quad  v_\eta T_\eta  = 0, \quad \forall |\gamma_\eta | = 1.
  \end{align}
 
\begin{lemma}\label{Lem:BTs}
Consider the well-posed boundary conditions \eqref{eq:BC_General2} with $|\gamma_\eta | \le 1$. The boundary term $\mathrm{BTs}$ defined in \eqref{eq:boundaryterm_101}   is negative semi-definite, $\mathrm{BTs} \le 0$, for all $Z_\eta > 0$.
\end{lemma}
\begin{proof}
With $\mathbf{v}^T\mathbf{T}  = \left(\mathbf{R}\mathbf{v}\right)^T\left(\mathbf{R}\mathbf{T}\right) = \sum_{\eta = l,m,n}{v_\eta T_\eta }$, the boundary term $ \mathrm{BTs}\left(v , T \right)$ defined in \eqref{eq:boundaryterm_101} can be written as
  \begin{align}\label{eq:boundaryterm_2}
 \mathrm{BTs}\left(v , T \right)  
  =&\sum_{\xi = q, r, s}\int_{\widetilde{\Gamma}} \left(\left(J\sqrt{\xi_x^2 + \xi_y^2 + \xi_z^2}\right) \sum_{\eta = l,m,n} v_\eta T_\eta \right) \Big|_{\xi = 1}\frac{dqdrds}{d\xi} \nonumber \\
  -& \sum_{\xi = q, r, s} \int_{\widetilde{\Gamma}} \left(\left(J\sqrt{\xi_x^2 + \xi_y^2 + \xi_z^2}\right)\sum_{\eta = l,m,n} v_\eta T_\eta \right)\Big|_{\xi = 0} \frac{dqdrds}{d\xi} .
 %
  \end{align}
 Finally, the identity  \eqref{eq:simplify_3}  completes the proof of the lemma. 
\end{proof}

 Using the energy method we can now prove:
 \begin{theorem}\label{Theo:energy_estimate_BC_00}
 Consider the  transformed  equation of motion \eqref{eq:gen_hyp_transformed} subject to the boundary condition \eqref{eq:BC_General2}, with $|\gamma| \le 1$.  The solutions of the transformed  equation \eqref{eq:gen_hyp_transformed} subject to the boundary condition \eqref{eq:BC_General2}  satisfy 
  \begin{align}\label{eq:energy_estimate_bc}
\frac{d}{dt}E\left(t\right) = \mathrm{BTs} \le 0.
 \end{align}
 \end{theorem}
\begin{proof}
The proof of Theorem \ref{Theo:energy_estimate_BC_00} follows the same steps as in the proof of Theorem \ref{Theo:energy_estimate_curvilinear}, arriving at \eqref{eq:energy_transformed_3}.
We complete the proof using Lemma \ref{Lem:BTs}, and by ensuring that the boundary terms are never positive, $ \mathrm{BTs} \le 0$.
\end{proof}

\subsection{Interface conditions}
We introduce physical interface conditions acting at internal DG elements boundaries,
in elastic solids. These physical interface conditions will connect two adjacent elements elastic media.
One objective of this study is to use the physical conditions to patch DG elements together \cite{DuruGabrielIgel2017}. 
Consider the interface, as in Figure \ref{fig:split_blocks} and  denote the corresponding fields and material parameters in the positive/negative sides of the interface with the superscripts $+/-$. To define the interface conditions we rotate the particle velocity vector and the traction vector on the boundary  into the local orthogonal coordinates $\mathbf{n}, \mathbf{m}, \mathbf{l}$,  as in \eqref{eq:local_velocity_tractions}, having
\begin{align}\label{eq:local_positve_negative} 
v_\eta ^{\pm} = \left(\mathbf{R}\mathbf{v}^{\pm}\right)_\eta , \quad T_\eta ^{\pm} = \left(\mathbf{R}\mathbf{T}^{\pm}\right)_\eta , \quad \eta = l, m, n.
\end{align}
We define the jumping condition in the velocity fields as
\[
\lJump{{v}_\eta \rJump} =  v_{\eta}^{+} - v_{\eta}^{-}, \quad \eta = l, m, n.
\]
As before, to ensure well--posedness, the number of interface conditions must be equal to the number of ingoing characteristics at the interface. 

\begin{figure}[h!]
    \centering
    \includegraphics[width=0.5\textwidth]{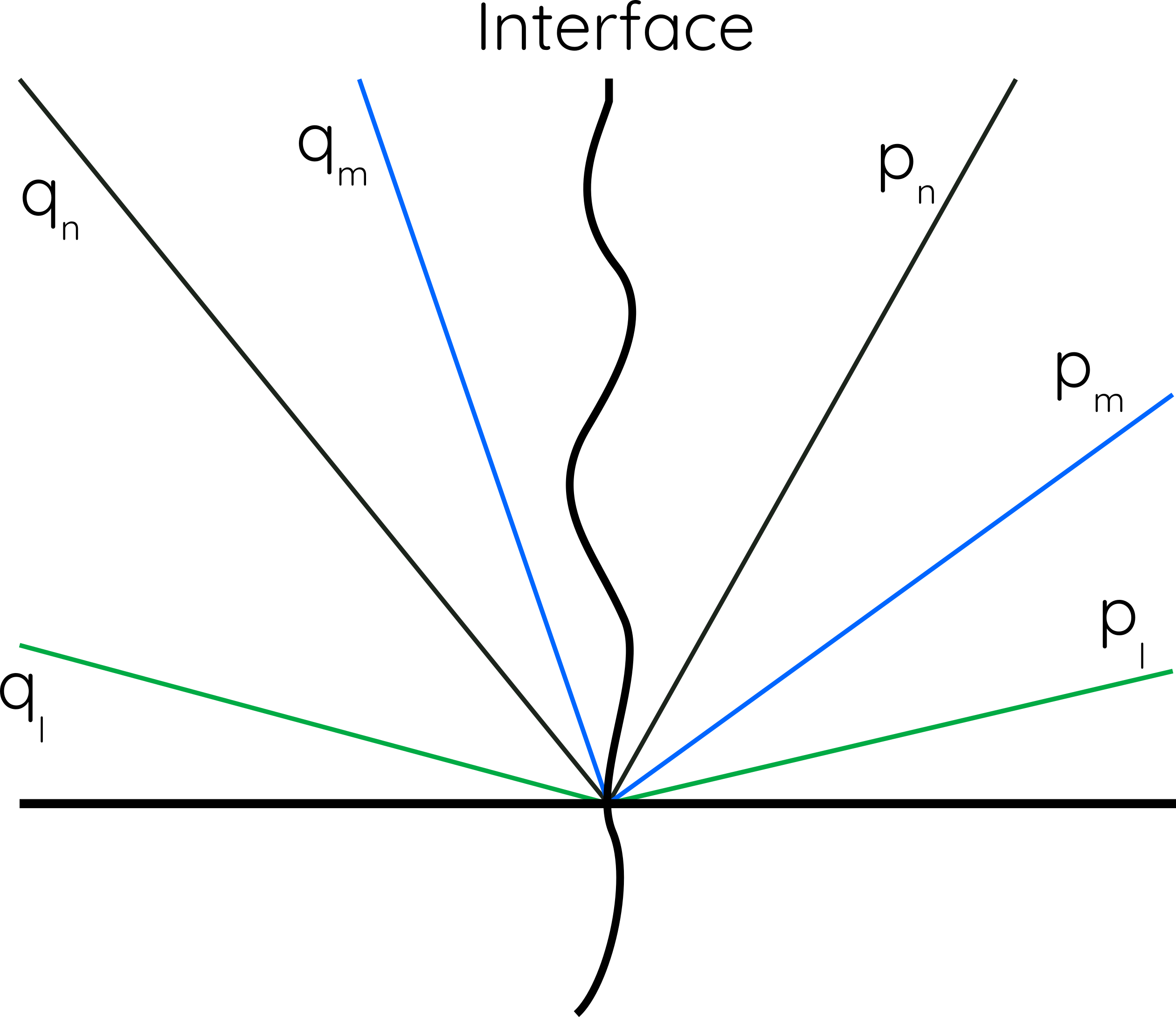}
    \caption{Sketch of characteristics propagating across an elastic-elastic interface.}
    \label{fig:characteristics_elastic_elastic}
\end{figure}

We consider two elastic solids in a locked  contact, with $Z_\eta ^{\pm} > 0$ for all $\eta = l, m, n$. The interface is locked,  that is, there is no opening/gap, no inter-penetration  and no slip.
As shown in figure \ref{fig:characteristics_elastic_elastic}, there are a total of  $6$ (3 characteristics going into the negative  element and 3 characteristics going into the positive  element)  ingoing characteristics at the interface. Therefore, we will need exactly 6 conditions specifying the relationships of the fields across the interface.
The interface conditions are force balance, and vanishing opening and slip velocities

\begin{equation}\label{eq:elastic_elastic_interface}
 T_{\eta}^{+}  = T_{\eta}^{-} = T_{\eta}, \quad \lJump{{v}_\eta \rJump} = 0, \quad  \eta = l, m, n.
\end{equation}

Now, we will show that the interface condition \eqref{eq:elastic_elastic_interface} with the equation of motion \eqref{eq:gen_hyp_transformed},  conserve the total mechanical energy. 
Introduce the interface term
\begin{align}\label{eq:interface_term}
  IT_s(v^{\pm}, T^{\pm}) = - \sum_{\xi = q, r, s}\int_{\widetilde{\Gamma}}\left(J\sqrt{\xi_x^2 + \xi_y^2 + \xi_z^2} \sum_{\eta = l,m,n}{T}_\eta  \lJump{{v}_\eta \rJump} \right) \frac{dqdrds}{d\xi} \equiv 0.
\end{align}
Thus, by  \eqref{eq:elastic_elastic_interface}  we have ${T}_\eta  \lJump{{v}_\eta \rJump}= 0$, thus, the interface term vanishes identically, $IT_s(v ^{\pm}, T ^{\pm}, Z ^{\pm}) \equiv 0$.
In particular we can prove the theorem:
  \begin{theorem}\label{Theo:linear_friction_energy_estimate}
Consider the modelling domain decomposed into two elements, $\Omega = \Omega_{-} \cup \Omega_{+}$, as in  Figure \ref{fig:split_blocks},  and denote the corresponding fields,  material parameters and energies in the positive/negative sides of the interface with the superscripts $+/-$.  Let each element be mapped to the reference element $\widetilde{\Omega} = [0,1]^3$, with the transformed  equation of motion \eqref{eq:gen_hyp_transformed}  defined on each element and subject to the interface conditions \eqref{eq:elastic_elastic_interface}. 
The sum of the energies satisfies
 \begin{align}\label{eq:energy_estimate_fault_30}
\frac{d}{dt}\left(E^{-}(t)+ E^{+}(t)\right)   = IT_s(v^{\pm}, T^{\pm}) \equiv 0.
 \end{align}
\end{theorem}
\begin{proof}
 Using the energy method, from  \eqref{eq:energy_transformed_1}--\eqref{eq:energy_transformed_3}, we again have
  \begin{align}\label{eq:energy_transformed_interface}
\frac{d}{dt}E\left(t\right) = \mathrm{BTs}\left(v , T \right),
 \end{align}
 where $\mathrm{BTs}\left(v , T \right)$ is the boundary term defined in \eqref{eq:boundaryterm_101}.
 Collecting contributions from both elements and ignoring all other boundaries, excepting the boundary at the interface, we have
 \begin{align}\label{eq:energy_estimate_fault_30_interface}
\frac{d}{dt}\left(E^{-}(t)+ E^{+}(t)\right)   = IT_s(v ^{\pm}, T ^{\pm}, Z ^{\pm}).
 \end{align}
 Using the fact  \eqref{eq:interface_term} that the interface term vanishes identically, $IT_s(v ^{\pm}, T ^{\pm}, Z ^{\pm}) \equiv 0$, completes the proof.
\end{proof}
 
 Theorems \eqref{Theo:energy_estimate_BC_00} and \eqref{Theo:linear_friction_energy_estimate} prove that the corresponding IBVPs are well-posed and asymptotically stable.
 The challenge, for the DG scheme, is how to incorporate the boundary conditions and the interface conditions in a seamless and provably stable manner,  for the discrete problem, in general heterogeneous media with complex geometries. To succeed, we will introduce hat-variables, so that we can simultaneously construct data for the velocity fields and traction fields, at internal and external element faces.
%
\section{Hat-variables and physics based fluxes}\label{sec:hat_variables}
We will now reformulate the boundary condition  \eqref{eq:BC_General2} and  interface condition \eqref{eq:elastic_elastic_interface}  by introducing transformed (hat-) variables so that we can simultaneously construct (numerical) boundary/interface data for particle velocities and tractions.  The hat-variables encode the solution of the IBVP on the boundary/interface. The hat-variables  will be constructed such that they preserve the amplitude of the outgoing characteristics and exactly satisfy the physical boundary conditions \cite{DuruandDunham2016}. To be more specific, the hat-variables are solutions of the Riemann problem constrained against physical boundary conditions \eqref{eq:BC_General2},  and  the interface condition \eqref{eq:elastic_elastic_interface}. We refer the reader to \cite{DuruGabrielIgel2017} for more detailed discussion. Once the hat-variables are available, we construct physics based numerical flux fluctuations by penalizing data against the incoming characteristics \eqref{eq:characteristics} at the element faces.
 
 \subsection{Boundary data}\label{subsec:boundary_hat_variables}
For $Z_\eta  > 0$, we define the characteristics
\begin{align}\label{eq:charac}
q_\eta  = \frac{1}{2}\left(Z_\eta  v_\eta  + T_\eta \right), \quad
p_\eta  = \frac{1}{2}\left(Z_\eta  v_\eta  -  T_\eta \right), \quad  \eta = l, m, n.
\end{align}
 Here, $q_\eta$ are the left going characteristics and $p_\eta$ are the right going characteristics.
We will  construct boundary data which satisfy the physical boundary conditions \eqref{eq:BC_General2} exactly and preserve the amplitude of the outgoing characteristics  $q_\eta $ \text{at}   $\xi = 0$, and  $p_\eta $ at $\xi = 1$.
 
To begin, define the hat-variables preserving the amplitude of outgoing characteristics
{\small
\begin{align}\label{eq:BC_hat_1bc}
{q}_\eta \left(\widehat{v}_\eta , \widehat{T}_\eta , Z_\eta \right) = {q}_\eta \left({v}_\eta , {T}_\eta , Z_\eta \right),  \quad \text{at} \quad \xi = 0,
\quad \text{and} \quad
{p}_\eta \left(\widehat{v}_\eta , \widehat{T}_\eta , Z_\eta \right) = {p}_\eta \left({v}_\eta , {T}_\eta , Z_\eta \right), \quad \text{at} \quad \xi = 1.
\end{align}
}
Since hat-variables also satisfy the physical boundary condition \eqref{eq:BC_General2}, we must have
{\small
\begin{align}\label{eq:BC_hat_2bc}
\frac{Z_{\eta}}{2}\left({1-\gamma_\eta }\right)\widehat{v}_\eta  -\frac{1+\gamma_\eta }{2} \widehat{T}_\eta  = 0, \quad \text{at} \quad \xi = 0,
\quad \text{and} \quad
\frac{Z_{\eta}}{2} \left({1-\gamma_\eta }\right)\widehat{v} _\eta + \frac{1+\gamma_\eta }{2}\widehat{T}_\eta  = 0, \quad \text{at} \quad \xi = 1.
\end{align}
}
The algebraic problem for the hat-variables, defined by equations \eqref{eq:BC_hat_1bc} and \eqref{eq:BC_hat_2bc},  has a unique solution, namely
{\small
\begin{align}\label{eq:data_hat}
\widehat{v}_\eta   = \frac{(1+\gamma_\eta )}{Z_{\eta}}q_\eta , \quad \widehat{T}_\eta   = {(1-\gamma_\eta )}q_\eta ,  \quad \text{at} \quad \xi = 0,
\quad \text{and} \quad
\widehat{v}_\eta  = \frac{(1+\gamma_\eta )}{Z_{\eta}}p_\eta ,  \quad \widehat{T}_\eta   = {-(1-\gamma_\eta )}p_\eta , \quad \text{at} \quad \xi = 1.
\end{align}
}
The expressions in \eqref{eq:data_hat} define a rule to update particle velocities and tractions on the physical boundaries $\xi = 0, 1$. That is 
{\small
\begin{align}\label{eq:boundary_data_hat}
v_\eta   = \widehat{v}_\eta  , \quad {T}_\eta  = \widehat{T}_\eta  , \quad \text{at} \quad \xi = 0,
\quad \text{and} \quad
v_\eta   = \widehat{v}_\eta   , \quad {T}_\eta  = \widehat{T}_\eta  , \quad \text{at} \quad \xi = 1.
\end{align}
}

The computation of the hat-variables for the external boundaries is summarized in Algorithm \ref{external_boundary_update} below.

\begin{algorithm}
\caption{Generalized Riemann solver for external element faces}\label{external_boundary_update}
\begin{algorithmic}[1]
\Procedure{ To compute the Riemann states $\widehat{T}_\eta $, $\widehat{v}_\eta $ with the inputs variables ${T}_\eta $, ${v}_\eta $, ${Z}_\eta $, $\gamma_\eta $, and $\xi = 0, 1$.}{}
\BState \emph{loop}: over the elements at  the external boundaries 
\BState \emph{loop}: for each $\eta = l, m, n$
\State \If {$Z_\eta  > 0$ and $\xi = 0$} from \eqref{eq:charac}  compute $q_\eta $, and
\[
\widehat{v}_\eta   = \frac{(1+\gamma_\eta )}{Z_{\eta}}q_\eta , \quad \widehat{T}_\eta   = {(1-\gamma_\eta )}q_\eta ,
\]
\EndIf
\If {$Z_\eta   > 0$ and $\xi = 1$}  from \eqref{eq:charac}  compute $p_\eta $, and
\[
\widehat{v}_\eta  = \frac{(1+\gamma_\eta )}{Z_{\eta}}p_\eta ,  \quad \widehat{T}_\eta   = {-(1-\gamma_\eta )}p_\eta ,
 \]
\EndIf
\EndProcedure
\end{algorithmic}
\end{algorithm}
  
 By construction, the hat-variables $\widehat{v}_\eta , \widehat{T}_\eta $,  satisfy the following algebraic identities:
\begin{subequations}\label{eq:identity_bc}
\small
\begin{equation}\label{eq:identity_1_bc}
{q}_\eta \left(\widehat{v}_\eta , \widehat{T}_\eta , Z_\eta \right) = {q}_\eta \left({v}_\eta , {T}_\eta , Z_\eta \right), \quad \text{at} \quad \xi = 0,
 \quad
{p}_\eta \left(\widehat{v}_\eta , \widehat{T}_\eta , Z_\eta \right) = {p}_\eta \left({v}_\eta , {T}_\eta , Z_\eta \right) , \quad \text{at} \quad \xi = 1,
\end{equation}
\begin{equation}\label{eq:identity_2_bc}
{q}_\eta ^2\left({v}_\eta , {T}_\eta , Z_\eta \right)-{p}_\eta ^2\left(\widehat{v}_\eta , \widehat{T}_\eta , Z_\eta \right) = Z_\eta \widehat{T}_\eta \widehat{v}_\eta ,  \quad \text{at} \quad \xi = 0, \quad {p}_\eta ^2\left({v}_\eta , {T}_\eta , Z_\eta \right) -{q}_\eta ^2\left(\widehat{v}_\eta , \widehat{T}_\eta , Z_\eta \right) = -Z_\eta \widehat{T}_\eta \widehat{v}_\eta ,  \quad \text{at} \quad \xi = 1,
\end{equation}
\begin{equation}\label{eq:identity_3_bc}
 \widehat{T}_\eta \widehat{v}_\eta  = \frac{1-\gamma_\eta ^2}{Z_{\eta}}{q}_\eta ^2\left({v}_\eta , {T}_\eta , Z_\eta \right) \ge 0,  \quad \text{at} \quad \xi = 0, \quad \widehat{T}_\eta \widehat{v}_\eta  = -\frac{1-\gamma_\eta ^2}{Z_{\eta}}{p}_\eta ^2\left({v}_\eta , {T}_\eta , Z_\eta \right) \le 0 , \quad \text{at} \quad \xi = 1.
\end{equation}
\end{subequations}
The algebraic identities \eqref{eq:identity_1_bc}--\eqref{eq:identity_3_bc} will be crucial in proving numerical stability.
Please see  \cite{DuruGabrielIgel2017} for  more details. 
\begin{lemma}\label{Lem:BTs_hat}
Consider the boundary term $\mathrm{BTs}$ defined in \eqref{eq:boundaryterm_101}, where $ v_{\eta} = \widehat{v}_{\eta}$,  $T_{\eta} = \widehat{T}_{\eta}$, with $\widehat{T}_{\eta}, \widehat{v}_{\eta}$ defined in \eqref{eq:data_hat}.  
The boundary term  $\mathrm{BTs}$ is never positive for all  $|\gamma_\eta | \le 1$, that is $\mathrm{BTs}\left(\widehat{v}_{\eta}, \widehat{T}_{\eta}\right)  \le 0$ where
{\small
 \begin{align}\label{eq:energy_estimate_BTs_hat}
\mathrm{BTs}\left(\widehat{v}_{\eta}, \widehat{T}_{\eta}\right)  &=   \sum_{\xi = q, r, s}\int_{\widetilde{\Gamma}}\left(\left(\sqrt{\xi_x^2 + \xi_y^2 + \xi_z^2}J \sum_{\eta = l,m,n}\widehat{T}_\eta \widehat{v}_\eta \right)\Big|_{\xi = 1}\right)\frac{dqdrds}{d\xi} \nonumber \\
& - \sum_{\xi = q, r, s}\int_{\widetilde{\Gamma}}\left(\left(\sqrt{\xi_x^2 + \xi_y^2 + \xi_z^2}J \sum_{\eta = l,m,n}\widehat{T}_\eta \widehat{v}_\eta \right)\Big|_{\xi = 0}\right)\frac{dqdrds}{d\xi}.
 \end{align}
 }
\end{lemma}
\begin{proof}
 With $Z_{\eta} > 0$, for any  $|\gamma_\eta | \le 1$, from \eqref{eq:identity_3_bc} we have
{
 \begin{align}
 \widehat{T}_\eta \widehat{v}_\eta  = \frac{1-\gamma_\eta ^2}{Z_{\eta}}{q}_\eta ^2\left({v}_\eta , {T}_\eta , Z_\eta \right) \ge 0,  \quad \text{at} \quad \xi = 0, \quad \widehat{T}_\eta \widehat{v}_\eta  = -\frac{1-\gamma_\eta ^2}{Z_{\eta}}{p}_\eta ^2\left({v}_\eta , {T}_\eta , Z_\eta \right) \le 0 , \quad \text{at} \quad \xi = 1.
  \end{align}
}
  \begin{equation}\label{eq:boundaryterm_00}
\begin{split}
 \mathrm{BTs}\left(v , T \right) =&-\sum_{\xi = q, r, s}\int_{\widetilde{\Gamma}}\left(\sqrt{\xi_x^2 + \xi_y^2 + \xi_z^2}J\sum_{\eta = l,m,n}\left({\frac{1-\gamma_\eta ^2}{Z_{\eta}}{p}_\eta ^2\left({v}_\eta , {T}_\eta , Z_\eta \right) }\right)\Big|_{\xi = 1}\right)\frac{dqdrds}{d\xi}\\
  &-\sum_{\xi = q, r, s}\int_{\widetilde{\Gamma}}\left(\sqrt{\xi_x^2 + \xi_y^2 + \xi_z^2}J\sum_{\eta = l,m,n}\left({ \frac{1-\gamma_\eta ^2}{Z_{\eta}}{q}_\eta ^2\left({v}_\eta , {T}_\eta , Z_\eta \right) }\right)\Big|_{\xi = 0}\right)\frac{dqdrds}{d\xi}.
\end{split}
  \end{equation}
We must have $\mathrm{BTs}\left(\widehat{v}_{\eta}, \widehat{T}_{\eta}\right)  \le 0$ for all $|\gamma_\eta | \le 1$ .
\end{proof}
 
Lemma \ref{Lem:BTs_hat}  is completely analogous to Lemma \ref{Lem:BTs}.  
We will now formulate a result equivalent to  Theorem \ref{Theo:energy_estimate_BC_00}.
  \begin{theorem}\label{Theo:energy_estimate_BC_hat}
 Consider the  transformed  equation of motion \eqref{eq:gen_hyp_transformed} subject to the boundary condition \eqref{eq:boundary_data_hat}. The solutions of the transformed  equation \eqref{eq:gen_hyp_transformed} with \eqref{eq:boundary_data_hat}  satisfy 
 \begin{align}\label{eq:energy_estimate_bc_hat}
\frac{d}{dt}E(t)  =  \mathrm{BTs}\left(\widehat{v}, \widehat{T}\right) \le 0.
 \end{align}
 \end{theorem}
 \subsection{Interface data}\label{subsec:boundary_hat_variables}
To begin,  define the outgoing  characteristics at the interface
\begin{align}\label{eq:characteristics_out}
q_\eta ^{+} = \frac{1}{2}\left(Z_\eta ^{+} v^{+}_\eta  + T_\eta ^{+}\right), \quad
p_\eta ^{-} = \frac{1}{2}\left(Z_\eta ^{-} v^{-}_\eta  -  T_\eta ^{-}\right), \quad  \eta = l, m, n,
\end{align}
where $Z_\eta ^{\pm} > 0$ are the impedances. 
We define the hat-variables  preserving the amplitude of  the outgoing characteristics at the interface
 \begin{equation}\label{eq:hat_variables}
\widehat{q}_\eta ^{+} \left(\widehat{v}_\eta ^{+}, \widehat{T}_\eta ^{+}, Z_\eta ^{+}\right)   = {q}_\eta ^{+} \left({v}_\eta ^{+}, {T}_\eta ^{+}, Z_\eta ^{+} \right), 
\quad
\widehat{p}_\eta ^{-} \left(\widehat{v}_\eta ^{-}, \widehat{T}_\eta ^{-}, Z_\eta ^{-}\right)  = {p}_\eta ^{-}\left({v}_\eta ^{-}, {T}_\eta ^{-}, Z_\eta ^{-} \right).
 \end{equation}
The hat-variables also satisfy the  interface conditions \eqref{eq:elastic_elastic_interface} exactly. For each setup, given ${q}_\eta ^{+},  {p}_\eta ^{-} $, the procedure will solve \eqref{eq:hat_variables}, and  \eqref{eq:elastic_elastic_interface} for the hat-variables, $\widehat{v}_\eta ^{\pm}, \widehat{T}_\eta ^{\pm}$. 
We consider  an interface separating two elastic solids. The hat-variables must satisfy \eqref{eq:hat_variables} and the interface conditions, \eqref{eq:elastic_elastic_interface}, that is, force balance, no opening and no slip conditions.
Combining the two equations in \eqref{eq:hat_variables}  and ensuring force balance, $\widehat{T}_\eta ^{-} = \widehat{T}_\eta ^{+} = \widehat{T}_\eta $,
we have
\begin{align}\label{eq:radiation_damping_line}
\widehat{T}_\eta  + \alpha_\eta  \lJump{\widehat{v}_\eta \rJump}  = \Phi_\eta ,    \quad \Phi_\eta  = \alpha_\eta \left(\frac{2}{Z_\eta ^{+}}q_\eta ^{+} - \frac{2}{Z_\eta ^{-}}p_\eta ^{-}\right), \quad  \alpha_\eta  = \frac{Z_\eta ^{+}Z_\eta ^{-}}{Z_\eta ^{+} + Z_\eta ^{-}} >0, \quad \eta = l, m, n.
\end{align}
 Furthermore, enforcing no opening and no slip conditions, $\lJump{\widehat{v}_\eta \rJump} = 0$,  in \eqref{eq:radiation_damping_line} gives
\begin{align}\label{eq:traction_jump_vel_hat}
\widehat{T}_\eta   = \Phi_\eta , \quad  \lJump{\widehat{v}_\eta \rJump} = 0,  \quad \eta = l, m, n.
\end{align}
We can now define explicitly the hat-variables corresponding to the particle velocities and tractions  
 \begin{equation}\label{eq:hat-variables_Elastic_Elastic}
 \widehat{T}_\eta ^{-} = \widehat{T}_\eta ^{+} = \widehat{T}_\eta , \quad \widehat{v}_\eta ^{+} = \frac{2p_\eta ^{-}+\widehat{T}_\eta }{Z_\eta ^{-}} , \quad 
 \widehat{v}_\eta ^{-} = \frac{2q_\eta ^{+}-\widehat{T}_\eta }{Z_\eta ^{+}} , \quad \quad \eta = l, m, n. 
 \end{equation}

 Note that we have equivalently redefined the physical interface condition \eqref{eq:elastic_elastic_interface} as  follows
 \begin{equation}\label{eq:Hat_Functions_BC}
 v_{\eta}^{\pm} = \widehat{v}_{\eta}^{\pm}, \quad T_{\eta}^{\pm} = \widehat{T}_{\eta}, \quad \eta = l, m, n.
 \end{equation}
The procedure to compute  the hat-variables for the internal element boundaries is given in Algorithm \ref{internal_boundary_update}.
\begin{algorithm}
\caption{Generalized Riemann solver for internal element faces}\label{internal_boundary_update}
\begin{algorithmic}[1]
\Procedure{ To compute the Riemann states $\widehat{T}_\eta ^{\pm}$, $\widehat{v}_\eta ^{\pm}$ with the inputs variables ${T}_\eta ^{\pm}$, ${v}_\eta ^{\pm}$ and ${Z}_\eta ^{\pm}$.}{}
\BState \emph{loop}: over the elements points on the interface
\BState \emph{loop}: for each $\eta = l, m, n$
\State \If {$Z_\eta ^{-} > 0$ and $Z_\eta ^{+} > 0$} from \eqref{eq:characteristics} and \eqref{eq:radiation_damping_line} compute $\Phi_\eta $, and
\begin{align*}
\widehat{T}_\eta ^{\pm}=\widehat{T}_\eta  = \Phi_\eta , \quad  \lJump{\widehat{v}_\eta \rJump} = 0, \quad
 \widehat{v}_\eta ^{+} = \frac{2p_\eta ^{-}+ \Phi_\eta  }{Z_\eta ^{-}}, 
   \quad
 \widehat{v}_\eta ^{-} = \frac{2q_\eta ^{+}- \Phi_\eta }{Z_\eta ^{+}} .
 \end{align*}
\EndIf
\EndProcedure
\end{algorithmic}
\end{algorithm}
 
By construction, the hat-variables $\widehat{v}_\eta ^{\pm}$, $\widehat{T}_\eta ^{\pm}$ satisfy the following algebraic identities:
\begin{subequations}\label{eq:identity}
\small
\begin{equation}\label{eq:identity_1}
{q}_\eta  \left(\widehat{v}_\eta ^{+}, \widehat{T}_\eta ^{+}, Z_\eta ^{+}\right)   = {q}_\eta  \left({v}_\eta ^{+}, {T}_\eta ^{+}, Z_\eta ^{+} \right), 
\quad
{p}_\eta  \left(\widehat{v}_\eta ^{-}, \widehat{T}_\eta ^{-}, Z_\eta ^{-}\right)  = {p}_\eta \left({v}_\eta ^{-}, {T}_\eta ^{-}, Z_\eta ^{-} \right), 
\end{equation}
\begin{equation}\label{eq:identity_2}
\left(q_\eta ^2  \left({v}_\eta ^{+}, {T}_\eta ^{+}, Z_\eta ^{+}\right)\right)-{p}_\eta ^2 \left(\widehat{v}_\eta ^{+}, \widehat{T}_\eta ^{+}, Z_\eta ^{+}\right) = Z_\eta ^{+}\widehat{T}_\eta \widehat{v}_\eta ^{+}, \quad p^2_\eta \left({v}_\eta ^{-}, {T}_\eta ^{-}, Z_\eta ^{-} \right) -{q}^2_\eta \left(\widehat{v}_\eta ^{-}, \widehat{T}_\eta ^{-}, Z_\eta ^{-}\right) = -Z_\eta ^{-}\widehat{T}_\eta \widehat{v}_\eta ^{-},
\end{equation}
\begin{equation}\label{eq:identity_3}
\frac{1}{Z_\eta ^+}\left(q^2_\eta \left({v}_\eta ^{+}, {T}_\eta ^{+}, Z_\eta ^{+} \right)- {p}^2_\eta \left(\widehat{v}_\eta ^{+}, \widehat{T}_\eta ^{+}, Z_\eta ^{+}\right)\right)  + \frac{1}{Z_\eta ^-}\left(p^2_\eta \left({v}_\eta ^{-}, {T}_\eta ^{-}, Z_\eta ^{-} \right) -{q}^2_\eta \left(\widehat{v}_\eta ^{-}, \widehat{T}_\eta ^{-}, Z_\eta ^{-}\right)\right) =  \widehat{T}_\eta  \lJump{\widehat{v}_\eta \rJump} \equiv 0,
\end{equation}
\end{subequations}
 The identities \eqref{eq:identity_1}--\eqref{eq:identity_2} can be easily verified, see  \cite{DuruGabrielIgel2017}, and will be useful in the prove of numerical stability.
 
 We can now formulate a result equivalent to Theorem \ref{Theo:linear_friction_energy_estimate}.
  \begin{theorem}
Consider the modeling domain decomposed into two elements, $\Omega = \Omega_{-} \cup \Omega_{+}$, as in  Figure \ref{fig:split_blocks}. Denote the corresponding fields, material parameters and energies in the positive/negative sides of the interface with the superscripts $+/-$.  Let the transformed  equation of motion \eqref{eq:gen_hyp_transformed}  be defined on each element, with the two adjacent elements connected at the interface through the interface condition  \eqref{eq:Hat_Functions_BC}.
The sum of the energies satisfies
 \begin{align}\label{eq:energy_estimate_fault_30}
\frac{d}{dt}\left(E^{-}+ E^{+}\right)   =  IT_s(\widehat{v}^{\pm}, \widehat{T}^{\pm}) \equiv 0.
 \end{align}
\end{theorem}

 \subsection{Physics based flux fluctuations}
%
%
The next step is to construct fluctuations by penalizing data, that is hat-variables, against the ingoing characteristics only.
  If $Z_\eta > 0$, then we have
    \begin{align*}
  &{G}_\eta = \frac{1}{2} {Z}_\eta \left({v}_\eta - \widehat{{v}}_\eta \right)+ \frac{1}{2}\left({T}_\eta  - \widehat{{T}}_\eta \right)\Big|_{\xi = 1},   \quad \widetilde{G}_\eta :=\frac{1}{{Z}_\eta}{G}_\eta  = \frac{1}{2} \left({v}_\eta - \widehat{{v}}_\eta \right)+ \frac{1}{2 {Z}_\eta}\left({T}_\eta  - \widehat{{T}}_\eta \right)\Big|_{\xi = 1},
  \\
  &{G}_\eta  = \frac{1}{2} {Z}_\eta\left({v}_\eta  - \widehat{{v}}_\eta \right)- \frac{1}{2}\left({T}_\eta  - \widehat{{T}}_\eta \right)\Big|_{\xi = 0}, \quad \widetilde{G}_\eta:=\frac{1}{{Z}_\eta }{G}_\eta  = \frac{1}{2} \left({v}_\eta - \widehat{{v}}_\eta \right)- \frac{1}{2 {Z}_\eta }\left({T}_\eta  - \widehat{{T}}_\eta \right)\Big|_{\xi = 0}.
  \end{align*}
  The fluctuations are computed in the transformed coordinates $l,m,n$ . We will now  rotate them to the physical coordinates $x,y,z$, having
  \begin{align}\label{eq:rotate_back_forth}
  {\mathbf{G}} := \begin{pmatrix}
{G}_{x} \\
{G}_{y} \\
{G}_{z}
\end{pmatrix}
 = \mathbf{R}^T\begin{pmatrix}
{G}_{n} \\
{G}_{m} \\
{G}_{l}
\end{pmatrix},
\quad 
\widetilde{\mathbf{G}}:= 
  \begin{pmatrix}
\widetilde{G}_{x} \\
\widetilde{G}_{y} \\
\widetilde{G}_{z}
\end{pmatrix} = \mathbf{R}^T\begin{pmatrix}
\widetilde{G}_{n} \\
\widetilde{G}_{m} \\
\widetilde{G}_{l}
\end{pmatrix}.
\end{align}
  Note that
{
\small
\begin{equation}\label{eq:identity_pen}
\begin{split}
 \left(\mathbf{v}^T \mathbf{G} - \mathbf{T}^T \widetilde{\mathbf{G}} + \mathbf{v}^T\mathbf{T}\right)\Big|_{\xi = 0} &= \sum_{\eta = l,m,n} \left( v_\eta  G_\eta    - \frac{1}{Z_\eta  }T_\eta   G_\eta + v_\eta  T_\eta \right)\Big|_{\xi = 0} \\
 &= \sum_{\eta = l,m,n} \frac{1}{Z_\eta  }\left(|G_\eta |^2 + p^2_\eta \left(v_\eta , T_\eta , Z_\eta \right) - {q}^2_\eta \left(\widehat{v}_\eta , \widehat{T}_\eta , Z_\eta \right)\right)\Big|_{\xi = 0}\\
&= \sum_{\eta = l,m,n} \left(\frac{1}{Z_\eta  }|G_\eta |^2  + \widehat{T}_\eta \widehat{v}_\eta \right)\Big|_{\xi = 0},\\
\left(\mathbf{v}^T \mathbf{G} + \mathbf{T}^T \widetilde{\mathbf{G}} - \mathbf{v}^T\mathbf{T}\right)\Big|_{\xi = 1} &= \sum_{\eta = l,m,n} \left(v_\eta  G_\eta   + \frac{1}{Z_\eta  }T_\eta   G_\eta - v_\eta  T_\eta \right)\Big|_{\xi = 1} \\
&= \sum_{\eta = l,m,n} \frac{1}{Z_\eta  }\left(|G_\eta |^2 + q^2_\eta \left(v_\eta , T_\eta , Z_\eta \right) - {p}^2_\eta \left(\widehat{v}_\eta , \widehat{T}_\eta , Z_\eta \right)\right)\Big|_{\xi = 1}\\
&= \sum_{\eta = l,m,n} \left(\frac{1}{Z_\eta  }|G_\eta |^2  - \widehat{T}_\eta \widehat{v}_\eta \right)\Big|_{\xi = 1}.
\end{split}
\end{equation}
}

The physics based flux fluctuations  obeying the eigen--structure of the elastic wave equation are given by
  {
  \small
   \begin{align}
&\mathbf{FL} = \left[{G}_x, {G}_y,  {G}_z,         -{n_x}\widetilde{{G}}_x, -{n_y}\widetilde{{G}}_y , -{n_z}\widetilde{{G}}_z, -\left({n_y}\widetilde{{G}}_x + {n_x}\widetilde{{G}}_y\right), -\left({n_z}\widetilde{{G}}_x + {n_x}\widetilde{{G}}_z\right), -\left({n_z}\widetilde{{G}}_y + {n_y}\widetilde{{G}}_z\right)   \right]^T,\\ \nonumber
&\mathbf{FR} = \left[{G}_x, {G}_y,  {G}_z,         {n_x}\widetilde{{G}}_x, {n_y}\widetilde{{G}}_y , {n_z}\widetilde{{G}}_z, \left({n_y}\widetilde{{G}}_x + {n_x}\widetilde{{G}}_y\right), \left({n_z}\widetilde{{G}}_x + {n_x}\widetilde{{G}}_z\right), \left({n_z}\widetilde{{G}}_y + {n_y}\widetilde{{G}}_z\right)   \right]^T.
      \end{align}
  Note that
  \begin{equation}\label{eq:scalar_product_flux}
  \mathbf{Q}^T\mathbf{FL} = \mathbf{v}^T \mathbf{G} - \mathbf{T}^T \widetilde{\mathbf{G}} , \quad \mathbf{Q}^T\mathbf{FR} = \mathbf{v}^T \mathbf{G} + \mathbf{T}^T \widetilde{\mathbf{G}}.
  \end{equation}
   }
   The elemental  weak form reads
{\small
\begin{align}\label{eq:transf_elemental_weak_form_velocity}
\int_{\widetilde{\Omega}} \boldsymbol{\phi}^T\widetilde{\mathbf{P}}^{-1} \frac{\partial }{\partial t} \mathbf{Q} dqdrds  &= \int_{\widetilde{\Omega}} \boldsymbol{\phi}^T\left(\div \mathbf{F} \left(\mathbf{Q} \right) + \mathbf{B}\left(\grad\mathbf{Q}\right) \right) dqdrds \nonumber \\
- & \sum_{\xi = q, r, s}\int_{\widetilde{\Gamma}}\sqrt{\xi_x^2 + \xi_y^2 + \xi_z^2}\left(\left[\boldsymbol{\phi}^T\mathbf{FL}\right]_{\xi = 0}  + \left[\boldsymbol{\phi}^T \mathbf{FR}\right]_{\xi = 1} \right) J\frac{dqdrds}{d\xi}.
\end{align}
}
Here, the variable $\xi = q, r, s$,  indicates the directions where the flux is computed. The flux fluctuation vectors, $\mathbf{FL}, \mathbf{FR}$ enforce weakly the boundary conditions \eqref{eq:boundary_data_hat} and interface interface conditions \eqref{eq:Hat_Functions_BC}   at the element faces.
Note, however, that since we have not introduced any numerical approximations, the fluctuations vanish identically, ${G}_\eta ^{\pm} = \widetilde{G}_\eta ^{\pm} =0 $, for the exact solutions that satisfy the PDE and the boundary and interface conditions,  we have $\mathbf{FL}\equiv 0$, $\mathbf{FR} \equiv 0$.

Introduce the fluctuation term
  \begin{align}\label{eq:weak_fluctuate}
  F_{luc}\left({G},  Z\right) = - \sum_{\xi = q, r, s}\int_{\widetilde{\Gamma}}\left(\left(\sqrt{\xi_x^2 + \xi_y^2 + \xi_z^2}J\sum_{\eta = l,m,n}\frac{1}{Z_{\eta}}|{G}_\eta |^2\right)\Big|_{\xi = 0} \right)\frac{dqdrds}{d\xi} \nonumber\\
  - \sum_{\xi = q, r, s}\int_{\widetilde{\Gamma}}\left(\left(\sqrt{\xi_x^2 + \xi_y^2 + \xi_z^2}J\sum_{\eta = l,m,n}\frac{1}{Z_{\eta}}|{G}_\eta |^2\right)\Big|_{\xi = 1} \right)\frac{dqdrds}{d\xi}.
  \end{align}
   For the two elements model we introduce the external boundary terms 
 \begin{equation}\label{eq:weak_BT1}
 \begin{split}
\mathrm{BTs}\left(\widehat{v}_{\eta}^{\pm}, \widehat{T}_{\eta}^{\pm}\right)  = 
& \sum_{\xi =  r, s}\int_{\widetilde{\Gamma}}\left(\left(\sqrt{\xi_x^2 + \xi_y^2 + \xi_z^2}J \sum_{\eta = l,m,n}\widehat{T}_\eta^{\pm}  \widehat{v}_\eta^{\pm}  \right)\Big|_{\xi = 1}\right)\frac{dqdrds}{d\xi} \\
& - \sum_{\xi =  r, s}\int_{\widetilde{\Gamma}}\left(\left(\sqrt{\xi_x^2 + \xi_y^2 + \xi_z^2}J \sum_{\eta = l,m,n}\widehat{T}_\eta^{\pm}  \widehat{v}_\eta^{\pm}  \right)\Big|_{\xi = 0}\right)\frac{dqdrds}{d\xi}  \\
&\pm \int_{\widetilde{\Gamma}}\left(\left(\sqrt{q_x^2 + q_y^2 + q_z^2}J \sum_{\eta = l,m,n}\widehat{T}_\eta^{\pm}  \widehat{v}_\eta^{\pm}  \right)\Big|_{q = 0}\right){drds} ,
\end{split}
\end{equation}
 and the interface term
  \begin{align*}
  IT_s(\widehat{v}^{\pm}, \widehat{T}^{\pm}) = - \sum_{\eta = l,m,n}\int_{\widetilde{\Gamma}}\sqrt{q_x^2 + q_y^2 + q_z^2}{\widehat{T}}_\eta  \lJump{{\widehat{v}}_\eta \rJump} J {drds} \equiv 0.
\end{align*}
      \begin{remark}
   Note that by construction the interface terms vanish identically $ IT_s(\widehat{v}^{\pm}, \widehat{T}^{\pm}) \equiv 0$, by inspection the fluctuation terms are negative semi-definite $F_{luc}\left({G},  Z\right) \le 0$, and vanishes identically $F_{luc}\left({G},  Z\right) \equiv 0$ for the exact solutions satisfying the PDE, the boundary and interface conditions. By  Lemma \ref{Lem:BTs_hat} the boundary terms are negative semi-definite $\mathrm{BTs}\left(\widehat{v}_{\eta}^{-}, \widehat{T}_{\eta}^{-}\right) \le 0$.
  \end{remark} 
   We can now prove the theorem.
 \begin{theorem}\label{theo:weak_form_statbility}
The elemental weak form \eqref{eq:transf_elemental_weak_form_velocity} satisfies the energy equation
   {
   \small
 \begin{equation}\label{eq:weak_energy}
\begin{split}
&\frac{d}{dt} \left({E}^{-}(t) + {E}^{+}(t)\right)=   IT_s\left(\widehat{v}^{\pm},  \widehat{T}^{\pm}\right)  +  BT_s\left(\widehat{v}^{-}, \widehat{T}^{-}\right) + BT_s\left(\widehat{v}^{+}, \widehat{T}^{+}\right) +   F_{luc}\left({G} ^{-},  Z^{-}\right) +   F_{luc}\left({G} ^{+},  Z^{+}\right) \le 0.
\end{split}
\end{equation}
}
\end{theorem}
\begin{proof}
In \eqref{eq:transf_elemental_weak_form_velocity} replace the test function $\boldsymbol{\phi}$ by the solution $\mathbf{Q}$,  and integrate the conservative flux term only by parts. The volume term vanishes, having 
{
\begin{equation}\label{eq:transf_elemental_weak_form_velocity_proof}
\begin{split}
\frac{d }{d t}\int_{\widetilde{\Omega}}\frac{1}{2} [\mathbf{Q} ^T\widetilde{\mathbf{P}}^{-1}  \mathbf{Q}] dqdrds  = 
- & \sum_{\xi = q, r, s}\int_{\widetilde{\Gamma}}\left(\left(J\sqrt{\xi_x^2 + \xi_y^2 + \xi_z^2}\left(\mathbf{v}^T \mathbf{G} - \mathbf{T}^T \widetilde{\mathbf{G}} + \mathbf{v}^T\mathbf{T}\right)\right)_{\xi = 0} \right) \frac{dqdrds}{d\xi}\\
- & \sum_{\xi = q, r, s}\int_{\widetilde{\Gamma}}\left(\left(J\sqrt{\xi_x^2 + \xi_y^2 + \xi_z^2}\left(\mathbf{v}^T \mathbf{G} + \mathbf{T}^T \widetilde{\mathbf{G}} - \mathbf{v}^T\mathbf{T}\right)\right)_{\xi = 1}\right) \frac{dqdrds}{d\xi}.
\end{split}
\end{equation}
}
Using \eqref{eq:identity_pen} in the right hand side of \eqref{eq:transf_elemental_weak_form_velocity_proof} gives
{
\begin{equation}\label{eq:transf_elemental_weak_form_velocity_proof0}
\begin{split}
\frac{d }{d t}\int_{\widetilde{\Omega}}\frac{1}{2} [\mathbf{Q} ^T\widetilde{\mathbf{P}}^{-1}  \mathbf{Q}] dqdrds  &= 
- \sum_{\xi = q, r, s}\int_{\widetilde{\Gamma}}\left(\left(J\sqrt{\xi_x^2 + \xi_y^2 + \xi_z^2}\sum_{\eta = l,m,n} \left(\frac{1}{Z_\eta  }|G_\eta |^2  + \widehat{T}_\eta \widehat{v}_\eta \right)\right)_{\xi = 0} \right) \frac{dqdrds}{d\xi}\\
&-  \sum_{\xi = q, r, s}\int_{\widetilde{\Gamma}}\left(\left(J\sqrt{\xi_x^2 + \xi_y^2 + \xi_z^2}\sum_{\eta = l,m,n} \left(\frac{1}{Z_\eta  }|G_\eta |^2  - \widehat{T}_\eta \widehat{v}_\eta \right)\right)_{\xi = 1}\right) \frac{dqdrds}{d\xi}.
\end{split}
\end{equation}
}
Collecting contributions from the two elements on both sides of the interface,  and using the identities \eqref{eq:identity_bc} and \eqref{eq:identity} gives the energy equation \eqref{eq:weak_energy}.
\end{proof}

\begin{remark}
It is particularly noteworthy that we have completely avoided explicit eigen-decomposition of the compound coefficient matrices
\[
 \widetilde{\mathbf{A}}_n = n_x \widetilde{\mathbf{A}}_x + n_y \widetilde{\mathbf{A}}_y + n_z \widetilde{\mathbf{A}}_z,  \quad \widetilde{\mathbf{A}}_\xi = {\mathbf{P}}{\mathbf{A}}_\xi
\]
which will be   necessary for some DG methods  \cite{DumbserKaeser2006, PeltiesdelaPuenteAmpueroBrietzkeKaser2012} based on classical Godunov numerical flux \cite{Godunov1959}. 
Here $\mathbf{P}$ is the parameter matrix defined in \eqref{eq:material_coeff},  ${\mathbf{A}}_\xi$ are the non-dimensional coefficient matrices defined in  \eqref{eq:elastic_coeff} and $ \mathbf{n} = [n_x, n_y, n_z]^T$ is the outward unit normal on the boundary defined in \eqref{eq:unit_normal}.
The eigen-decomposition of the matrices $ \widetilde{\mathbf{A}}_n$ on each element face and for every element can be cumbersome in general curvilinear elements, and in particular in heterogeneous anisotropic elastic medium.
\end{remark}
 In the next section, we will introduce the discontinuous Galerkin approximation and prove numerical stability by  deriving discrete energy estimates analogous to \eqref{eq:weak_energy}.


\section{The discontinuous Galerkin approximation}\label{sec:DG_method}
Inside the transformed  element  $(q, r, s) \in \widetilde{\Omega} = [0, 1]^3$, approximate the elemental solution and the conservative flux term by  polynomial interpolants of degree $P$,  and write 
\begin{equation}\label{eq:variables_elemental}
\bar{\mathbf{Q}}(q,r,s, t) =  \sum_{i = 1}^{P+1} \sum_{j = 1}^{P+1}\sum_{k = 1}^{P+1}\bar{\mathbf{Q}}_{ijk}(t) \boldsymbol{\phi}_{ijk}(q, r, s), \quad \mathbf{F} \left(\bar{\mathbf{Q}} \right) =   \sum_{i = 1}^{P+1} \sum_{j = 1}^{P+1}\sum_{k = 1}^{P+1} \mathbf{F}_{ijk} \left(\bar{\mathbf{Q}}(t)\right) \boldsymbol{\phi}_{ijk}\left(q, r, s\right),
\end{equation}
where $\mathbf{F}_{ijk} \left(\bar{\mathbf{Q}}(t)\right)  = \mathbf{F} \left(\bar{\mathbf{Q}}_{ijk}(t)\right)$,   $\bar{\mathbf{Q}}_{ijk}(t)$ are the evolving elemental degrees of freedom to be determined, 
 and $ \boldsymbol{\phi}_{ijk}(q, r, s)$ are the $ijk$-th interpolating polynomial. We consider tensor products of  nodal basis with $ \boldsymbol{\phi}_{ijk}(q, r, s) = \mathcal{L}_i(q)\mathcal{L}_j (r)\mathcal{L}_k(s)$,  where $\mathcal{L}_i(q)$, $\mathcal{L}_j (r)$, $\mathcal{L}_k(s)$, are one dimensional nodal interpolating Lagrange polynomials of degree $P$, with 
\begin{equation*}\label{eq:delta}
\begin{split}
&\mathcal{L}_i(q_m)= \left \{
\begin{array}{rl}
1 \quad {}  \quad {}& \text{if} \quad i = m ,\\
0 \quad {}  \quad {}& \text{if} \quad i \ne m.
\end{array} \right.
\end{split}
\end{equation*}
The interpolating nodes $q_m$, $m = 1, 2, \dots, P+1$, are the nodes of a Gauss quadrature with
\begin{equation}\label{eq:quad_rule_3D}
 \sum_{i = 1}^{P+1} \sum_{j = 1}^{P+1}  \sum_{k = 1}^{P+1}f(q_i, r_j , s_k)h_ih_\eta h_k \approx \int_{\widetilde{\Omega}}f(q, r,s) dq dr ds,
\end{equation}
where $h_i > 0$, $h_j >0$, $h_k>0$, are the quadrature weights.
We will only use quadrature rules such  that for all polynomial integrand $f(q)$ of degree $\le 2P-1$, the corresponding one space dimensional rule is exact, 
 $\sum_{i = 1}^{P+1} f(q_i)h_i = \int_{0}^{1}f(q) dq.$
Admissible candidates are the Gauss-Legendre-Lobatto quadrature rule with GLL nodes, the Gauss-Legendre quadrature rule with GL nodes and the Gauss-Legendre-Radau quadrature rule with GLR nodes. 
While both endpoints, $q = 0, 1$, are part of  GLL quadrature nodes, the GLR quadrature contains only the first endpoint $q = 0$ as a node. Lastly,  for the GL quadrature, both endpoints, $q = 0, 1$, are not quadrature nodes. Note that when an endpoint is not a quadrature node, $q_1 \ne 0$ or $q_{P+1} \ne 1$, extrapolation is needed to compute numerical fluxes at the element boundaries, $q = 0, 1$. We also remark that the GLL quadrature rule is exact for polynomial integrand of degree $2P-1$, GLR quadrature rule is exact for polynomial integrand of degree $2P$, and GL quadrature rule is exact for polynomial integrand of  degree $2P+1$.

 \subsection{The spectral difference approximation}
 Introduce the square  matrices $ H, A \in \mathbb{R}^{(P+1) \times (P+1)}$, defined by
{\small
\begin{equation}
 H = \mathrm{diag}[h_1, h_2, \cdots, h_{P+1}], \quad A_{ij} = \sum_{m = 1}^{P+1} h_m \mathcal{L}_i(q_m)  {\mathcal{L}_j^{\prime}(q_m)} = \int_{0}^{1}\mathcal{L}_i(q)  {\mathcal{L}_j^{\prime}(q)} dq.
\end{equation}
}
Note that the matrix
\begin{equation}\label{eq:derivative_operator_1d}
D = H^{-1} A \approx \frac{\partial}{\partial q},
\end{equation}
is a one space dimensional spectral difference approximation of the first derivative.

Using the fact that the quadrature rule is exact  for all polynomial  integrand of degree $\le 2P-1$ 
  implies that
\begin{equation}\label{eq:sbp_property_a}
A + A^T =  B,
\quad
B_{ij} = \mathcal{L}_i(1)  \mathcal{L}_j (1)- \mathcal{L}_i(0)  \mathcal{L}_j (0).
\end{equation}
Equations \eqref{eq:derivative_operator_1d} and \eqref{eq:sbp_property_a} are the discrete equivalence of the integration-by-parts property.
If both endpoints $q = 0,1$ are quadrature nodes and we consider nodal bases  then we have $B = \text{diag}[-1, 0,0, \dots, 0, 1].$ 
The matrix $B$ projects the nodal degrees of freedom to element faces.

The one space dimensional derivative operator \eqref{eq:derivative_operator_1d} can be extended to higher space dimensions using the Kronecker products $\otimes$, having
\begin{equation}\label{eq:derivative_operator_3d}
\mathbf{D}_q = \left(I_{9}\otimes D\otimes I\otimes I\right), \quad \mathbf{D}_r = \left(I_{9}\otimes I\otimes D\otimes I\right), \quad \mathbf{D}_s = \left(I_{9}\otimes I\otimes I\otimes D\right),
\end{equation}
\begin{equation}\label{eq:norm_operator_3d}
\mathbf{H}_q = \left(I_{9}\otimes H\otimes I\otimes I\right), \quad \mathbf{H}_r = \left(I_{9}\otimes I\otimes H\otimes I\right), \quad \mathbf{H}_s = \left(I_{9}\otimes I\otimes I\otimes H\right), \quad \mathbf{H} = \mathbf{H}_q\mathbf{H}_r\mathbf{H}_s.
\end{equation}
Here, $I$ is the $(P+1)\times(P+1)$ identity matrix, and $I_9$ is the $9\times 9$ identity matrix. Note that since $\mathbf{H}_q, \mathbf{H}_r, \mathbf{H}_s$ are diagonal matrices then $\mathbf{H}$ is also a diagonal matrix. The matrix $\mathbf{H}$ is independent of the order of the matrix products,  that is,  $\mathbf{H} = \mathbf{H}_q\mathbf{H}_r\mathbf{H}_s = \mathbf{H}_r\mathbf{H}_q\mathbf{H}_s  = \mathbf{H}_s\mathbf{H}_r\mathbf{H}_q$.

We also introduce the projection matrices
{
\begin{displaymath}
\mathbf{e}_q(\xi) = \left(I_{9}\otimes \boldsymbol{e}(\xi)\otimes I\otimes I\right), \quad \mathbf{e}_r(\xi) = \left(I_{9}\otimes I\otimes \boldsymbol{e}(\xi)\otimes I\right), \quad \mathbf{e}_s(\xi) = \left(I_{9}\otimes I\otimes I\otimes \boldsymbol{e}(\xi)\right), \quad  \mathbf{B}_{\eta}(\psi, \xi) = \mathbf{e}_{\eta}(\psi) \mathbf{e}_{\eta}^T(\xi),
\end{displaymath}
}
where
\begin{displaymath}
e(\xi) = [\mathcal{L}_1(\xi), \mathcal{L}_2(\xi), \dots, \mathcal{L}_m(\xi), \mathcal{L}_{m+1}(\xi)]^T.
\end{displaymath}
 \subsection{The semi--discrete approximation}
We will now make a classical Galerkin approximation, by choosing test functions in the same space as the basis functions.  Thus, replacing ${\mathbf{Q}}(q,r,s, t)$ by $\bar{\mathbf{Q}}(q,r,s, t)$ and $\mathbf{F} \left({\mathbf{Q}} \right) $ by $\mathbf{F} \left(\bar{\mathbf{Q}} \right) $ in \eqref{eq:transf_elemental_weak_form_velocity}, and approximating all integrals with the corresponding  quadrature rules yields the semi-discrete equation,
  \begin{align}\label{eq:gen_hyp_transformed_discrete}
\widetilde{\mathbf{P}}^{-1} \frac{d }{d t} \bar{\mathbf{Q}} = \grad_{D}\bullet {\mathbf{F} \left(\bar{\mathbf{Q}} \right)} + \sum_{\xi= q, r, s}\mathbf{B}_\xi\left(\grad_{D}\bar{\mathbf{Q}}\right) - \mathbf{Flux}\left(\bar{\mathbf{Q}}\right),
\end{align}
for the evolving degrees of freedom, $ \bar{\mathbf{Q}}(t) = [\bar{\mathbf{Q}}_{ijk}(t)]$. 
The numerical flux fluctuation term $\mathbf{Flux}\left(\bar{\mathbf{Q}}\right)$ implements the boundary conditions  \eqref{eq:BC_General2}, and  the interface conditions \eqref{eq:elastic_elastic_interface}, 
at the element faces, and it is defined by
{
\begin{align}\label{eq:num_flux_fluctuation}
\mathbf{Flux}\left(\bar{\mathbf{Q}}\right):= \sum_{\xi = q, r, s}{\mathbf{H}_\xi^{-1}}\left(\mathbf{e}_\xi(0)  \left[\mathbf{J}\sqrt{\xi_x^2 + \xi_y^2 + \xi_z^2}\mathbf{FL}(\bar{\mathbf{Q}}(t))\right]\Big|_{\xi = 0}  + \mathbf{e}_\xi(1)  \left[\mathbf{J}\sqrt{\xi_x^2 + \xi_y^2 + \xi_z^2}\mathbf{FR}(\bar{\mathbf{Q}}(t))\right]\Big|_{\xi = 1}  \right), 
\end{align}
}
with
\begin{equation}\label{eq:disc_scalar_product_flux}
\begin{split}
\bar{\mathbf{Q}} ^T {\mathbf{H}} \mathbf{Flux}\left(\bar{\mathbf{Q}}\right) 
& = \sum_{\psi = 1}^{P+1}\sum_{\theta = 1}^{P+1}\left[\left(J\sqrt{{\xi}_x^2 + {\xi}_y^2 + {\xi}_z^2}  \left(\mathbf{v}^T \mathbf{G} - \mathbf{T}^T \widetilde{\mathbf{G}}\right)\right)\Big|_{\xi = 0} \right]_{\psi\theta}h_{\psi}h_{\theta}\\
&+\sum_{\psi = 1}^{P+1}\sum_{\theta = 1}^{P+1}\left[\left(J\sqrt{{\xi}_x^2 + {\xi}_y^2 + {\xi}_z^2}\left(\mathbf{v}^T \mathbf{G} + \mathbf{T}^T\mathbf{G} \right)\right)\Big|_{\xi = 1} \right]_{\psi\theta}h_{\psi}h_{\theta}.
\end{split}
\end{equation}
The elemental degrees of freedom have been arranged row-wise as a single vector of length $9(P+1)^3$. Note the close similarity between the semi-discrete approximation \eqref{eq:gen_hyp_transformed_discrete} and the continuous analogue \eqref{eq:gen_hyp_transformed}. The discrete operator  $\grad_{D} = \left(\mathbf{D}_q, \mathbf{D}_r, \mathbf{D}_s\right)^T$ is also analogous to the continuous gradient operator $\grad = \left(\partial/\partial q, \partial/\partial r, \partial/\partial s\right)^T$. In  $\grad_{D}$ we have replaced the continuous derivative operators in $\grad$ with their discrete counterparts, 
\[
\frac{\partial}{\partial q} \to \mathbf{D}_q, \quad  \frac{\partial}{\partial r} \to \mathbf{D}_r, \quad \frac{\partial}{\partial s} \to \mathbf{D}_s,
\]
where  the spatial derivative operators, $\mathbf{D}_q, \mathbf{D}_r, \mathbf{D}_s$, are given in \eqref{eq:derivative_operator_3d}. 
\subsection{Numerical stability}
To prove the stability of the semi-discrete approximation \eqref{eq:gen_hyp_transformed_discrete}, we will  derive a semi--discrete energy estimate analogous to \eqref{eq:weak_energy}.  To begin, approximate the continuous energy  
in each element by the quadrature rule \eqref{eq:quad_rule_3D}, having
{
\small
\begin{align}\label{eq:disc_energy}
\mathcal{E}(t) := \frac{1}{2}\bar{\mathbf{Q}} ^T {\mathbf{H}} \widetilde{\mathbf{P}}^{-1}  \bar{\mathbf{Q}} =  \sum_{i=1}^{P+1}\sum_{j=1}^{P+1}\sum_{k=1}^{P+1}\frac{1}{2} [\bar{\mathbf{Q}}^T\left(q_i, r_j , s_k\right) \widetilde{\mathbf{P}}^{-1}\left(q_i, r_j , s_k\right)  \bar{\mathbf{Q}}\left(q_i, r_j , s_k\right)] h_ih_j h_k.
\end{align}
}
We also approximate the surface integrals in the boundary and interface terms
\eqref{eq:interface_term}, \eqref{eq:energy_estimate_BTs_hat} and \eqref{eq:weak_fluctuate}. 
We will now prove the discrete anti-symmetric property \eqref{eq:anti_symmetry} and a discrete analogue of Lemma \ref{Lem:Anti-symmetry_00}.
We have
\begin{lemma}\label{Lem:Discrete_Anti-symmetry}
Consider the semi-discrete approximations  of the transformed equation of motion \eqref{eq:gen_hyp_transformed_discrete}, in curvilinear coordinates, with the flux terms and non-conservative product terms given by  \eqref{eq:flux_ncp}. For a polynomial approximation of degree $P$, if the quadrature rule is exact  for all polynomial integrand $f(\xi)$ of degree $\le 2P-1$,  $\sum_{m = 1}^{P+1} f(\xi_m)h_m = \int_{0}^{1}f(\xi) d\xi$, then the corresponding spatial discrete  operators satisfy the discrete anti-symmetric property
\begin{align}\label{semi_discrete_Anti-symmetry_0}
    \bar{\mathbf{Q}}^T \mathbf{B}_\xi \left(\grad_{D}\bar{\mathbf{Q}}\right) - \left(\mathbf{D}_{\xi} \bar{\mathbf{Q}}\right)^T  \mathbf{F}_\xi \left(\bar{\mathbf{Q}} \right) = 0,
\end{align}
and
{
\begin{equation}\label{semi_discrete_sbp_vector}
\begin{split}
&\bar{\mathbf{Q}} ^T {\mathbf{H}}  \grad_{D}\bullet {\mathbf{F} \left(\bar{\mathbf{Q}} \right)} =  - \sum_{\xi = q, r, s}\left(\mathbf{D}_{\xi} \bar{\mathbf{Q}}\right)^T {\mathbf{H}} \mathbf{F}_\xi \left(\bar{\mathbf{Q}} \right) \\
&+  \sum_{\xi = q, r, s}\sum_{\psi = 1}^{P+1}\sum_{\theta = 1}^{P+1}\left[\left(J\sqrt{\xi_x^2 + \xi_y^2 + \xi_z^2}\mathbf{v}^T\mathbf{T} \right)_{\xi = 1} - \left(J\sqrt{\xi_x^2 + \xi_y^2 + \xi_z^2}\mathbf{v}^T\mathbf{T} \right)_{\xi = 0}\right]_{\psi\theta}h_{\psi}h_{\theta}.
\end{split}
\end{equation}
}
\end{lemma}
\begin{proof}
Equation \eqref{semi_discrete_Anti-symmetry_0} follows directly from  our choice of the anti-symmetric form \eqref{eq:gen_hyp_antisymmetry}.
Now consider 
  \begin{align}\label{eq:discrete_flux_integration_by_parts_0}
 \bar{\mathbf{Q}} ^T {\mathbf{H}}  \grad_{D}\bullet {\mathbf{F} \left(\bar{\mathbf{Q}} \right)} &=  \sum_{\xi = q, r, s}\bar{\mathbf{Q}} ^T {\mathbf{H}}   \mathbf{D}_{\xi} \mathbf{F}_\xi \left(\bar{\mathbf{Q}}  \right), 
\end{align}
 and use the discrete integrate-by-parts principle, \eqref{eq:derivative_operator_1d} and \eqref{eq:sbp_property_a}, we have
   \begin{equation}\label{eq:discrete_flux_integration_by_parts_1}
   \begin{split}
\bar{\mathbf{Q}} ^T {\mathbf{H}}  \grad_{D}\bullet {\mathbf{F} \left(\bar{\mathbf{Q}} \right)} =  &\sum_{\xi = q, r, s}\left(-\left(\mathbf{D}_{\xi}\bar{\mathbf{Q}} \right)^T {\mathbf{H}}    \mathbf{F}_\xi \left(\bar{\mathbf{Q}}  \right)
 +  \bar{\mathbf{Q}}^T{\mathbf{H}}{\mathbf{H}}_{\xi}^{-1}\left(\mathbf{B}_{\xi}\left(1, 1\right) - \mathbf{B}_{\xi}\left(0, 0\right)\right)\mathbf{F}_\xi \left( \bar{\mathbf{Q}} \right) \right).
 \end{split}
\end{equation}
Using the fact that
{
\begin{displaymath}
\begin{split}
\bar{\mathbf{Q}}^T{\mathbf{H}}{\mathbf{H}}_{\xi}^{-1}\left(\mathbf{B}_{{\xi}}\left(1, 1\right) - \mathbf{B}_{{\xi}}\left(0, 0\right)\right)\mathbf{F}_{\xi} \left( \bar{\mathbf{Q}} \right) 
& = \sum_{\psi = 1}^{P+1}\sum_{\theta = 1}^{P+1}\left[\left(J\sqrt{{\xi}_x^2 + {\xi}_y^2 + {\xi}_z^2}\mathbf{v}^T\mathbf{T} \right)_{{\xi} = 1}\right]_{\psi\theta}h_{\psi}h_{\theta}\\
&-\sum_{\psi = 1}^{P+1}\sum_{\theta = 1}^{P+1}\left[\left(J\sqrt{{\xi}_x^2 + {\xi}_y^2 + {\xi}_z^2}\mathbf{v}^T\mathbf{T} \right)_{{\xi} = 0}\right]_{\psi\theta}h_{\psi}h_{\theta}
\end{split}
\end{displaymath}
}
 completes the proof.
\end{proof}
A direct consequence of the discrete anti-symmetric  form \eqref{semi_discrete_Anti-symmetry_0} is the  discrete identity
{
\small
\begin{align}\label{semi_discrete_Anti-symmetry}
  \bar{\mathbf{Q}}^T{\mathbf{H}} \mathbf{B}_\xi \left(\grad_{D}\bar{\mathbf{Q}}\right) - \left(\mathbf{D}_{\xi} \bar{\mathbf{Q}}\right)^T {\mathbf{H}} \mathbf{F}_\xi \left(\bar{\mathbf{Q}} \right) = \sum_{i=1}^{P+1}\sum_{j=1}^{P+1}\sum_{k=1}^{P+1} \left( \bar{\mathbf{Q}}^T \mathbf{B}_\xi \left(\grad_{D}\bar{\mathbf{Q}}\right) - \left(\mathbf{D}_{\xi} \bar{\mathbf{Q}}\right)^T  \mathbf{F}_\xi \left(\bar{\mathbf{Q}} \right)\right)_{ijk} h_{i}h_{j}h_{k} = 0.
\end{align}
}
   
   For the two elements model we introduce the external boundary terms 
 \begin{align*}
 \mathcal{BT}_s\left(\widehat{v}_{\eta}^{\pm}, \widehat{T}_{\eta}^{\pm}\right)  
 & =  \sum_{\xi =  r, s}\sum_{\psi = 1}^{P+1}\sum_{\theta = 1}^{P+1}\left(\left(\sqrt{\xi_x^2 + \xi_y^2 + \xi_z^2}J\sum_{\eta = l, m, n} \widehat{T}_\eta^{\pm} \widehat{v}_\eta^{\pm}  \right)\Big|_{\xi = 1}\right)_{\psi \theta} h_{\psi} h_{\theta}\\
  &- \sum_{\xi = r, s}\sum_{\psi = 1}^{P+1}\sum_{\theta = 1}^{P+1}\left(\left(\sqrt{\xi_x^2 + \xi_y^2 + \xi_z^2}J\sum_{\eta = l, m, n}\widehat{T}_\eta^{\pm} \widehat{v}_\eta^{\pm}  \right)\Big|_{\xi = 0}\right)_{\psi \theta} h_{\psi} h_{\theta}\\
&\pm\sum_{\psi = 1}^{P+1}\sum_{\theta = 1}^{P+1}\left(\left(\sqrt{q_x^2 + q_y^2 + q_z^2}J\sum_{\eta = l, m, n}\widehat{T}_\eta^{\pm} \widehat{v}_\eta^{\pm}  \right)\Big|_{q = 0}\right)_{\psi \theta} h_{\psi} h_{\theta} ,
 \end{align*}
  the interface term
  \begin{align*}
  \mathcal{IT}_s\left(\widehat{v}^{\pm}, \widehat{T}^{\pm} \right) = - \sum_{\psi = 1}^{P+1}\sum_{\theta = 1}^{P+1}\left(\sqrt{q_x^2 + q_y^2 + q_z^2} J \sum_{\eta = l,m,n}{\widehat{T}}_\eta  \lJump{{\widehat{v}}_\eta \rJump} \right)_{\psi \theta} h_{\psi} h_{\theta} \equiv 0,
\end{align*}
      and the fluctuation term
    \begin{align*}
  \mathcal{F}_{luc}\left({G},  Z\right) = &- \sum_{\xi = q, r, s}\sum_{\psi = 1}^{P+1}\sum_{\theta = 1}^{P+1}\left(\left(\sqrt{\xi_x^2 + \xi_y^2 + \xi_z^2}J\sum_{\eta = l, m, n}\frac{1}{Z_{\eta}}|{G}_\eta |^2 \right)\Big|_{\xi = 0}\right)_{\psi \theta} h_{\psi} h_{\theta}\\
  &- \sum_{\xi = q, r, s}\sum_{\psi = 1}^{P+1}\sum_{\theta = 1}^{P+1}\left(\left(\sqrt{\xi_x^2 + \xi_y^2 + \xi_z^2}J\sum_{\eta = l, m, n}\frac{1}{Z_{\eta}}|{G}_\eta |^2 \right)\Big|_{\xi = 1}\right)_{\psi \theta} h_{\psi} h_{\theta}.
  \end{align*}
   \begin{remark}
   As in the continuous setting, it is noteworthy that by construction  the interface terms vanish identically $ \mathcal{IT}_s(\widehat{v}^{\pm}, \widehat{T}^{\pm}) \equiv 0$, by inspection the fluctuation terms are negative definite $ \mathcal{F}_{luc}\left({G},  Z\right) < 0$, and will vanish in the limit of mesh resolution. By  Lemma \ref{Lem:BTs_hat} the boundary terms are negative semi-definite $\mathcal{BT}_s\left(\widehat{v}_{\eta}^{-}, \widehat{T}_{\eta}^{-}\right) \le 0$.
\end{remark}
   We can now prove the numerical stability of the  semi-discrete approximation \eqref{eq:gen_hyp_transformed_discrete}.
\begin{theorem}\label{theo:semi_discrete_statbility}
Consider the semi-discrete DG approximation \eqref{eq:gen_hyp_transformed_discrete}. If the quadrature rule  $\sum_{m = 1}^{P+1} f(\xi_m)h_m = \int_{0}^{1}f(\xi) d\xi$  is exact  for all polynomial integrand $f(\xi)$ of degree $\le 2P-1$, then the numerical solution 
 satisfies the energy equation
   {
   \small
 \begin{equation}\label{eq:weak_energy_discrete}
\begin{split}
&\frac{d}{dt} \left(\mathcal{E}^{-}(t) + \mathcal{E}^{+}(t)\right)= \mathcal{IT}_s\left(\widehat{v}^{\pm},  \widehat{T}^{\pm}\right)  +   \mathcal{BT}_s\left(\widehat{v}^{-}, \widehat{T}^{-}\right)  +  \mathcal{BT}_s\left(\widehat{v}^{+}, \widehat{T}^{+}\right) + \mathcal{F}_{luc}\left({G}^{-},  Z^{-}\right)  + \mathcal{F}_{luc}\left({G}^{+},  Z^{+}\right)\le 0.
\end{split}
\end{equation}
}
\end{theorem}
\begin{proof}
We will use the energy method. From the left, multiply  \eqref{eq:gen_hyp_transformed_discrete} by $\bar{\mathbf{Q}} ^T {\mathbf{H}}$
  \begin{align}\label{eq:gen_hyp_transformed_discrete_proof_1}
\bar{\mathbf{Q}} ^T {\mathbf{H}} \widetilde{\mathbf{P}}^{-1} \frac{d }{d t} \bar{\mathbf{Q}} =\bar{\mathbf{Q}} ^T {\mathbf{H}}  \grad_{D}\bullet {\mathbf{F} \left(\bar{\mathbf{Q}} \right)} + \sum_{\xi= q, r, s}\bar{\mathbf{Q}} ^T {\mathbf{H}} \mathbf{B}_\xi \left(\grad_{D}\bar{\mathbf{Q}}\right) -\bar{\mathbf{Q}} ^T {\mathbf{H}} \mathbf{Flux}\left(\bar{\mathbf{Q}}\right).
\end{align}
On the left hand side  of \eqref{eq:gen_hyp_transformed_discrete_proof_1} we recognize time derivative of the semi-discrete energy $\mathcal{E}(t) $ defined in \eqref{eq:disc_energy}.
With  $\sum_{m = 1}^{P+1} f(q_m)h_m = \int_{0}^{1}f(q) dq$, then the summation-by-parts principle, \eqref{eq:sbp_property_a}, \eqref{eq:derivative_operator_1d} holds.
  On the right hand side of \eqref{eq:gen_hyp_transformed_discrete_proof_1}, we use Lemma \ref{Lem:Discrete_Anti-symmetry}, that is replace the conservative flux term with \eqref{semi_discrete_sbp_vector}, we have
{
 \begin{equation}\label{eq:gen_hyp_transformed_discrete_proof_2}
 \begin{split}
\frac{d }{d t}\mathcal{E}(t)   &= \sum_{\xi= q, r, s}\left[\bar{\mathbf{Q}}^T{\mathbf{H}} \mathbf{B}_\xi \left(\grad_{D}\bar{\mathbf{Q}}\right) - \left(\mathbf{D}_{\xi} \bar{\mathbf{Q}}\right)^T {\mathbf{H}} \mathbf{F}_\xi \left(\bar{\mathbf{Q}} \right)\right] \\
&+  \sum_{\xi = q, r, s}\sum_{\psi = 1}^{P+1}\sum_{\theta = 1}^{P+1}\left[\left(J\sqrt{\xi_x^2 + \xi_y^2 + \xi_z^2}\mathbf{v}^T\mathbf{T} \right)_{\xi = 1} - \left(J\sqrt{\xi_x^2 + \xi_y^2 + \xi_z^2}\mathbf{v}^T\mathbf{T} \right)_{\xi = 0}\right]_{\psi\theta}h_{\psi}h_{\theta}
 -\bar{\mathbf{Q}} ^T {\mathbf{H}} \mathbf{Flux}\left(\bar{\mathbf{Q}}\right).
 \end{split}
\end{equation}
}
Using the fact \eqref{semi_discrete_Anti-symmetry}, the volume terms vanish, remaining only the surface terms.
Simplifying the flux term $\bar{\mathbf{Q}} ^T {\mathbf{H}} \mathbf{Flux}\left(\bar{\mathbf{Q}}\right)$ with \eqref{eq:disc_scalar_product_flux} having 
 \begin{equation}\label{eq:transf_elemental_weak_form_velocity_semi_disc}
 \begin{split}
\frac{d }{d t}\mathcal{E}(t)  = 
- & \sum_{\xi = q, r, s}\sum_{\psi = 1}^{P+1}\sum_{\theta = 1}^{P+1}\left(\left(J\sqrt{\xi_x^2 + \xi_y^2 + \xi_z^2}\left(\mathbf{v}^T \mathbf{G} - \mathbf{T}^T \widetilde{\mathbf{G}} + \mathbf{v}^T\mathbf{T}\right)\right)_{\xi = 0} \right)_{\psi\theta}h_{\psi}h_{\theta} \\
-&  \sum_{\xi = q, r, s}\sum_{\psi = 1}^{P+1}\sum_{\theta = 1}^{P+1}\left(\left(J\sqrt{\xi_x^2 + \xi_y^2 + \xi_z^2}\left(\mathbf{v}^T \mathbf{G} + \mathbf{T}^T \widetilde{\mathbf{G}} - \mathbf{v}^T\mathbf{T}\right)\right)_{\xi = 1}\right)_{\psi\theta}h_{\psi}h_{\theta} .
 \end{split}
\end{equation}
Using \eqref{eq:identity_pen} in the right hand side of \eqref{eq:transf_elemental_weak_form_velocity_semi_disc} gives
{
\begin{equation}\label{eq:transf_elemental_weak_form_velocity_semi_disc0}
\begin{split}
\frac{d }{d t}\mathcal{E}(t)   = 
- & \sum_{\xi = q, r, s}\sum_{\psi = 1}^{P+1}\sum_{\theta = 1}^{P+1}\left(\left(J\sqrt{\xi_x^2 + \xi_y^2 + \xi_z^2}\sum_{\eta = l,m,n} \left(\frac{1}{Z_\eta  }|G_\eta |^2  + \widehat{T}_\eta \widehat{v}_\eta \right)\right)_{\xi = 0} \right)_{\psi\theta}h_{\psi}h_{\theta} \\
- & \sum_{\xi = q, r, s}\sum_{\psi = 1}^{P+1}\sum_{\theta = 1}^{P+1}\left(\left(J\sqrt{\xi_x^2 + \xi_y^2 + \xi_z^2}\sum_{\eta = l,m,n} \left(\frac{1}{Z_\eta  }|G_\eta |^2  - \widehat{T}_\eta \widehat{v}_\eta \right)\right)_{\xi = 1}\right)_{\psi\theta}h_{\psi}h_{\theta}.
\end{split}
\end{equation}
}
Collecting contributions from both elements, and using  the identities \eqref{eq:identity_bc} and \eqref{eq:identity} give the energy equation \eqref{eq:weak_energy_discrete}.
\end{proof}

Note the similarities between Theorem \ref{theo:semi_discrete_statbility} and Theorem \ref{theo:weak_form_statbility}. The only differences are that we have replaced the solution $\mathbf{Q}$ by its polynomial interpolant $\bar{\mathbf{Q}}$, and the integrals by  quadrature rules.

By \eqref{eq:weak_energy_discrete}, the semi-discrete DG approximation \eqref{eq:gen_hyp_transformed_discrete} is asymptotically stable. This means that the solutions can never grow in time.
We can now integrate  \eqref{eq:gen_hyp_transformed_discrete} in time using any suitable explicit ODE  scheme, such as high order explicit Runge-Kutta methods. In this paper we will integrate the \eqref{eq:gen_hyp_transformed_discrete}  in time using the ADER scheme \cite{Toro1999, DumbserKaeser2006}. The ADER time discretisation is summarized in \ref{sec:ADER}. 


\section{Implementation in the ExaHyPE engine: ExaSeis}
We implement the here presented curvilinear method in the \emph{ExaHyPE-Engine} \cite{ExaHyPE2019} which is publicly available at \url{www.exahype.org}.
The ExaHyPE-Engine solves linear and non-linear hyperbolic partial differential equations in first order formulation.
In the linear case, ExaHyPE uses the ADER-DG method as summarised in \ref{sec:ADER}.
To solve for any given hyperbolic PDE in ExaHyPE, users are required to implement all PDE specific terms in a generated C++ interface.
For the method we present in this paper, these are the conservative and non-conservative flux terms of Eqs. \eqref{eq:anti_symmetry} and \eqref{eq:conervative_flux_and_nonconervative_product} and the eigenvalues in Eq. \eqref{eq:anisotropic_wavespeed}.
While the engine provides the Rusanov flux as default numerical flux, we here replace it with the numerical flux introduced in \eqref{eq:num_flux_fluctuation}, to properly treat boundary conditions and material interfaces.
As initial conditions for the experiments in Sec. \ref{sec:numerical_experiments} we implement problem specific material parameters and moment-tensor point sources.
We generate curvilinear meshes based on a simple k-d-tree approach, 
which is publicly available as extracted library at \url{https://gitlab.lrz.de/ExaHyPE-Seismic/curvi}.

The elastic wave propagation model with curvilinear meshes is part of the ExaSeis application collection (\url{https://gitlab.lrz.de/ExaHyPE-Seismic}).
ExaSeis includes several approaches related to computational seismology, allowing to study topography scattering effects or multi-physics dynamic rupture models.
It also allows to use perfectly matching layers (PML) \cite{DuruRannabauerGabrielKreissBader2019}.
To allow reproducibility of all results of Sec. \ref{sec:numerical_experiments} we provide a repository with prepared scripts at \url{https://gitlab.lrz.de/ExaHyPE-Seismic/ExaSeis-Benchmarks}.
\section{Numerical experiments}\label{sec:numerical_experiments}
We will now present  numerical simulations in complex 3D isotropic and anisotropic elastic solids. The experiments are designed to verify  accuracy and numerical stability of the method,  for both, body waves and elastic surface waves, as well as for high frequency scattered waves. We will also demonstrate the potential of our method for applications including geometrically complex free surface topography. 

We will consider three benchmark problems with different levels of difficulty \cite{Kristekova_etal2006, Kristekova_etal2009, Achenbach1973, Petersson_etal2016, Favretto-Cristini_etal2011, Komatitsch_etal2011,DuruFungWilliams2020}.  The benchmark problems are i) the 3D layer over a homogeneous half-space (LOH1) problem, ii) 3D elastic surface waves in an anisotropic medium, and iii) a 3D regional strong topography contrast simulation including scattered high frequency waves and a geologically informed complex surface geometry (Zugspitze model). These benchmark problems are designed to quantify and assess the accuracy of simulation codes for seismic surface and interface waves, and effective resolution of high frequency wave modes generated by scattering from the complex non-planar topography.

The numerical experiments will be conducted in 3D elastic media. 
We demonstrate spectral convergence by performing simulations on moderately fine 3D meshes where we vary the degree of the polynomial approximations, $P =3, 5, 7$. The simulated solutions are compared to analytical and reference solutions of community defined benchmark problems \cite{Kristekova_etal2006, Kristekova_etal2009, Favretto-Cristini_etal2011, Komatitsch_etal2011}. We will perform error analysis and compute error parameters relevant to computational seismology.  The reader is also referred to \cite{DuruGabrielIgel2017}, where detailed grid convergence studies for the physics based flux for elastodynamics are performed in 1D and 2D. 

In many of the  simulations, we have used PML boundary conditions to efficiently absorb outgoing waves. The PML modelling error is set to $\approx1\%$, which is a moderately small error tolerance. Details of the discretization and implementation of the PML is reported in \cite{DuruRannabauerGabrielKreissBader2019}. 
We have used GLL nodes for most of the experiments presented here. However, it is important to note that equivalent results have been obtained using GL nodes.

 We introduce the effective sub-element resolution, grid size $h = \Delta{x}/(P+1)$, comparable to a finite difference grid size. We will use the global time-step
\begin{align}\label{eq:time_step}
\Delta{t} = \frac{\mathrm{CFL}}{d} \frac{ h_{\text{min}}}{  c_{\mathrm{max}}}, \quad h_{\text{min}} = \frac{\Delta_{\mathrm{min}}}{\left(P+1\right)},
\end{align}
where $P$ is the  degree of the polynomial approximation, $d = 3$ is the spatial dimension,  $ \mathrm{CFL}  = 0.9$ and 
{
\small
\begin{align*}
\Delta_{\mathrm{min}}= \min\left(\Delta{x}, \Delta{y}, \Delta{z}\right), \quad c_{\mathrm{max}} = \max\left(\sqrt{\sum_{\xi=x,y,z}q_{\xi}^2}{c}_{n}, \sqrt{\sum_{\xi=x,y,z}r_{\xi}^2}{c}_{n}, \sqrt{\sum_{\xi=x,y,z}s_{\xi}^2}{c}_{n}\right), \quad  c_{n} =  \sqrt{\sum_{\xi=x,y,z}\left(n_{\xi}c_{p\xi}\right)^2}. 
\end{align*}
}
Note that $\sum_{\xi=x,y,z}n_{\xi}^2 = 1$,
and in an isotropic medium  the effective normal p-wave speed  is $c_n = c_p$.

We will consider seismic sources defined by the singular moment tensor point source 
{
\begin{align}\label{eq:pointsource}
\mathbf{f}(x,y,z,t) &=   \mathbf{M} \delta_x(x-x_0)\delta_y(y-y_0)\delta_z(z-z_0)g(t), \quad 
\mathbf{M} 
=
\begin{pmatrix}
M_{xx} & M_{xy} &M_{xz} \\
M_{xy} & M_{yy} &M_{yz} \\
M_{xz} & M_{yz} &M_{zz} 
\end{pmatrix}.
\end{align}
}
Here, $ \mathbf{M} $ is the symmetric second order moment tensor, $ \delta_{\eta}(\eta)$ are the one dimensional Dirac delta function, $(x_0, y_0, z_0)$ is the source location, and $g(t)$ is the source time function. 
\subsection{Layer over a half-space (LOH1)}
We consider the 3D LOH1 benchmark problem, which has an analytical solution \cite{Kristekova_etal2006,Kristekova_etal2009}. 
The LOH1 benchmark problem 
consists of a planar free surface and an internal planar interface separating a thin homogeneous 
soft layer and hard half-space. The material properties in the medium are given by
\begin{align*}
&\text{soft upper crust}: \quad \rho = 2600~\ \text{kg/m}^3, \quad c_p = 4000~\ \text{m/s}, \quad c_s = 2000~\ \text{m/s}, \quad x \le 1~\ \text{km},\\
&\text{hard lower crust}: \quad \rho = 2700~\ \text{kg/m}^3, \quad c_p = 6000~\ \text{m/s}, \quad c_s = 3464~\ \text{m/s}, \quad x  >1~\ \text{km}.
\end{align*}
The benchmark considers homogeneous initial conditions for all fields,  and generates waves by  adding  the double-couple moment tensor point source at depth, $x = 2~$km below the free surface.
The source time function is given by 
\begin{align}\label{eq:source_time_function}
 g(t) = \frac{t}{T^{2}}e^{-{t}/{T}},  \quad T = 0.1 ~\ \text{s},
\end{align}
and the moment tensor $\mathbf{M}$ which is zero except for the shear components $M_{yz}=M_{zy}=M_0$, and $M_0 =10^{18}$~Nm is the moment magnitude.
We place 9 receivers  on the free-surface, at $x = 0$, where the solutions are sampled. Table \ref{tbl:receivers_hhs1} shows the positions of the receivers relative to the epicenter. Receiver 6 and Receiver 9, which are about
$10.39~$km away from the epicenter, are the two farthest receivers from the source.
\begin{table}[h!]
  \caption{Receiver positions of the 3D HHS and LOH1 problems relative to the epicenter.}
  \centering{
    \begin{tabular}{|c|c|c|c|c|c|c|c|c|c|}
      \hline
      Receiver &1&2&3&4&5&6&7&8&9  \\ \hline
      y[km]& 0     & 0     & 0       & 0.490  & 3.919 & 7.348 & 0.577 & 4.612 & 8.647 \\\hline
      z[km]& 0.693 & 5.542 & 10.392  & 0.490  & 3.919 & 7.348 & 0.384 & 3.075 & 5.764 \\ \hline
    \end{tabular}
  }
  \label{tbl:receivers_hhs1}
\end{table}

Note that the domain is unbounded at depth and in the tangential directions.
%
%
%
%
To perform numerical simulations we consider the bounded computational cube $(x,y,z) = [0, 16.333~\text{km}]\times[-2.287~\text{km}, 14.046~\text{km}]\times[-2.287~\text{km}, 14.046~\text{km}]$, and the source is located at $(x_0, y_0, z_0) = (2~\text{km}, 0, 0)$.  The truncated boundaries of the domain were surrounding by the PML to absorb outgoing waves. First we discretize the computational domain uniformly with 25 DG elements in each spatial direction. We consider degrees $P =3, 5,7$ polynomial approximations, and run the simulation until the final time $t = 9~$s. 
\begin{figure}[h!]
\begin{subfigure}
    \centering
\stackunder[5pt]{\includegraphics[draft=false,width=0.325\columnwidth]{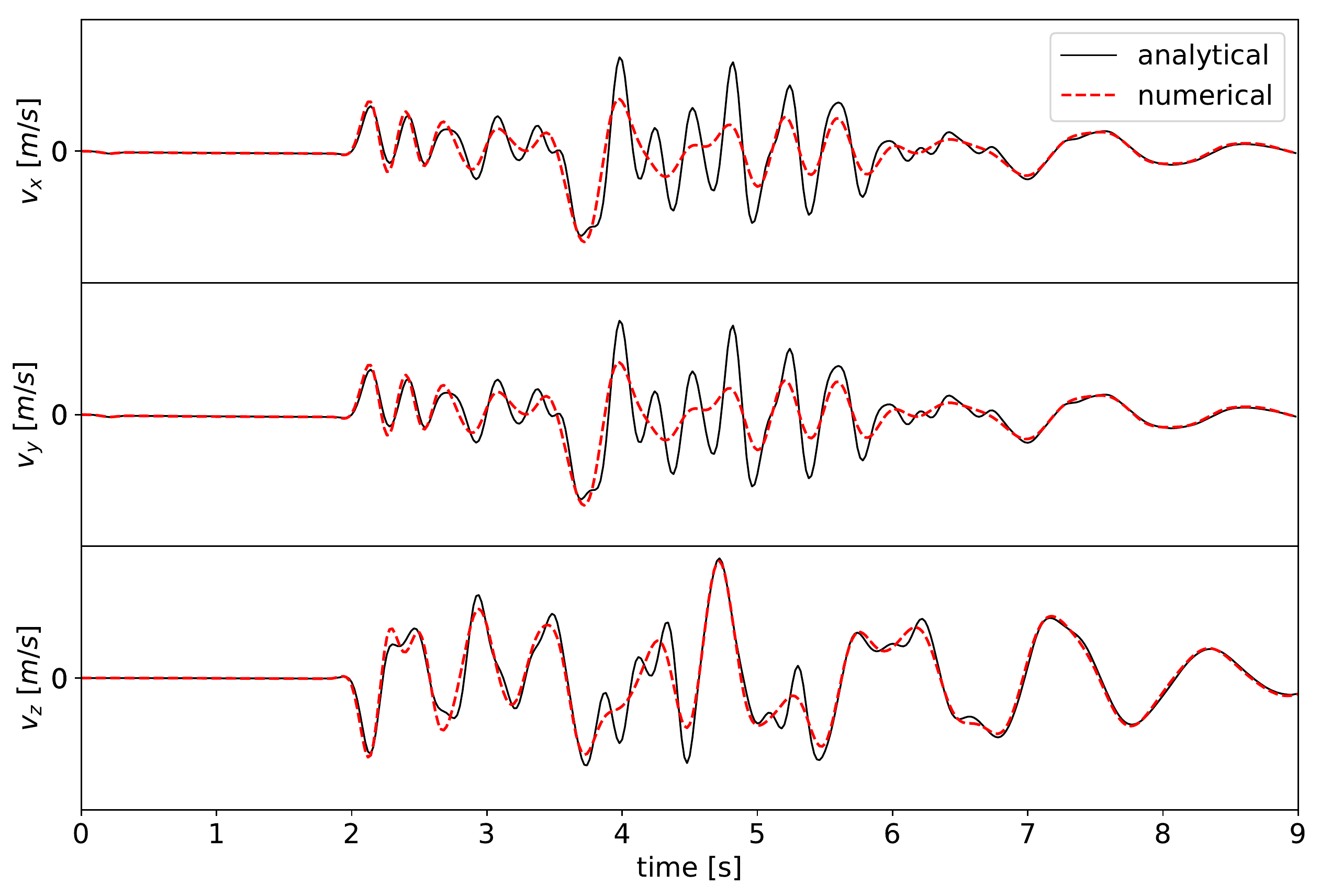}}{$P = 3$}%
\hspace{0.0cm}%
\stackunder[5pt]{\includegraphics[draft=false,width=0.325\columnwidth]{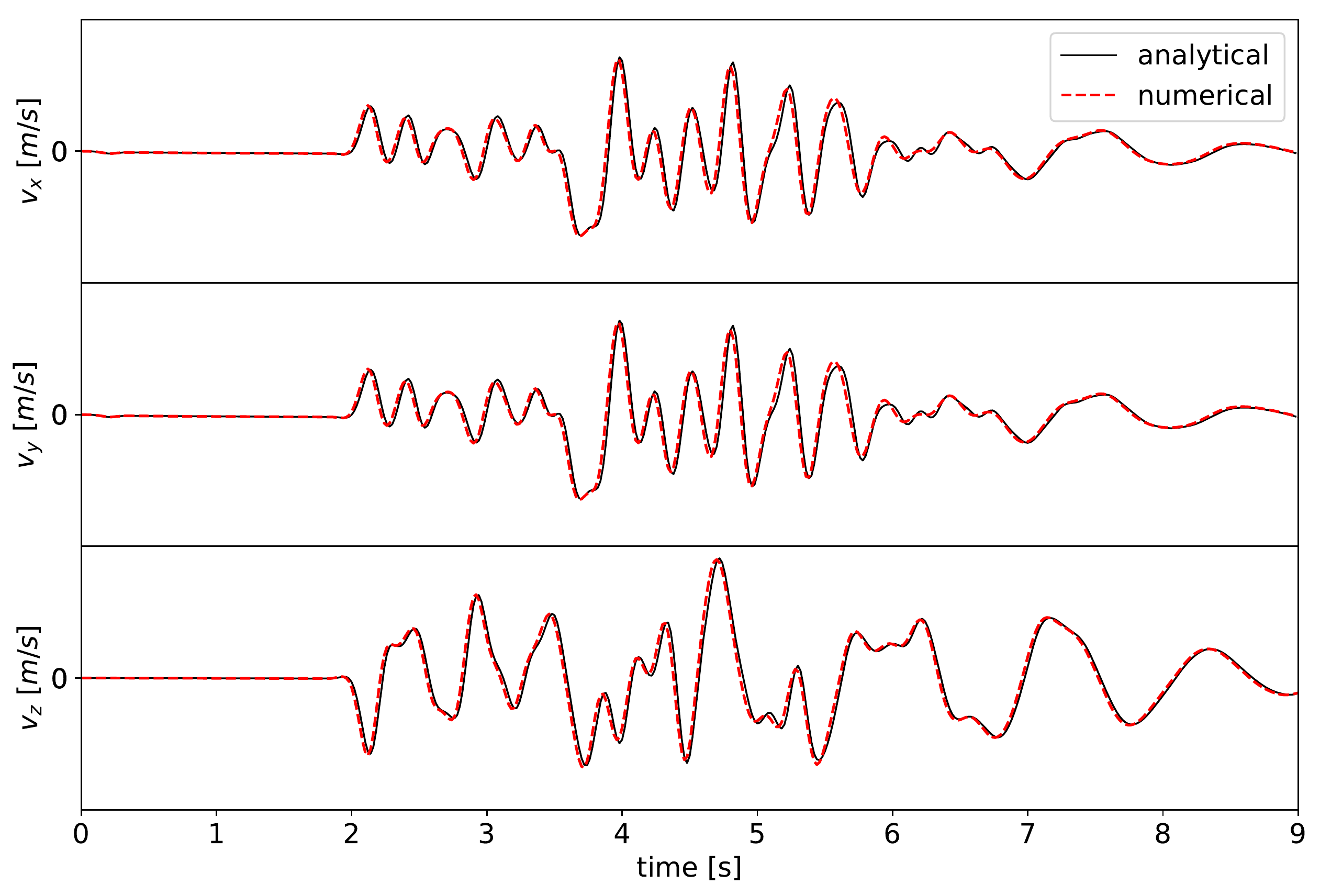}}{$P = 5$}%
\hspace{0.0cm}%
\stackunder[5pt]{\includegraphics[draft=false,width=0.325\columnwidth]{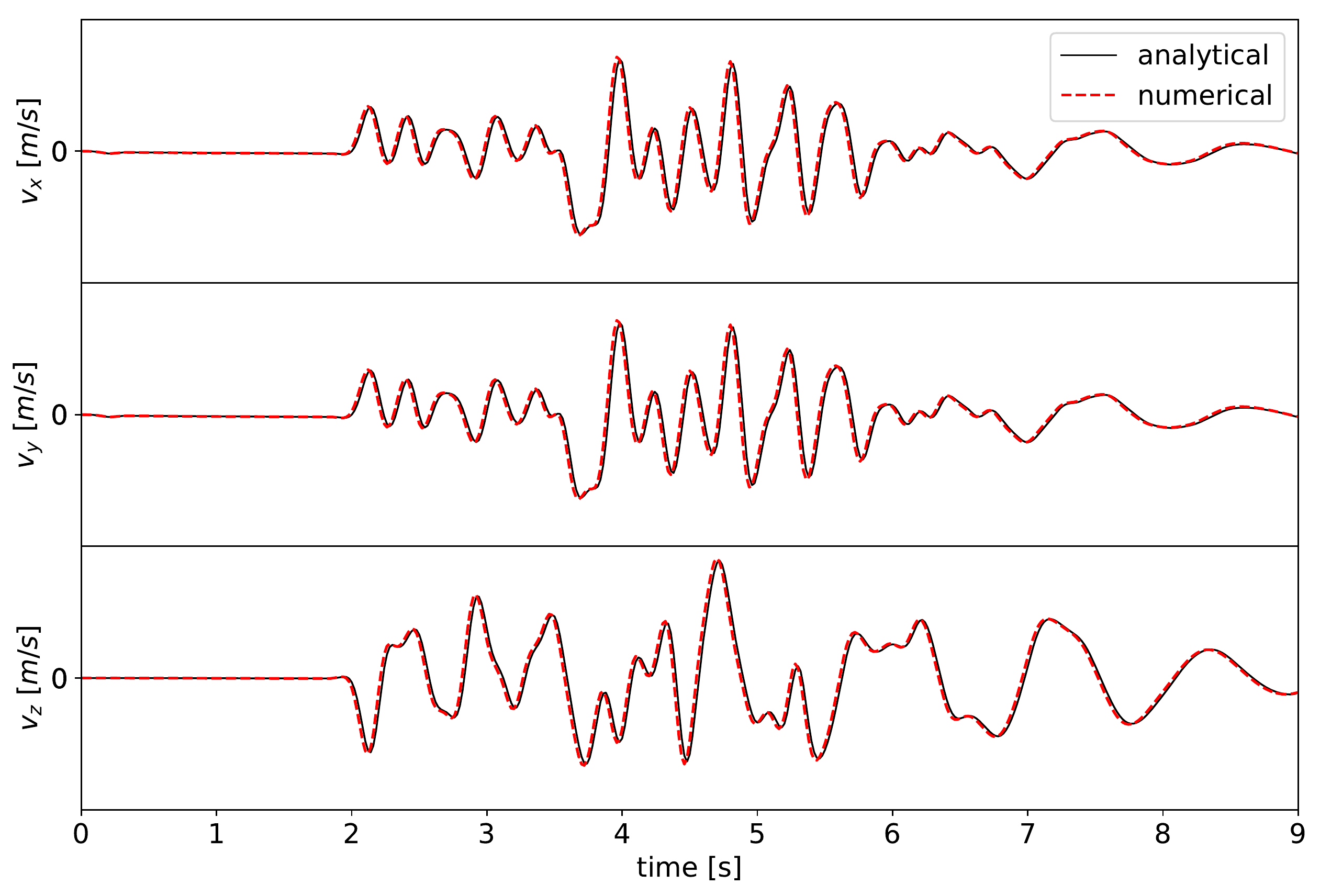}}{$P = 7$}%
     \end{subfigure}
    \caption{The 3D layer of half-space LOH1 benchmark problem. Comparing ExaSeis numerical solutions with the analytical solution for different polynomial degrees, $P = 3, 5, 7$, at receiver 6.}
    \label{fig:loh1_r6}
\end{figure}
\begin{figure}[h!]
\begin{subfigure}
    \centering
\stackunder[5pt]{\includegraphics[draft=false,width=0.325\columnwidth]{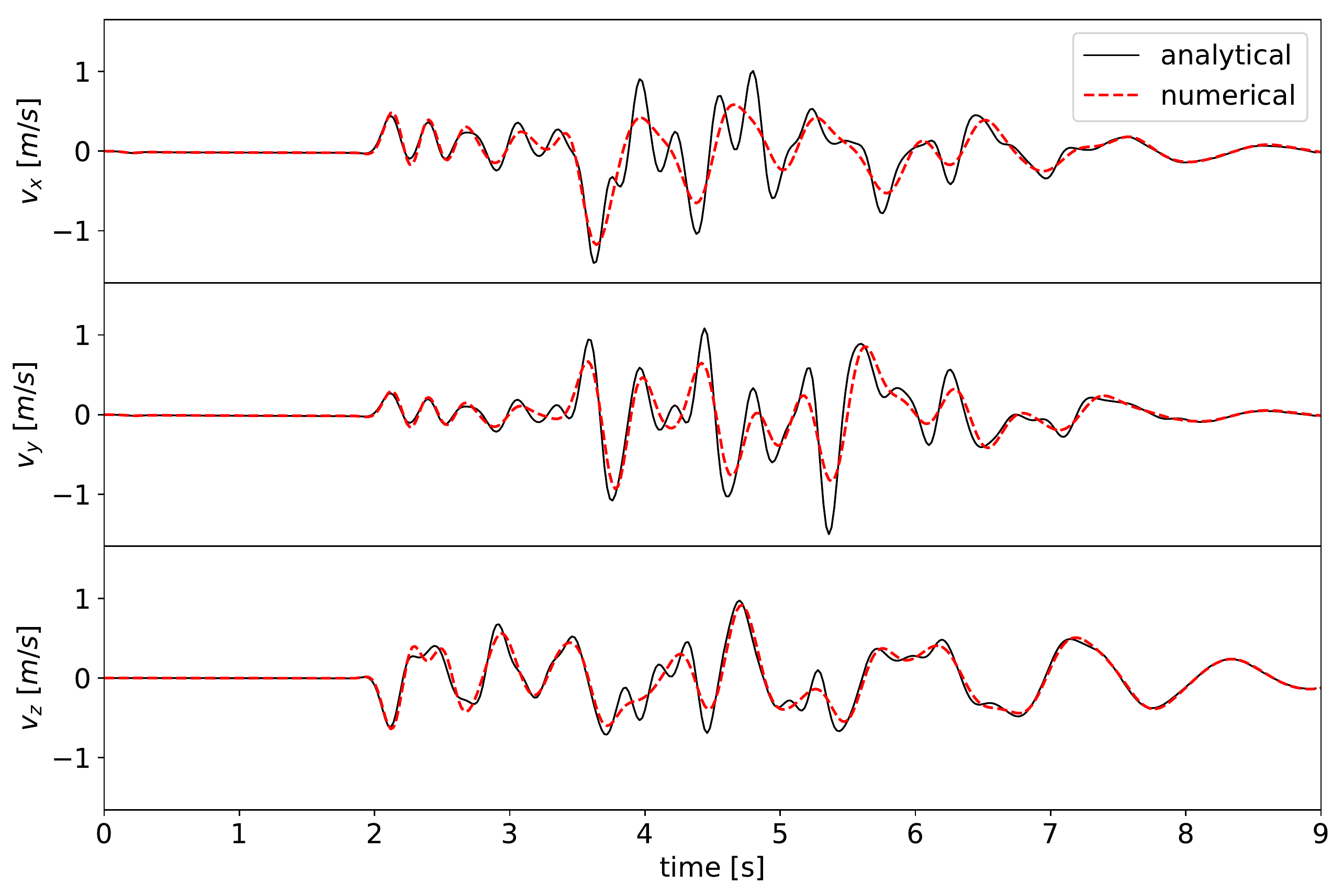}}{$P = 3$}%
\hspace{0.0cm}%
\stackunder[5pt]{\includegraphics[draft=false,width=0.325\columnwidth]{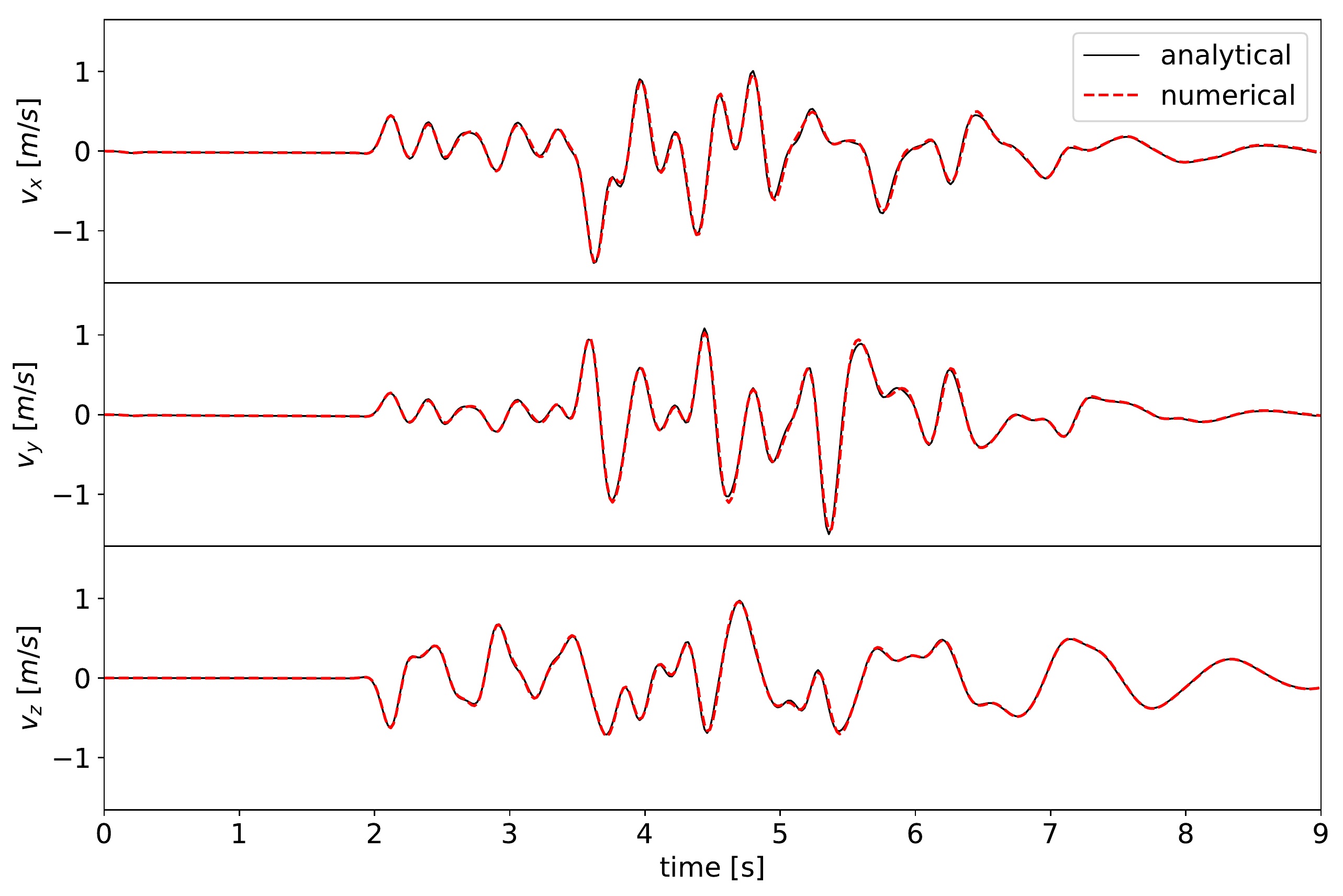}}{$P = 5$}%
\hspace{0.0cm}%
\stackunder[5pt]{\includegraphics[draft=false,width=0.325\columnwidth]{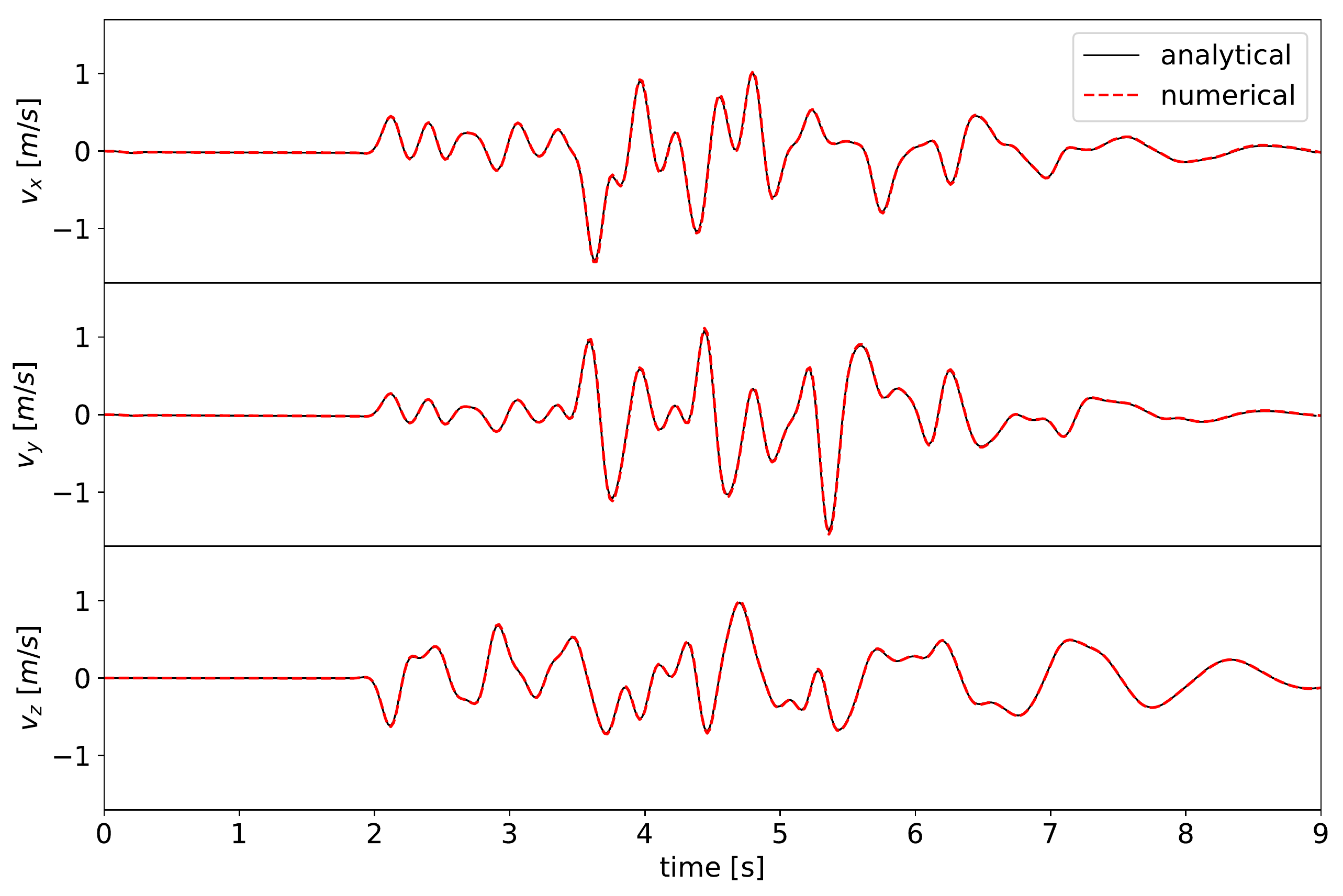}}{$P = 7$}%
     \end{subfigure}
    \caption{The 3D layer of half-space LOH1 benchmark problem. Comparing ExaSeis numerical solutions with the analytical solution for different polynomial degrees, $P = 3, 5, 7$, at receiver 9.}
    \label{fig:loh1_r9}
\end{figure}
The solutions are displayed in Figure \ref{fig:loh1_r6} for Receiver 6 and  in Figure \ref{fig:loh1_r9} for Receiver 9, and for $P = 3, 5, 7$. 
With increasing polynomial degree, $P = 3, 5, 7$,  the numerical  solutions converge spectrally to the analytical solution. For $P = 5, 7$, the numerical solutions match the exact solution excellently well, with about $\sim 1$\% relative error. 

\subsubsection{Error analysis}
We now analyse the numerical error for the LOH1 problem, and compute error parameters relevant to seismological applications.
To quantitatively assess the accuracy of the numerically simulated seismograms, we compare the numerical seismograms for the  LOH1 benchmark simulations with the exact reference solution using the time-frequency (TF) misfit criteria proposed by \cite{Kristekova_etal2006, Kristekova_etal2009}. We will briefly describe the technique here, and refer the reader to \cite{Kristekova_etal2006, Kristekova_etal2009} for a more elaborate discussion. The misfit criteria are based on the Time-Frequency Representation (TFR) of the seismogram, denoted as $W(t, f)$. A local TF envelope difference is defined
as:
\begin{displaymath}
\Delta{E}(t, f) = |W(t,f)|-|W_{ref}(t,f)|,
\end{displaymath}
and local TF phase difference is defined as
\begin{displaymath}
\Delta{P}(t, f) = |W_{ref}(t,f)| \frac{\mathrm{Arg}[W(t,f)] - \mathrm{Arg}[W_{ref}(t,f)]}{\pi},
\end{displaymath}
where $|W(t, f)|$ and $|Wref (t, f)|$ are the TFR of the numerical and reference data, respectively. To obtain a single-valued measure of the EM
or PM between these two seismograms,  we use the following normalized formulars
\begin{displaymath}
\mathrm{EM} =  \sqrt{ \frac{\sum_f\sum_t |\Delta{E}(t, f)|^2}{\sum_f\sum_t |{W}_{ref}(t, f)|^2}}, \quad \mathrm{PM} =   \sqrt{\frac{\sum_f\sum_t |\Delta{P}(t, f)|^2}{\sum_f\sum_t |{W}_{ref}(t, f)|^2}}.
\end{displaymath}
We use  the software package {\it Obspy} \cite{ObsPy2015} to post process the seismograms, by passing the seismograms through a band limited filter, and
  evaluate the accuracy level of the numerical seismograms  at all stations.
\begin{table}[h!]
  \caption{Accuracy levels}
  \centering{
    \begin{tabular}{|c|c|c|c|}
      \hline
      Accuracy level & Envelop Misfit & Phase Misfit  \\ \hline
      A & $\le 5\%$     & $\le 5\%$      \\ \hline
      B & $\le 10\%$     & $\le 10\% $     \\ \hline
      C & $\le 20\%$     & $ \le 20\%$      \\ \hline
    \end{tabular}
  }
  \label{tbl:loh1_misfit}
\end{table}
 Table \ref{tbl:loh1_misfit} classifies the quality of the numerical seismograms, where $A$, with EM, PM $\le 5\%$, indicates the highest quality numerical seismograms  \cite{Kristekova_etal2006, Kristekova_etal2009}.

\begin{figure}[h!]
\begin{subfigure}
    \centering
\stackunder[5pt]{\includegraphics[width=0.325\textwidth]{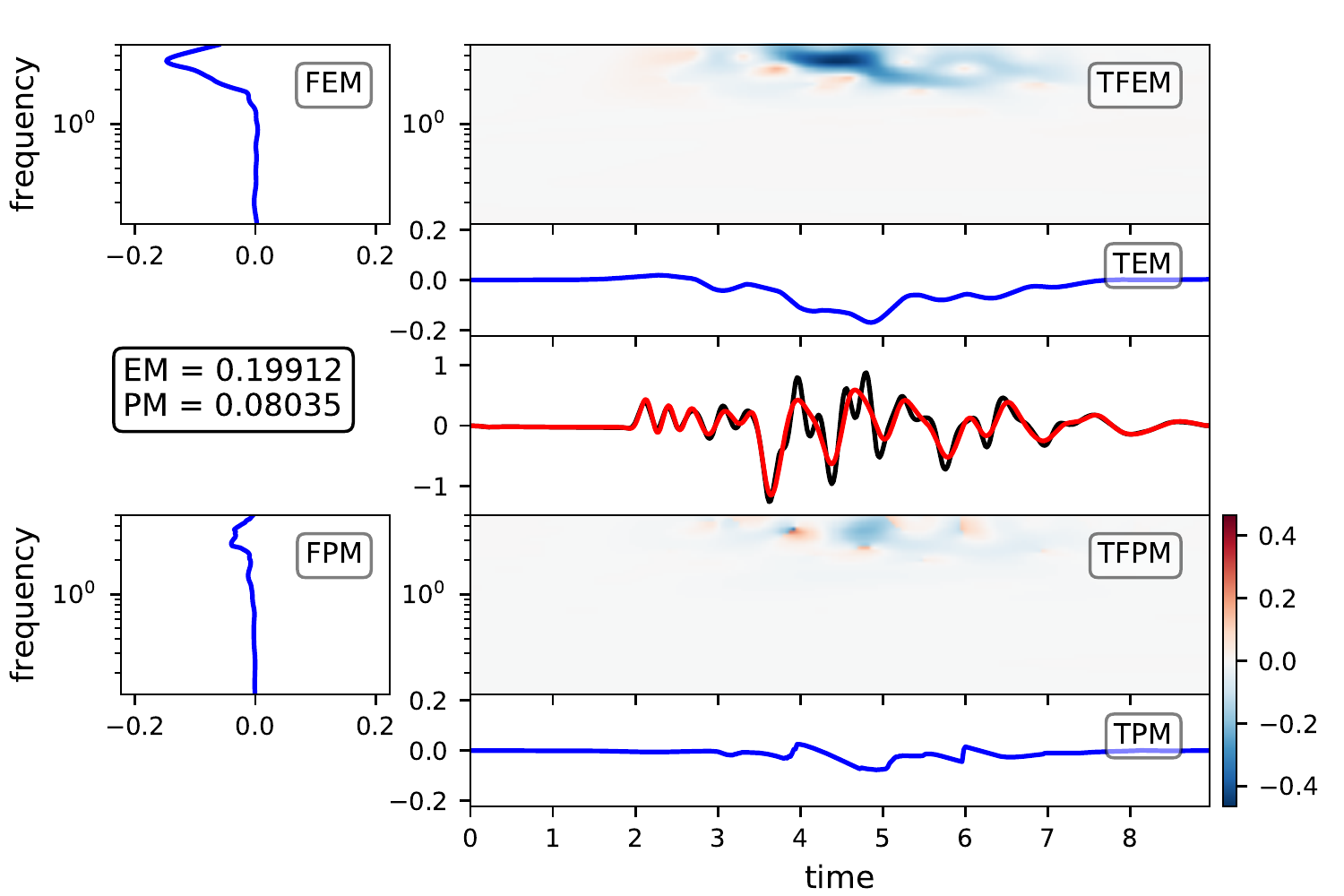}}{$P = 3$}%
\hspace{0.0cm}%
\stackunder[5pt]{\includegraphics[width=0.325\textwidth]{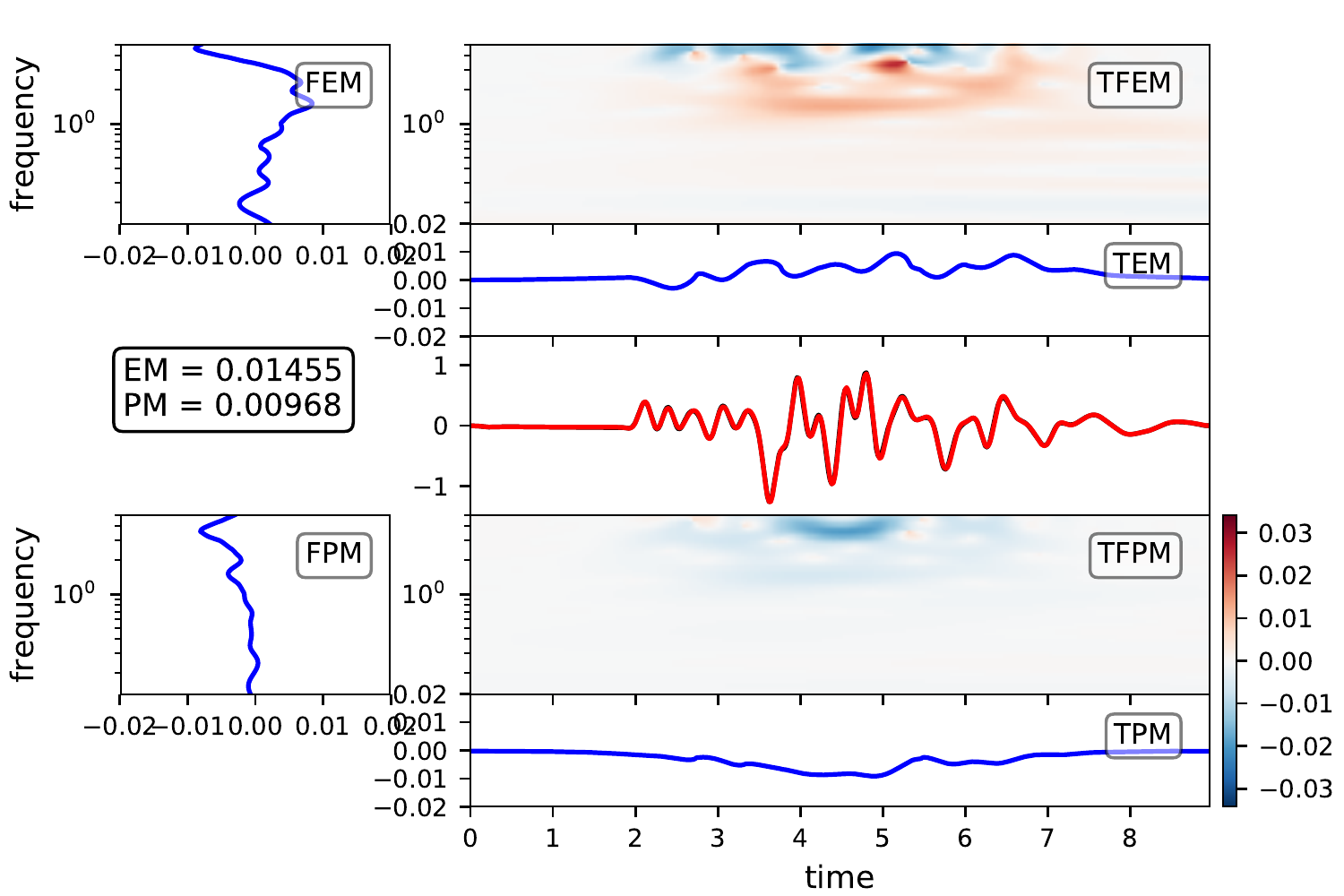}}{$P = 5$}%
\hspace{0.0cm}%
\stackunder[5pt]{\includegraphics[width=0.325\textwidth]{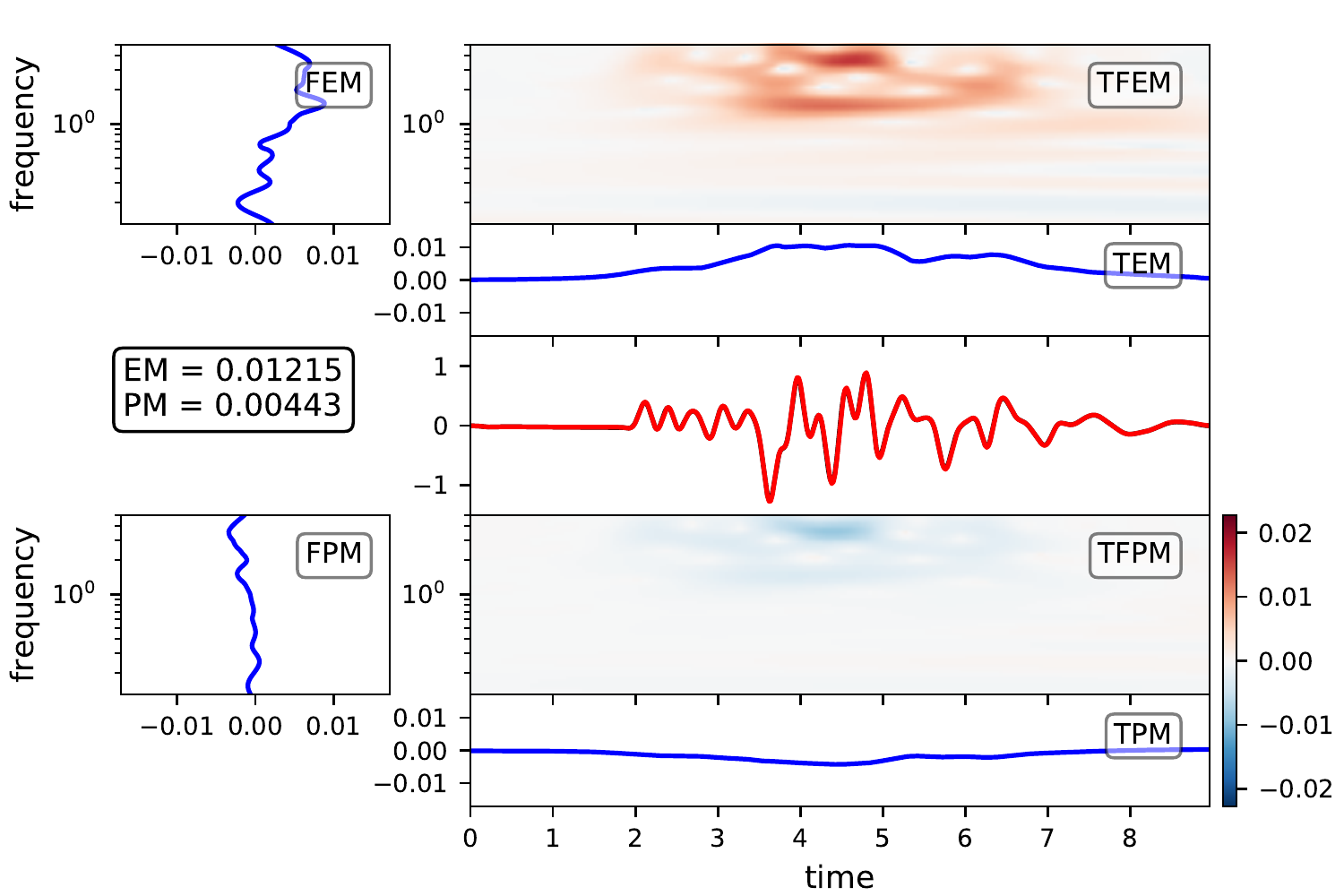}}{$P =7$}%
\hspace{0.0cm}%
     \end{subfigure}
    \caption{The 3D layer of half-space LOH1 benchmark problem. Time-frequency misfit analysis of the ExaSeis numerical and the analytical solutions for $v_x[\mathrm{m/s}]$ at Receiver 9 (0.13 -- 5 Hz).}
     \label{fig:loh1_tf_misfit_r9_vx}
\end{figure}

The results of the analysis, for Receiver 6 with $P = 3, 5, 7$, are displayed in Figure \ref{fig:loh1_tf_misfit_r9_vx}.
\begin{table}[h!]
  \caption{Accuracy for  polynomial degree $P=3,5,7$.}
  \centering{
    \begin{tabular}{|c|c|c|c|c|}
      \hline
     polynomial degree & Envelop Misfit & Phase Misfit & Accuracy level \\ \hline
      $P=3 $ & $\le 20.6\%$     & $\le 8\%$ & C      \\ \hline
      $P=5 $ & $\le 1.8\%$     & $\le 1 \% $  & A     \\ \hline
      $P= 7$ & $\le 1.5\%$     & $ \le 0.6\%$ & A  \\ \hline
    \end{tabular}
  }
  \label{tbl:loh1_misfit_P}
\end{table}
As shown in Figure \ref{fig:loh1_tf_misfit_r9_vx} and Table \ref{tbl:loh1_misfit_P}, at this resolution the seismograms for polynomial approximations of degree $ P = 5, 7$, belong to the A-class, are of the highest quality. We can also see that the seismograms for polynomial approximations of degree $ P = 3$ belong to the C-class. For $ P = 3$, the accuracy can be improved by increasing the mesh resolution (h-refinement).
\subsubsection{Non-conforming adaptive mesh}
\begin{figure}[ht!]%
	\begin{center}
	\includegraphics[draft=false,width=0.35\columnwidth]{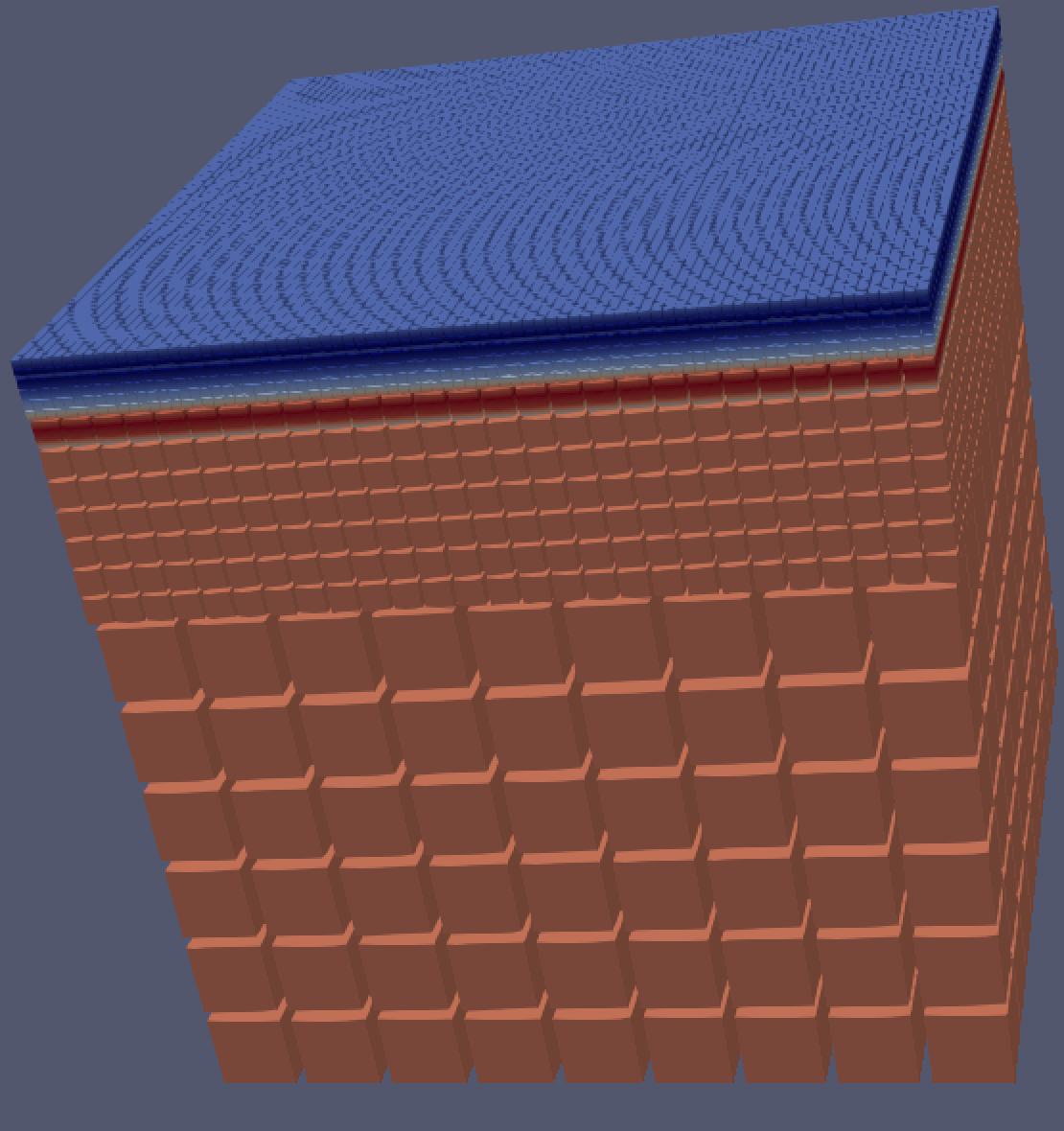} 
	\includegraphics[draft=false,width=0.6\columnwidth]{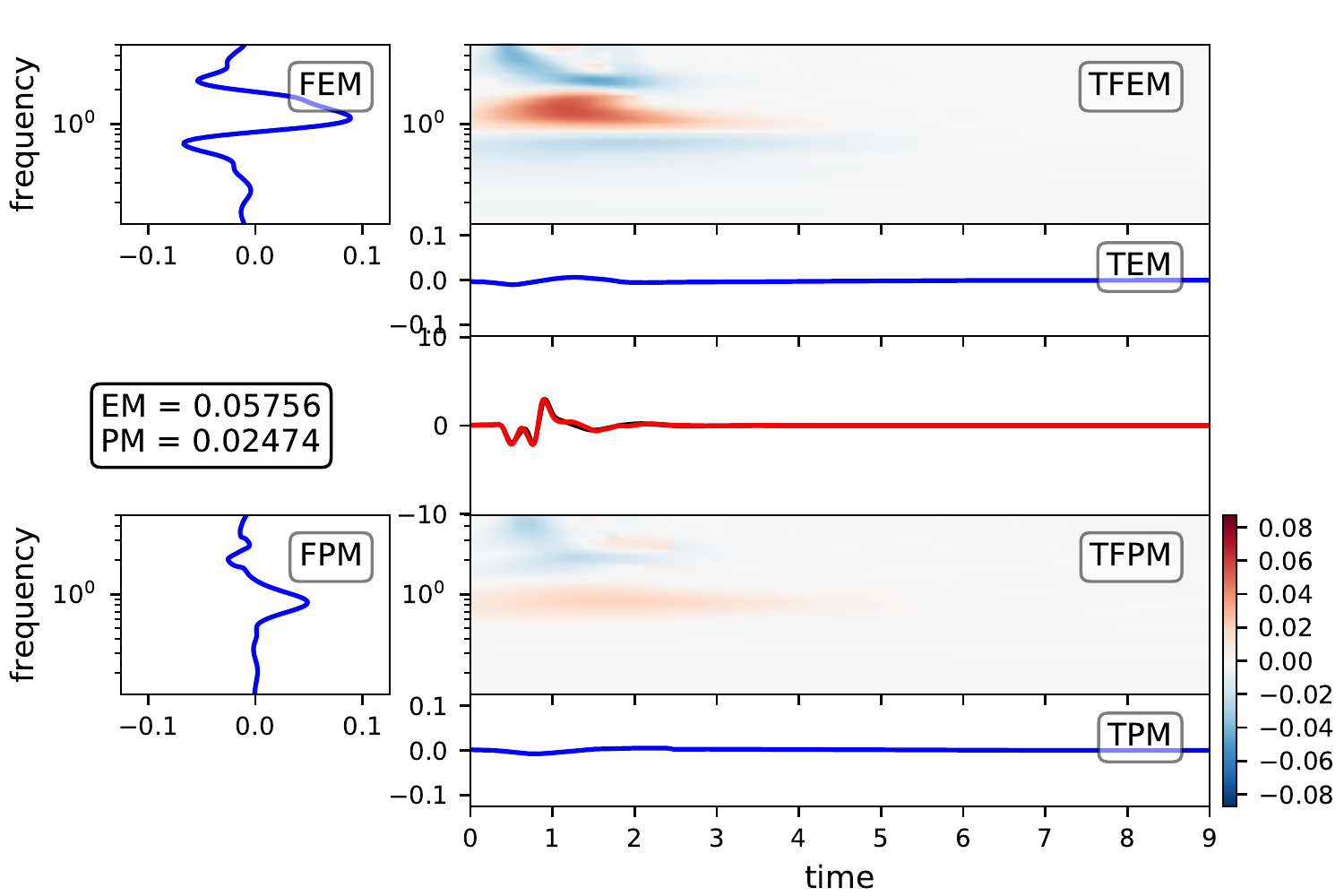} 
	\end{center}
	\caption{The 3D layer of half-space LOH1 benchmark problem. Time and frequency error misfit of $v_z[\mathrm{m/s}]$ at receiver 3 (0.13 -- 5 Hz) }%
	\label{fig:Adaptive_LOH1}%
\end{figure}

Further, we investigate static adaptive mesh refinement, and stability for non conforming elements. 
Now the domain is discretized adaptively with uniform elements of $\Delta{x} = 17/9$~km in the bedrock when $x > 4$~km,  $\Delta{x} = 17/(9\times 3)$~km in the source region when $1< x \le 4$~km, and $\Delta{x} = 17/(9\times 9)$~km in the topmost layer $0\le x \le 1$~km (see the left panel in Figure \ref{fig:Adaptive_LOH1}). Note that in the topmost layer, closest to the surface, where surface waves are present the mesh is 3 times finer. We used GL nodes and degree $P = 3$ polynomial approximation. 
The effective grid resolution comparable to a finite difference grid size is $h = 53$~m in the topmost layer $0\le x \le 1$~km, $h = 3\times 53$~m in the source region $1< x \le 4$~km, $h = 9\times 53$~m in the bedrock $x > 4$~km. We evolve the wave fields with the time-step $\Delta{t} = 0.0026$~s  until the final time of 9~s. 

As above, we post-process the seismograms, and compute the EM and the PM. For Receiver 3 and the $z$ component of the particle velocity  $v_z$,  these quantities are displayed in the right panel of Figure \ref{fig:Adaptive_LOH1}. Note that EM$\sim 5\%$, and PM$\sim 2\%$, which are  reasonably moderate error tolerances.

\subsection{Surface waves in anisotropic media}
We consider simulations of elastic waves propagating in an orthotropic anisotropic medium \cite{Favretto-Cristini_etal2011, Komatitsch_etal2011}. Wave propagation in anisotropic elastic media  is more complex than in isotropic elastic media. From both analytical  and numerical perspectives, anisotropy introduces additional challenges, since a wave mode can propagate with different wave speeds in different directions. In fact, analytical solutions exist for only few problems. Numerical simulations become inevitable for studying wave propagation in anisotropic elastic media. We will  study the accuracy  of the method for surface waves in an anisotropic crystal, where analytical solutions exists. 
%
%
%
The crystal is apatite and the set-up is the same as in \cite{Komatitsch_etal2011}. The density of the medium and elastic constants are given in Table \ref{tab:Apatite}, with $c_{55} = c_{44}$,   $c_{66} = \left(c_{11}-c_{12}\right)/2$.
 \begin{table}[h!]
    \centering
    \caption{The density and elastic constants for the crystal Apatite.  With $c_{55} = c_{44}$   $c_{66} = \left(c_{11}-c_{12}\right)/2$ }
    \begin{tabular}{l||l|l|l|l|l|l|l|l|l|l}
    $\rho$ & $c_{11}$ & $c_{12}$ & $c_{13}$& $c_{23}$ & $c_{22}$ & $c_{33}$ & $c_{44}$   \\
    \cline{1-8}
    $3190$ kg/m$^{3}$ & $167$ GPa & $13.1$ GPa & $66$ GPa & $66$ GPa & $167$ GPa & $140$ GPa & $66.3$ GPa  \\
   \end{tabular}
   \label{tab:Apatite}
   \end{table}
  From \eqref{eq:anisotropic_wavespeed}, we can determine the wave speeds
  \begin{align*}
c_{px} = 7235~\ \text{m/s} ,\quad c_{shx} = 4559 ~\ \text{m/s}, \quad c_{svx} = 6945~\ \text{m/s},
\\
\nonumber
c_{py} = 7235~\ \text{m/s} ,\quad c_{shy} = 6945~\ \text{m/s} , \quad c_{svy} = 6945~\ \text{m/s},
\\
\nonumber
c_{pz} =  6624 ~\ \text{m/s} ,\quad c_{shz} = 4559 ~\ \text{m/s}, \quad c_{svz} = 4559 ~\ \text{m/s}.
\end{align*}
%
%
The computational domain is the cube $\left(x,y,z\right) \in [0, 20~\text{cm}] \times [0, 20~\text{cm}] \times [0, 20~\text{cm}]$.
Waves  are excited by  adding  the point source
\begin{align}\label{eq:point_source}
f(x,y,z,t) =   \delta_x(x-x_0)\delta_y(y-y_0)\delta_z(z-z_0)g(t),
\end{align}
to the first component of the momentum equation, where the source time function is
\[
g(t) = \cos\left[2\pi\left(t-t_0\right)f_0\right]e^{-2\left(t-t_0\right)^2f_0^2}, \quad t_0 = 3/(2f_0) + 5\times 10^{-6} ~\ \text{s}, \quad f_0 = 250 ~\ \text{kHz},
\]
 $ \delta_{\eta}(\eta-\eta_0)$ are the one dimensional Dirac delta function. The source is located at $x_0= 10~\text{cm}$, $y_0= 10~\text{cm}$,  $z_0 =0$, the surface of the crystal.
At the surface $z = 0$, we impose the free surface boundary condition and the absorbing boundary condition at the other 5 boundaries of the domain.

\begin{figure}[h!]
    \centering
    \includegraphics[width=0.525\textwidth]{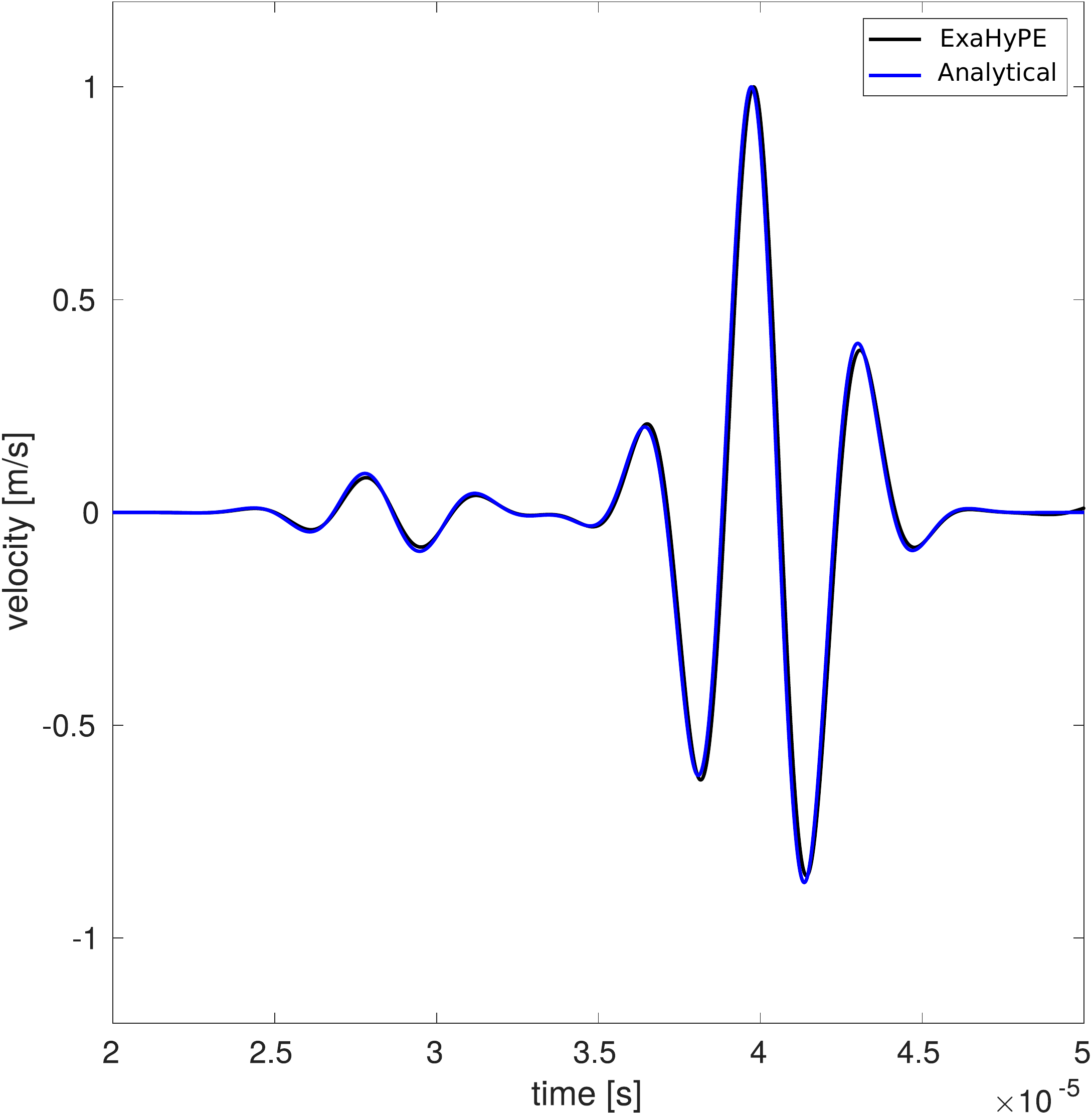}
    \caption{ A comparison of anisotropic wave propagation of the ExaSeis numerical and the analytical solution in apatite, an orthotropic anisotropic elastic medium at a single receiver (see text for details of the benchmark setup).}
    \label{fig:seismograms_anisotropic}
\end{figure}

We discretize the computational domain uniformly with 27 elements in each spatial direction, with a polynomial approximation of degree $P = 5$. 
We place a receiver at $x_0, = 10~\text{cm}$, $y_0, = 10~\text{cm}$,  $z_0 =15~\text{cm}$, and  evolve the wave fields, with the global time step $\Delta{t} = 0.015~\mu$s, until the final  time, $t = 50~\mu$s. 
In Figure \ref{fig:seismograms_anisotropic}, we compare the analytical solution with the numerically derived solution.
The numerical solution matches the analytical solution very well, with less than $1\%$ relative error.
\subsection[The Zugspitze model]{Wave propagation in complex geometries}
We will now demonstrate the potential of the method in modeling elastic wave propagation in a 3D domain with geometrically complex free surface topography.
 We consider a 3D setup incorporating the strong topography contrasts of Mount Zugspitze, Germany. 
\begin{figure}[h!]
\center
\includegraphics[width=\textwidth, angle =0]{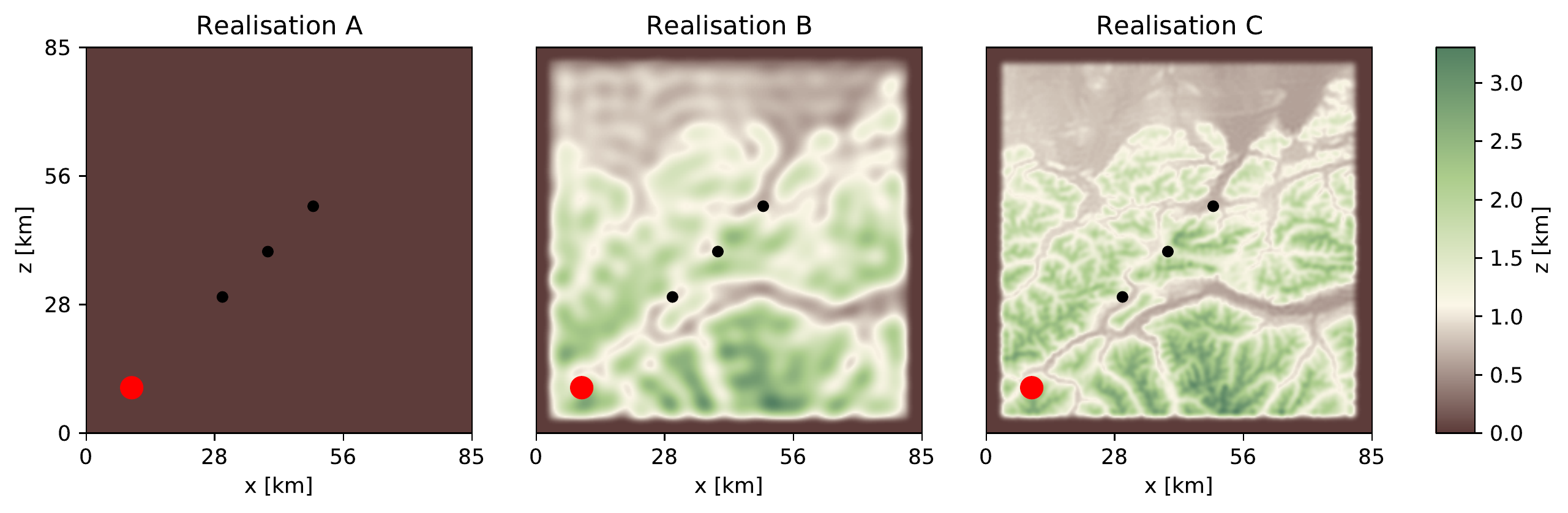}
\caption{The Mount Zugspitze model set up for Realisations A, B and C. The red dot ({\red \bf $\bullet$}) at ($x=10$ km, $z=10$ km) depicts the epicenter  of a buried moment tensor point source and the black dots ({\black \bf$\bullet$})  indicate Station 1: ($x=30$ km, $z=30$ km), Station 2: ($x=40$ km, $z=40$ km)  and Station 3: ($x=50$ km, $z=50$ km), which are the receiver stations where the solutions are sampled. Station 2 is collocated with the top ($x=40$ km, $z=40$ km) of Mount Zugspitze.}
  \label{fig:topographydata_0}
\end{figure}

\begin{figure}[h!]
\center
\includegraphics[width=0.45\textwidth, angle =0]{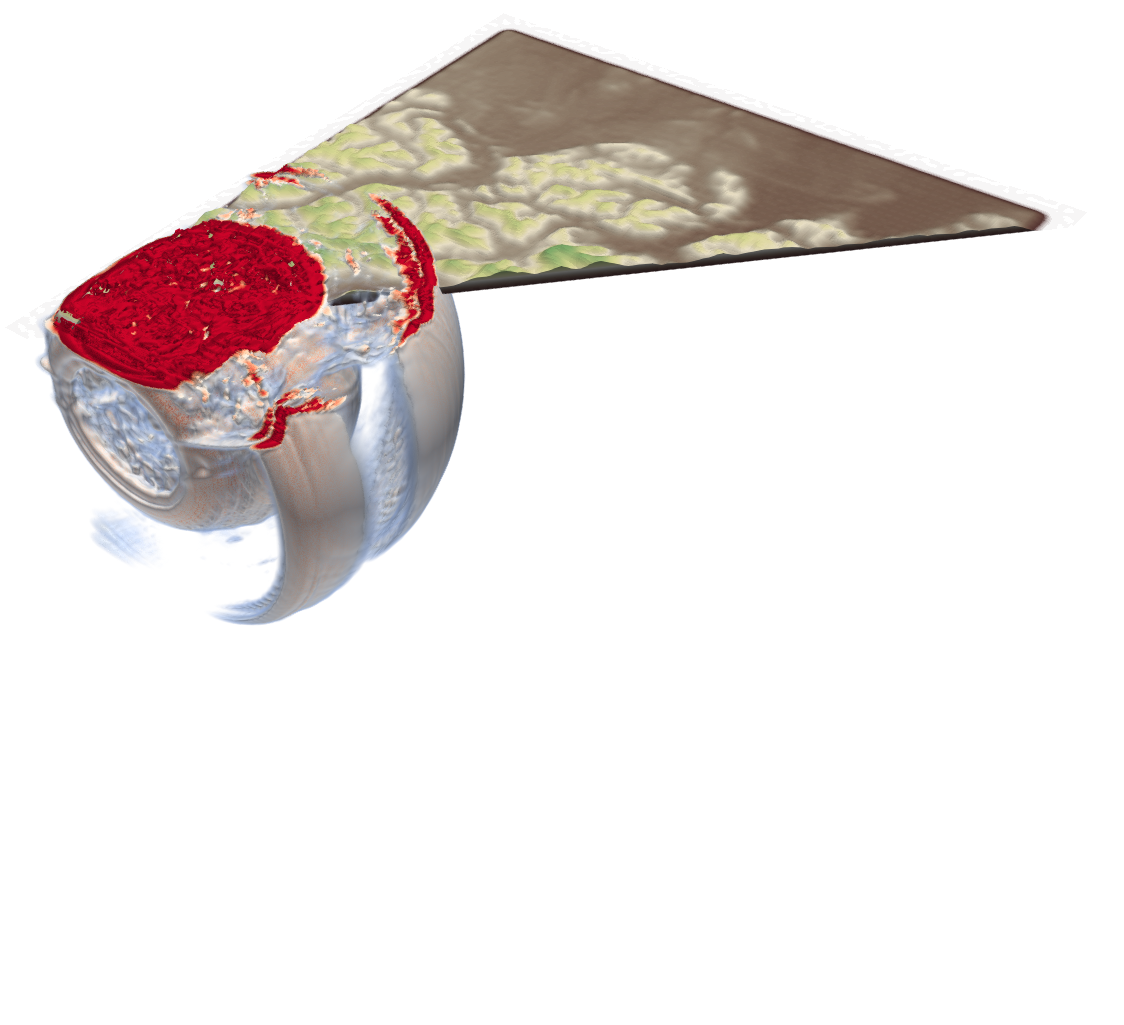}
\includegraphics[width=0.45\textwidth, angle =0]{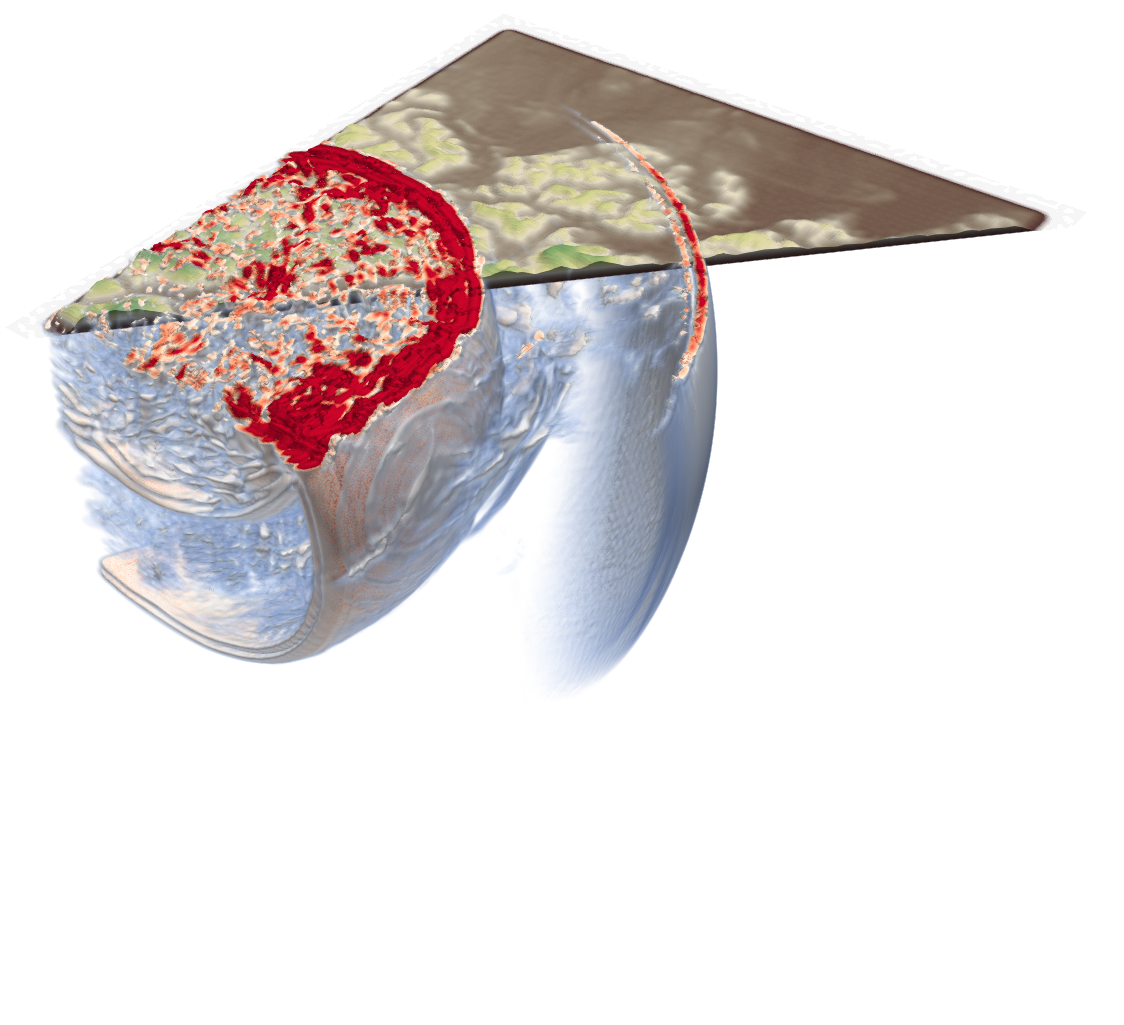}
\includegraphics[width=0.45\textwidth, angle =0]{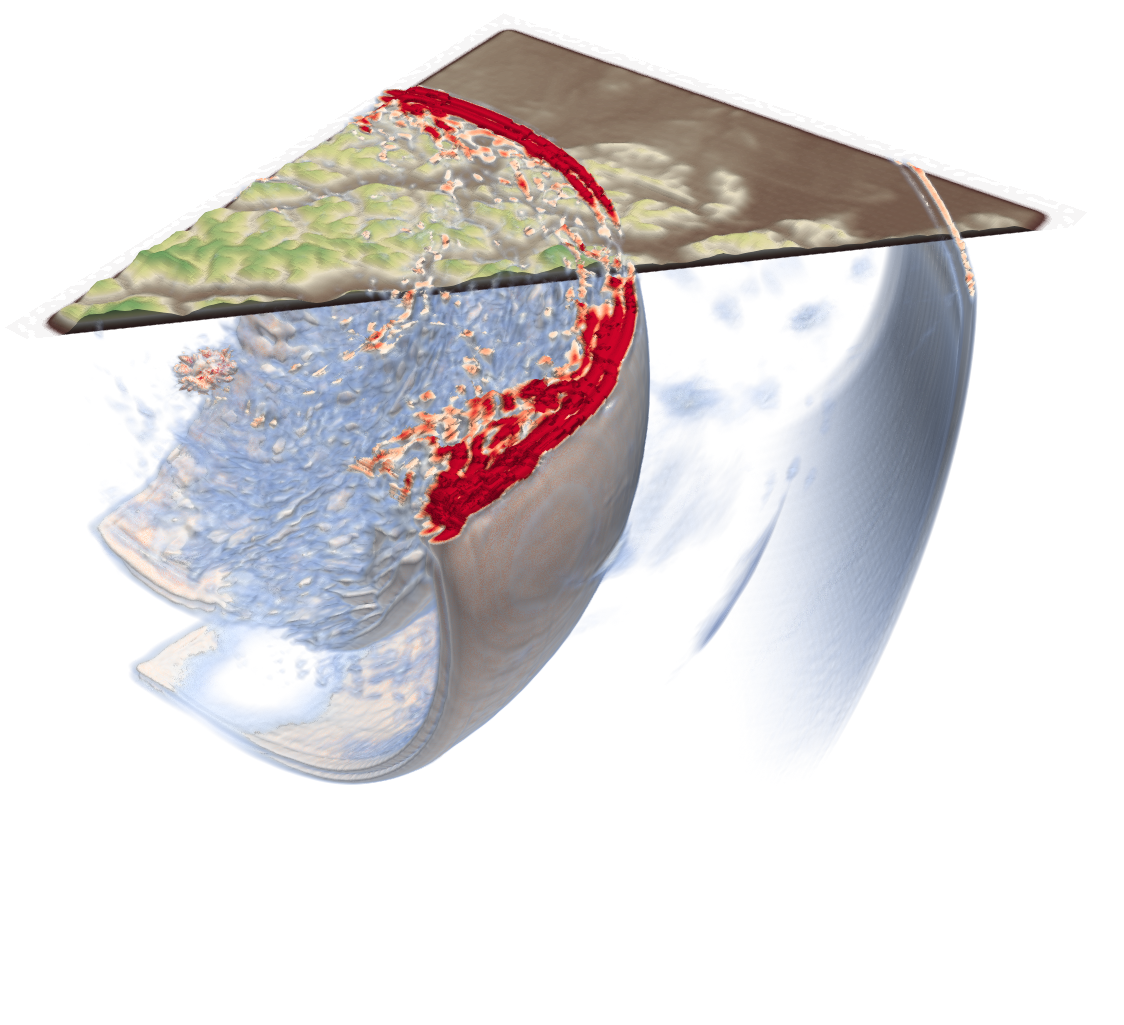}
\includegraphics[width=0.45\textwidth, angle =0]{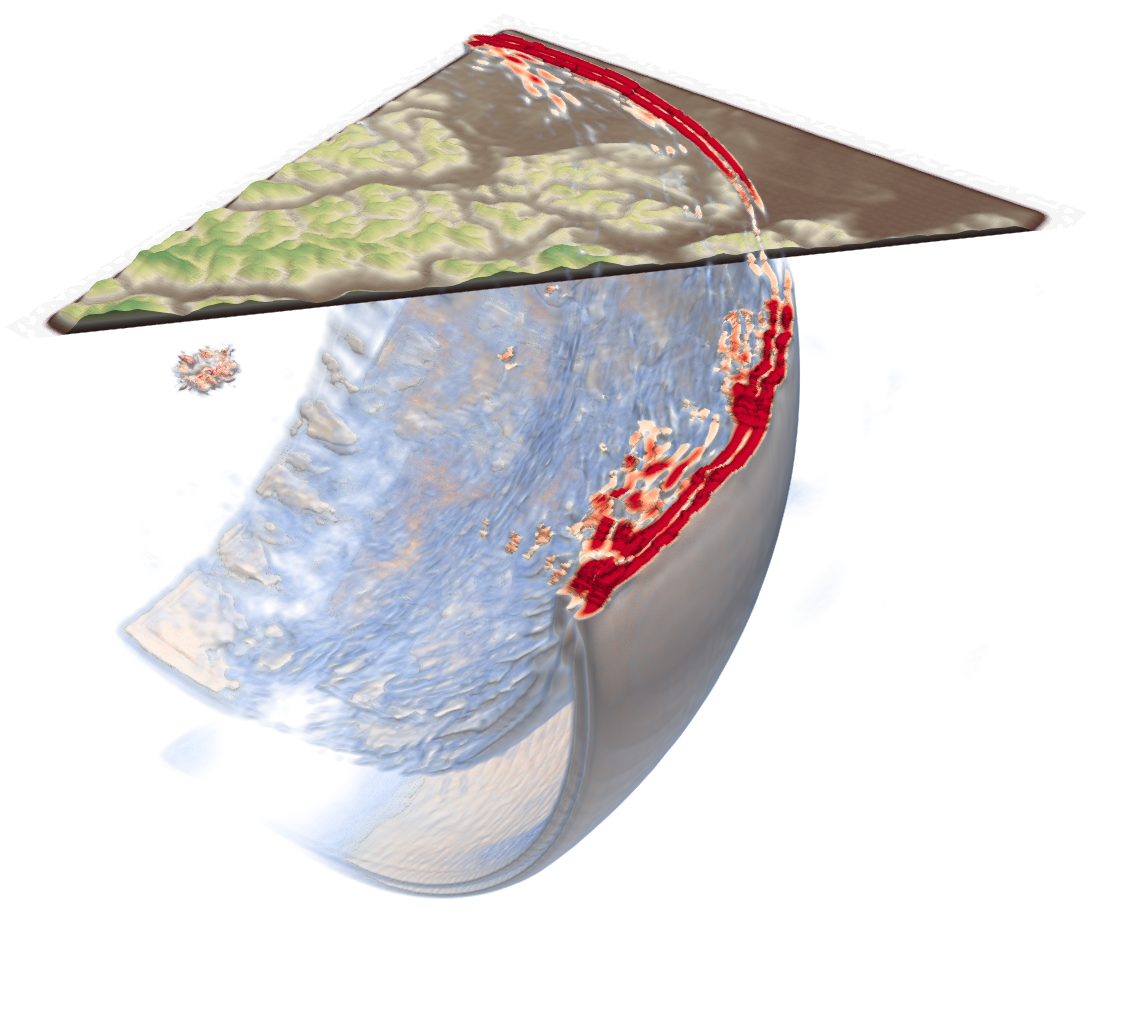}
  \caption{3D snapshots of the absolute velocity of the propagating seismic wavefield for the Zugspitze model at t=5~s, t=10~s, t=15~s and t=20~s simulated with ExaSeis.}
  \label{fig:topographydata}
\end{figure}
 Mount Zugspitze is reaching about $3$~km above sea level and is the highest peak of the Wetterstein mountain range as well as the highest mountain in Germany.
The topography of this region is complex, with large deviations from planarity across the Earth's surface, and constitute a significant challenge for numerical solvers. Effective numerical resolution of the scattered wavefield is very important for application purposes, and for example may affect the resolution of regional seismic tomography performed  by the AlpArray project \cite{Hetenyi2018}.
In particular, due to the high frequency wave modes generated by scattering  from the complex non-planar topography, accurate and stable numerical simulation of seismic wave propagation in this region is both numerically and computationally arduous. 
The Zugspitze model has no analytical solutions. We will verify accuracy by making comparisons with the reference data in \cite{DuruFungWilliams2020} simulated with WaveQLab3D \cite{DuruandDunham2016}, a petascale finite difference elastic wave solver.


The topography data spans a rectangular surface area $(x, z) \in [0, 80~\text{km}]\times [0, 80~\text{km}] $, with $100~$m resolution of the elevation data. Down dip we truncate the domain at $y = 80~\ $km, so that the computational domain corresponds to the modulated cuboid $(x,y,z)\in [0, 80~\text{km}]\times [\widetilde{y}, 80~\text{km}]\times [0, 80~\text{km}]$, where $\widetilde{y}(x, z)$ is the elevation data. 
We have processed the topography data by passing it through a band limited filter, see Figure \ref{fig:topographydata_0}, obtaining 3 realisations: Realisation A,  Realisation B and Realisation C. Realisation A is a flat topography and includes low frequencies,  Realisation B includes intermediate frequencies and Realisation C includes high frequencies.
 
At each truncated boundary in $x$-axis and $z$-axis, and down dip at $y = 80$~m, we have included a 5~km absorbing layer where PML boundary conditions \cite{DuruRannabauerGabrielKreissBader2019} are implemented  to prevent artificial reflections  from the computational boundaries from contaminating the solution. A stable numerical implementation of the PML for 3D linear elastodynamics is nontrivial and allows the generation of high quality seismograms.

We consider the homogenous material properties 
$$\rho = 2670~\ \text{kg/m}^3, \quad c_p = 6000~\ \text{m/s}, \quad c_s = 3464~\ \text{m/s}.$$

The domain is discretized uniformly with $241$ elements in each spatial direction and  degree $P =3$ polynomial approximation. The effective sub-cell grid resolution comparable to a finite difference grid size is $h = \Delta{x}/(P+1) = 93.36~$m.  We evolve the solution with the time step \eqref{eq:time_step} until the final time $t = 30$~s. 

%

\begin{figure}[h!]
\begin{subfigure}
    \centering
\stackunder[5pt]{\includegraphics[draft=false,width=0.75\columnwidth]{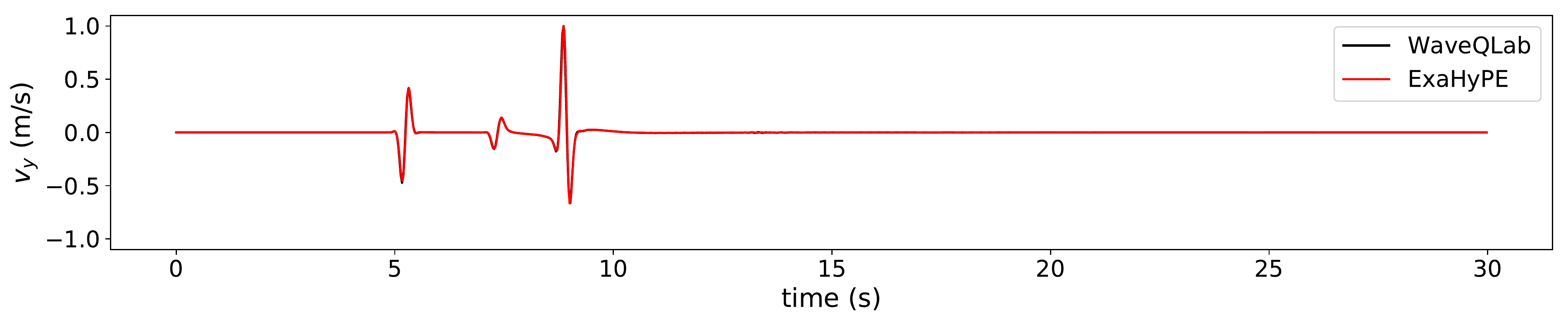}}{Realisation A}%
\hspace{0.0cm}%
\stackunder[5pt]{\includegraphics[draft=false,width=0.75\columnwidth]{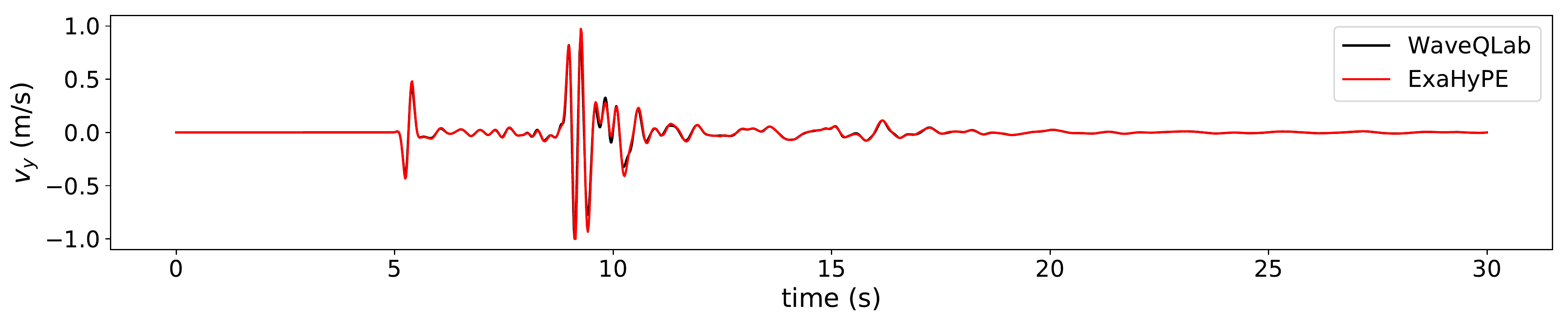}}{Realisation B}%
\hspace{0.0cm}%
\stackunder[5pt]{\includegraphics[draft=false,width=0.75\columnwidth]{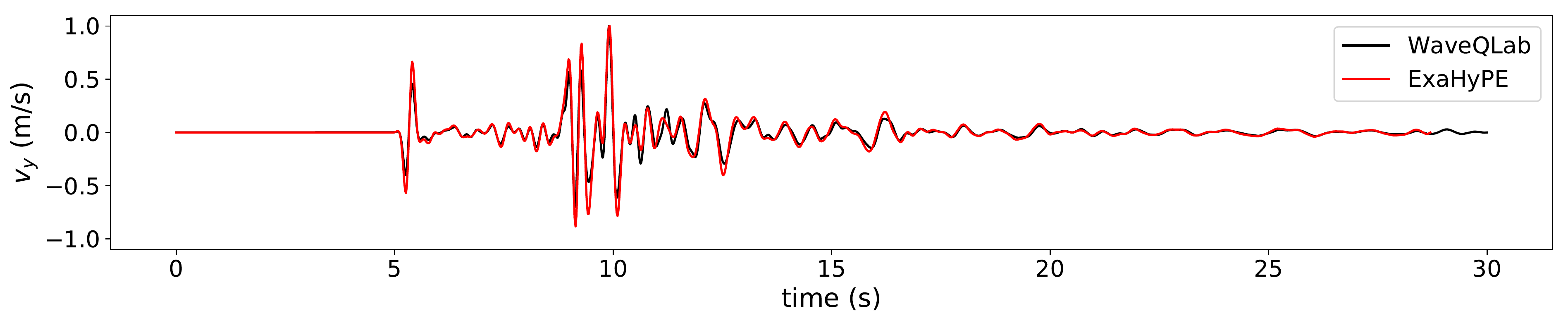}}{Realisation C}%
     \end{subfigure}
    \caption{Numerical verification of the time series in $v_y[\mathrm{m/s}]$ (synthetic seismogram) recorded at Station 1 ($x=30$~km, $z=30$~km) comparing the ExaSeis application of ExaHyPE to the finite difference solver WaveQLab3D.}
    \label{fig:zugs_r1}
\end{figure}
\begin{figure}[h!]
\begin{subfigure}
    \centering
\stackunder[5pt]{\includegraphics[draft=false,width=0.75\columnwidth]{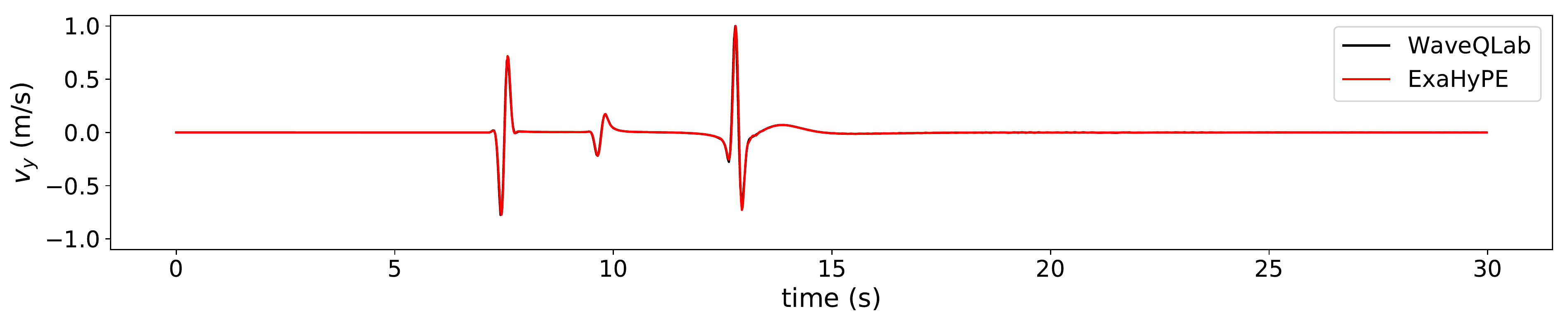}}{Realisation A}%
\hspace{0.0cm}%
\stackunder[5pt]{\includegraphics[draft=false,width=0.75\columnwidth]{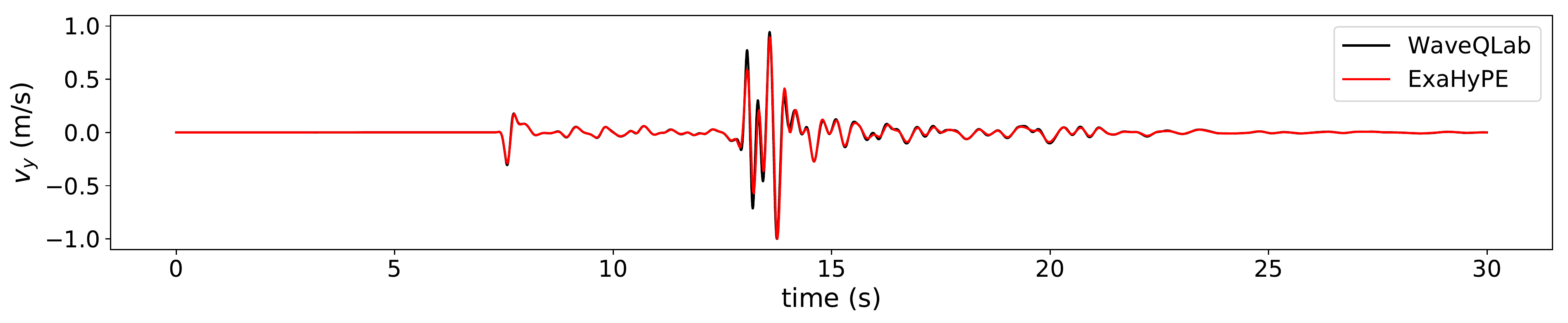}}{Realisation B}%
\hspace{0.0cm}%
\stackunder[5pt]{\includegraphics[draft=false,width=0.75\columnwidth]{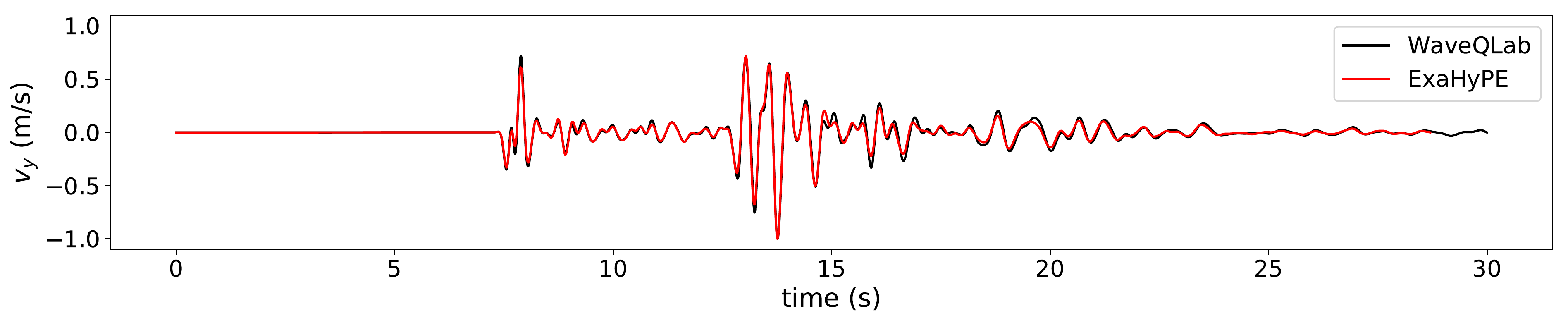}}{Realisation C}%
     \end{subfigure}
    \caption{Numerical verification of the time series in $v_y[\mathrm{m/s}]$ (synthetic seismogram) recorded at Station 2 ($x=40$~km, $z=40$~km) located at the top of Mount Zugspitze comparing the ExaSeis application of ExaHyPE to the finite difference solver WaveQLab3D. }
    \label{fig:zugs_r2}
\end{figure}
\begin{figure}[h!]
\begin{subfigure}
    \centering
\stackunder[5pt]{\includegraphics[draft=false,width=0.75\columnwidth]{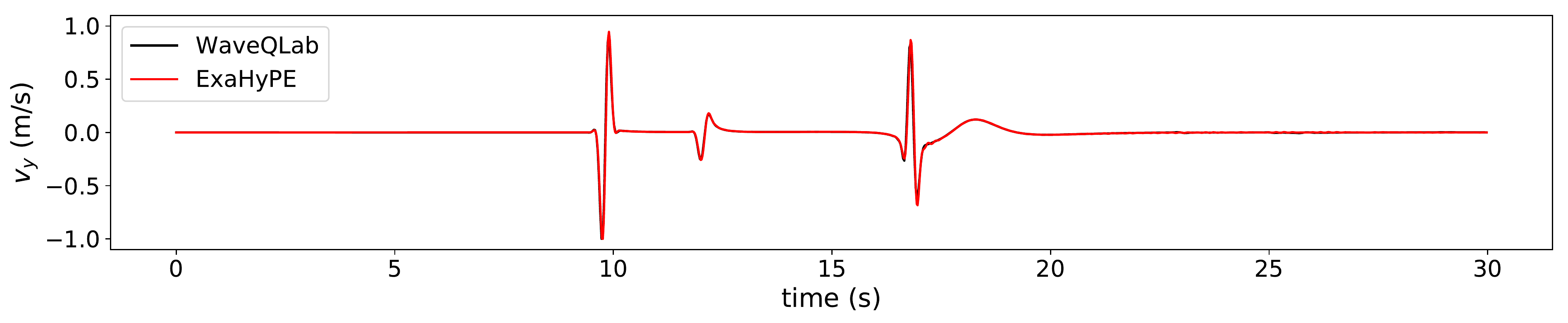}}{Realisation A}%
\hspace{0.0cm}%
\stackunder[5pt]{\includegraphics[draft=false,width=0.75\columnwidth]{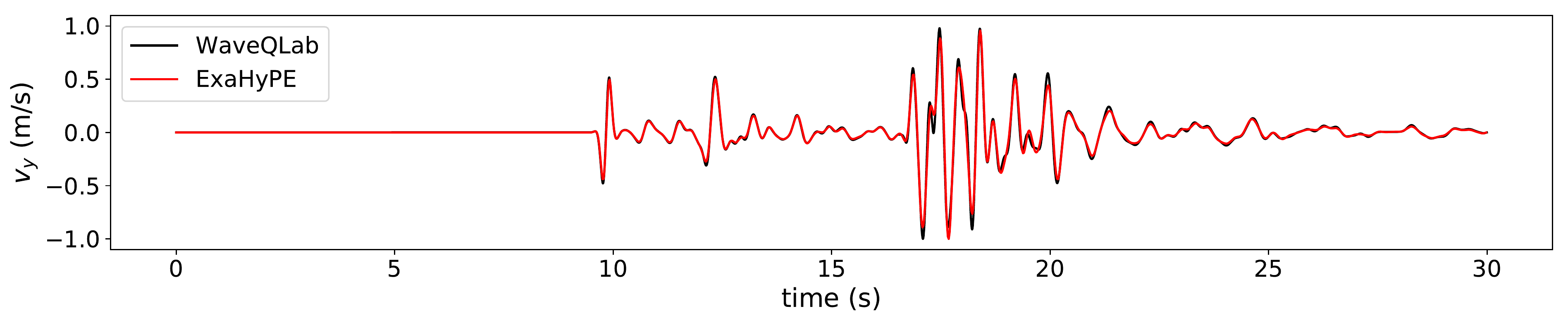}}{Realisation B}%
\hspace{0.0cm}%
\stackunder[5pt]{\includegraphics[draft=false,width=0.75\columnwidth]{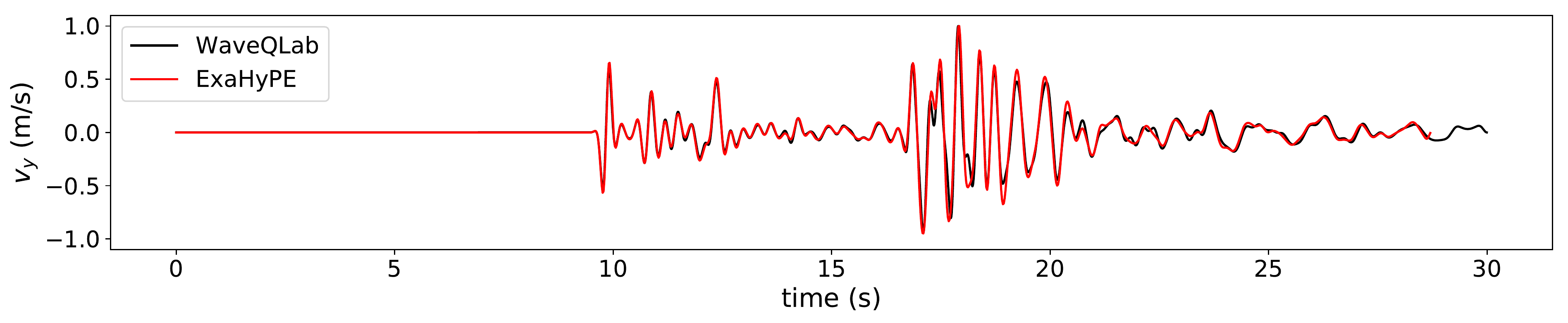}}{Realisation C}%
     \end{subfigure}
    \caption{Numerical verification of the time series in $v_y[\mathrm{m/s}]$ (synthetic seismogram) recorded at Station 3 ($x=50$~km, $z=50$~km) comparing the ExaSeis application of ExaHyPE to the finite difference solver WaveQLab3D.}
    \label{fig:zugs_r3}
\end{figure}

 Snapshots of numerical simulations are shown in Figure \ref{fig:topographydata} and the seismograms are shown in Figures \ref{fig:zugs_r1}--\ref{fig:zugs_r3}, for the 3 receiver stations and the 3 topography realisations. For a low frequency topography model, Realisation A, note that the seismograms are dominated by the direct arrivals and waves propagate coherently.
 The seismograms for Realisations B and C are dominated by high frequency scattered waves, which are particularly present in the coda waves.
%
We compare seismograms from ExaSeis and WaveQLab3D with $h=100$~m finite difference grid spacing.  We observe a near perfect agreement of the seismograms at all frequencies, see Figures \ref{fig:zugs_r1}--\ref{fig:zugs_r3}. However at high frequencies there are tiny differences, which will diminish as we increase the resolution.

\section{Summary and outlook}
In this paper, we present a  new energy-stable discontinuous Galerkin approximation of the elastic wave equation in general and geometrically complex 3D media, using the physics-based numerical penalty-flux \cite{DuruGabrielIgel2017}. As opposed to the Godunov flux, the physics-based flux does not require a complete eigenvector-eigenvalue decomposition of the spatial coefficient matrices. By construction, our numerical flux is upwind and yields a discrete energy estimate analogous to the continuous energy estimate. The discrete energy estimate holds for conforming and non-conforming curvilinear elements. The ability to handle non-conforming curvilinear meshes allows for flexible adaptive mesh refinement strategies. The numerical scheme has been implemented in the ExaSeis application in ExaHyPE \cite{ExaHyPE2019}, a simulation engine for hyperbolic PDEs on adaptive structured meshes.
 %
 Numerical experiments are presented in 3D isotropic and anisotropic media demonstrating stability and accuracy. Finally, we present numerical verification in a regional geophysical wave propagation problem in an Earth model with  geometrically complex free-surface topography.
 
 We utilize a recent extension of the method to discretize PML boundary conditions for elastodynamics \cite{DuruRannabauerGabrielKreissBader2019}.
We expect, that the proposed numerical method will extend to model linear and nonlinear boundary and interface wave phenomena. In a forthcoming paper we will extend the method to nonlinear friction problems, and present numerical simulations of nonlinear dynamic earthquake ruptures on dynamically adaptive meshes, embedded in 3D geometrically complex solid earth models.
\section*{Acknowledgments}
The work presented in this paper was enabled by funding from the European Union's Horizon 2020 research and innovation program under grant agreements no. 671698 (ExaHyPE), no. 823844 (ChEESE) and no. 852992 (TEAR).
The authors also acknowledge support by the German Research Foundation (DFG) (grants no. GA 2465/2-1, GA 2465/3-1), by KAUST-CRG (grant no. ORS-2017-CRG6 3389.02) and by KONWIHR (project NewWave). 
A.L. is supported by the Swiss Federal Institute of Technology grant (project ETH-10 17-2).
Computing resources were provided by the Institute of Geophysics of LMU Munich \cite{Oeser2006}, the Leibniz Supercomputing Centre (SuperMUC-NG project pr63qo) and  the KAUST Shaheen Supercomputing Laboratory (project k1488).
\newline{}
\hfill{\includegraphics[angle=90, width=0.15\textwidth]{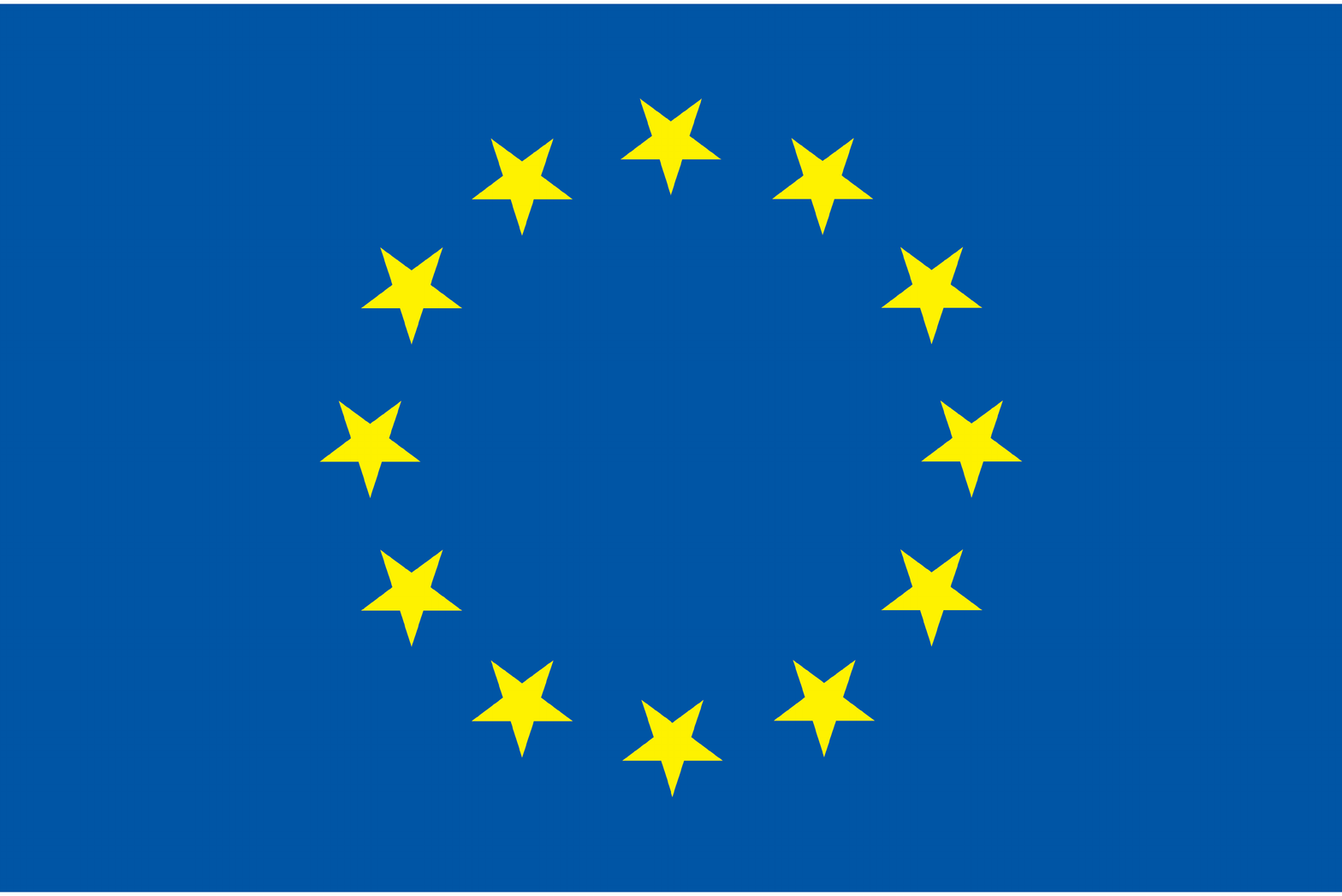}}
\\
The first author KD would like to thank Dimitri Komatitsch (1970--2019) for his help, fruitful discussions and for providing the analytical solution for elastic surface waves in an anisotropic medium.
\appendix
 \section{The Arbitrary DERivative (ADER) time integration}\label{sec:ADER}
 In this section we will summarize the ADER time-stepping scheme. For more elaborate discussions, we refer the reader to \cite{Toro1999, DumbserKaeser2006, delaPuenteAmpueroKaser2009}.
 To begin,  we rewrite the semi-discrete approximation  \eqref{eq:gen_hyp_transformed_discrete} as a system first order ODEs
\begin{equation}\label{eq:ODE}
\begin{split}
\frac{d \bar{\mathbf{Q}} }{dt} = \underbrace{D\bar{\mathbf{Q}}}_{\text{PDE}} + \underbrace{F\bar{\mathbf{Q}}}_{\text{Num. flux} \to 0},
\end{split}
\end{equation}
with
  \begin{align}\label{eq:gen_hyp_transformed_discrete_ADER}
D\bar{\mathbf{Q}}:= \widetilde{\mathbf{P}} \left(\grad_{D} \bullet \mathbf{F} \left( \bar{\mathbf{Q}} \right) + \sum_{\xi= q, r, s}\mathbf{B}_\xi \left(\grad_{D} \bar{\mathbf{Q}} \right)\right), \quad F\bar{\mathbf{Q}} := - \widetilde{\mathbf{P}} \mathbf{Flux}\left(\bar{\mathbf{Q}}\right),
\end{align}
where $D\bar{\mathbf{Q}}$ is the discrete spatial operator is split into the derivative term, emanating from the PDE, and $F\bar{\mathbf{Q}}$ is the numerical flux fluctuation term, incorporating the boundary and interface conditions.


The numerical flux fluctuation is a very small term, $F\bar{\mathbf{Q}}  \approx 0$, and will vanish $F\bar{\mathbf{Q}} \to 0$ in the limit of mesh refinement $\Delta{t} \to 0$. 

We now introduce the discrete time variables $t_k \le t \le t_{k+1}$, $\Delta{t}_k = t_{k+1}- t_{k}$, and the pseudo time variable $\tau = t -t_k$ such that $0 \le \tau \le \Delta{t}_k$, and  ${d}/{d\tau} = {d}/{dt} $.
 Going from the current time level $\tau =0$ to the next time level $\tau = \Delta{t}_k$, we integrate the ODE \eqref{eq:ODE}, exactly having
\begin{equation}\label{eq:one_step_1}
\begin{split}
\small
\bar{\mathbf{Q}}(\Delta{t}) &= \bar{\mathbf{Q}}(0)  + \int_0^{\Delta{t}_k}{D\bar{\mathbf{Q}}} d\tau + \int_0^{\Delta{t}_k}{F\bar{\mathbf{Q}}} d\tau, \quad \\
& =  \bar{\mathbf{Q}}(0)  + D\int_0^{\Delta{t}_k}{\bar{\mathbf{Q}}} d\tau + F\int_0^{\Delta{t}_k}{\bar{\mathbf{Q}}} d\tau,
\end{split}
\end{equation}
where the second equality follows from linearity. If we can evaluate the integrals $\int_0^{\Delta{t}_k}{\bar{\mathbf{Q}}} d\tau$ in \eqref{eq:one_step_1} exactly, then the time integration in \eqref{eq:one_step_1} is exact. However, exact time integration is  possible only in the most trivial case where the right hand side of \eqref{eq:ODE} vanish identically for all components. Now, we will make an important approximation.  We assume that the time step $\Delta{t}$ is sufficiently small, such that
\begin{align}\label{eq:one_step_2}
 \frac{d\bar{\mathbf{Q}}(\tau)}{d\tau} \approx D\bar{\mathbf{Q}}, \quad F\bar{\mathbf{Q}} \approx 0,
\end{align}
are reasonable approximations. Next we construct the predictor, $ \widetilde{\bar{\mathbf{Q}}}(\tau)$, by Taylor expansions of the solution around $\tau = 0$ and replace the time derivatives with  spatial operator in  \eqref{eq:one_step_2}, we have
\begin{align}\label{eq:one_step_3}
 \widetilde{\bar{\mathbf{Q}}}(\tau) = \bar{\mathbf{Q}}(0) + \tau  \frac{d\bar{\mathbf{Q}}(0)}{d\tau}  + \frac{\tau^2}{2} \frac{d^2\bar{\mathbf{Q}}(0)}{d\tau^2}  + ...
\approx \sum_{m = 0}^{P}\frac{\tau^m}{m !} D^m \bar{\mathbf{Q}}(0),
\end{align}
where $P$ is the polynomial degree used in the spatial approximation. We can now approximate the integrals in \eqref{eq:one_step_1} using the predictor. The result of this integration is called the time average, $ \bar{\bar{\mathbf{Q}}}(0) $,
\begin{align}\label{eq:one_step_4}
\int_0^{\Delta{t}_k}{\bar{\mathbf{Q}}} d\tau \approx { \bar{\bar{\mathbf{Q}}}(0) = \int_{0}^{\Delta{t}_k} \widetilde{\bar{\mathbf{Q}}}(\tau) d\tau
= \sum_{m = 0}^{P}\frac{\Delta{t}_k^{(m+1)}}{(m +1) !} D^m \bar{\mathbf{Q}}(0)}.
\end{align}
By replacing the integrals in \eqref{eq:one_step_1} with the time average ${\bar{\bar{\mathbf{Q}}}(0)} $, we derive a high order accurate, explicit, one-step, time integration  scheme
\begin{align}\label{eq:one_step_5}
\small
\bar{\mathbf{Q}}(\Delta{t}) & =  \bar{\mathbf{Q}}(0)  + D\bar{\bar{\mathbf{Q}}}(0)  + F\bar{\bar{\mathbf{Q}}}(0).
\end{align}
Note that the numerical flux fluctuations, $F\bar{\bar{\mathbf{Q}}}(0)$,   is evaluated only once for any order of approximation. This is opposed to Runge--Kutta methods or standard Taylor series methods where  the numerical flux fluctuation is included in the spatial operator, $D^m \to (D + {F})^m$, to  approximate higher time derivatives in the Taylor series terms. For the ADER scheme, the fact that the numerical flux fluctuation is evaluated only once for any order implies that most of the computations are performed within the element to compute the predictor in \eqref{eq:one_step_3} and  the time average in \eqref{eq:one_step_4}. This has a huge impact in high performance computing applications, since we can design efficient communication avoiding parallel algorithms. Since the predictor $\widetilde{\bar{\mathbf{Q}}}(\tau)$ is defined in the entire time interval $0 \le \tau \le \Delta{t}_k$, the  ADER scheme, \eqref{eq:one_step_1}--\eqref{eq:one_step_5}, is also easily amenable to local time-stepping methods. When the ADER time stepping scheme is combined with the DG spatial approximation the fully discrete scheme is often referred as the ADERDG method \cite{DumbserKaeser2006, delaPuenteAmpueroKaser2009}. For a  DG polynomial approximation of degree $P$, a stable ADERDG method is $(P+1)$th order accurate in both space and time.
%
%
%
%
%
%
\bibliographystyle{elsarticle-num-names}
\bibliography{mybib}{}

\begin{thebibliography}{61}
\expandafter\ifx\csname natexlab\endcsname\relax\def\natexlab#1{#1}\fi
\providecommand{\url}[1]{\texttt{#1}}
\providecommand{\href}[2]{#2}
\providecommand{\path}[1]{#1}
\providecommand{\DOIprefix}{doi:}
\providecommand{\ArXivprefix}{arXiv:}
\providecommand{\URLprefix}{URL: }
\providecommand{\Pubmedprefix}{pmid:}
\providecommand{\doi}[1]{\href{http://dx.doi.org/#1}{\path{#1}}}
\providecommand{\Pubmed}[1]{\href{pmid:#1}{\path{#1}}}
\providecommand{\bibinfo}[2]{#2}
\ifx\xfnm\relax \def\xfnm[#1]{\unskip,\space#1}\fi
\bibitem[{Hill(1973)}]{ReedHill1973}
\bibinfo{author}{W.~H. R. T.~R. Hill}, \bibinfo{title}{Triangular mesh methods
  for the neutron transport equation}, \bibinfo{howpublished}{Technical Report
  LA-UR-73-479, Los Alamos National Laboratory, Los Alamos, New Mexico, USA},
  \bibinfo{year}{1973}.
\bibitem[{Cockburn and Shu(1989)}]{CockburnShu1989}
\bibinfo{author}{B.~Cockburn}, \bibinfo{author}{C.~W. Shu},
\newblock \bibinfo{title}{Tvb runge-kutta local projection discontinuous
  galerkin finite element method for conservation laws ii: general framework},
\newblock \bibinfo{journal}{Math Comput.} \bibinfo{volume}{52}
  (\bibinfo{year}{1989}) \bibinfo{pages}{411--435}.
\bibitem[{Cockburn et~al.(1990)Cockburn, Hou, and Shu}]{CockburnHouShu1990}
\bibinfo{author}{B.~Cockburn}, \bibinfo{author}{S.~Hou}, \bibinfo{author}{C.~W.
  Shu},
\newblock \bibinfo{title}{The runge-kutta local projection discontinuous
  galerkin finite element method for conservation laws iv},
\newblock \bibinfo{journal}{J. Comput. Phys.} \bibinfo{volume}{54}
  (\bibinfo{year}{1990}) \bibinfo{pages}{545--581}.
\bibitem[{Hesthaven and Warburton(2008)}]{HesthavenWarburton2008}
\bibinfo{author}{J.~Hesthaven}, \bibinfo{author}{T.~Warburton},
  \bibinfo{title}{Nodal Discontinuous Galerkin Methods: Algorithms, Analysis,
  and Applications}, \bibinfo{publisher}{Springer, New York},
  \bibinfo{year}{2008}.
\bibitem[{Hesthaven and Warburton(2002)}]{HesthavenWarburton2002}
\bibinfo{author}{J.~S. Hesthaven}, \bibinfo{author}{T.~Warburton},
\newblock \bibinfo{title}{Nodal high-order methods on unstructured grids: I.
  time-domain solution of maxwell's equations},
\newblock \bibinfo{journal}{J. Comput. Phys.} \bibinfo{volume}{181}
  (\bibinfo{year}{2002}) \bibinfo{pages}{186--221}.
\bibitem[{Dumbser and K{\"a}ser(2006)}]{DumbserKaeser2006}
\bibinfo{author}{M.~Dumbser}, \bibinfo{author}{M.~K{\"a}ser},
\newblock \bibinfo{title}{An arbitrary high-order discontinuous {G}alerkin
  method for elastic waves on unstructured meshes -- {II.} the
  three-dimensional isotropic case},
\newblock \bibinfo{journal}{Geophysical Journal International}
  \bibinfo{volume}{167} (\bibinfo{year}{2006}) \bibinfo{pages}{319--336}.
  \DOIprefix\doi{10.1111/j.1365-246X.2006.03120.x}.
\bibitem[{Burstedde et~al.(2010)Burstedde, Ghattas, Gurnis, Isaac, Stadler,
  Warburton, and Wilcox}]{Burstedde:2010}
\bibinfo{author}{C.~Burstedde}, \bibinfo{author}{O.~Ghattas},
  \bibinfo{author}{M.~Gurnis}, \bibinfo{author}{T.~Isaac},
  \bibinfo{author}{G.~Stadler}, \bibinfo{author}{T.~Warburton},
  \bibinfo{author}{L.~C. Wilcox},
\newblock \bibinfo{title}{Extreme-scale {AMR}},
\newblock in: \bibinfo{booktitle}{SC10: Proc. Int. Conf. HPC, Networking,
  Storage and Analysis}, \bibinfo{year}{2010}.
\bibitem[{Breuer et~al.(????)Breuer, Heinecke, Rettenberger, Bader, Gabriel,
  and Pelties}]{Breuer2014}
\bibinfo{author}{A.~Breuer}, \bibinfo{author}{A.~Heinecke},
  \bibinfo{author}{S.~Rettenberger}, \bibinfo{author}{M.~Bader},
  \bibinfo{author}{A.-A. Gabriel}, \bibinfo{author}{C.~Pelties},
\newblock in: \bibinfo{booktitle}{{Supercomputing. ISC 2014. Lecture Notes in
  Computer Science, vol 8488}}, ????
\bibitem[{Heinecke et~al.(2014)Heinecke, A.~Breuer, Bader, Gabriel, Pelties,
  A.~Bode, Vaidyanathan, Smelyanskiy, and Dubey}]{Heineckeetal2014}
\bibinfo{author}{A.~Heinecke}, \bibinfo{author}{S.~R. A.~Breuer},
  \bibinfo{author}{M.~Bader}, \bibinfo{author}{A.-A. Gabriel},
  \bibinfo{author}{C.~Pelties}, \bibinfo{author}{X.-K.~L. A.~Bode, W.~Barth},
  \bibinfo{author}{K.~Vaidyanathan}, \bibinfo{author}{M.~Smelyanskiy},
  \bibinfo{author}{P.~Dubey}, \bibinfo{title}{Petascale high order dynamic
  rupture earthquake simulations on heterogeneous supercomputers},
  \bibinfo{howpublished}{In: Proceedings of SC 2014 New Orleans, LA},
  \bibinfo{year}{2014}.
\bibitem[{Uphoff et~al.(2017)Uphoff, Rettenberger, Bader, Madden, Ulrich,
  Wollherr, and Gabriel}]{Uphoff:2017}
\bibinfo{author}{C.~Uphoff}, \bibinfo{author}{S.~Rettenberger},
  \bibinfo{author}{M.~Bader}, \bibinfo{author}{E.~H. Madden},
  \bibinfo{author}{T.~Ulrich}, \bibinfo{author}{S.~Wollherr},
  \bibinfo{author}{A.-A. Gabriel},
\newblock \bibinfo{title}{Extreme scale multi-physics simulations of the
  tsunamigenic 2004 {Sumatra} megathrust earthquake},
\newblock in: \bibinfo{booktitle}{SC '17: Proc. Int. Conf. HPC, Networking,
  Storage and Analysis}, \bibinfo{publisher}{ACM}, \bibinfo{year}{2017}.
  \URLprefix \url{https://dl.acm.org/citation.cfm?id=3126948}.
\bibitem[{Bielak et~al.(2010)Bielak, Graves, Olsen, Taborda,
  Ram{\'\i}rez-Guzm{\'a}n, Day, Ely, Roten, Jordan, Maechling
  et~al.}]{bielak2010shakeout}
\bibinfo{author}{J.~Bielak}, \bibinfo{author}{R.~W. Graves},
  \bibinfo{author}{K.~B. Olsen}, \bibinfo{author}{R.~Taborda},
  \bibinfo{author}{L.~Ram{\'\i}rez-Guzm{\'a}n}, \bibinfo{author}{S.~M. Day},
  \bibinfo{author}{G.~P. Ely}, \bibinfo{author}{D.~Roten},
  \bibinfo{author}{T.~H. Jordan}, \bibinfo{author}{P.~J. Maechling}, et~al.,
\newblock \bibinfo{title}{{The ShakeOut earthquake scenario: Verification of
  three simulation sets}},
\newblock \bibinfo{journal}{Geophysical Journal International}
  \bibinfo{volume}{180} (\bibinfo{year}{2010}) \bibinfo{pages}{375--404}.
\bibitem[{Chaljub et~al.(2010)Chaljub, Moczo, Tsuno, Bard, Kristek, K{\"a}ser,
  Stupazzini, and Kristekova}]{chaljub2010grenoble}
\bibinfo{author}{E.~Chaljub}, \bibinfo{author}{P.~Moczo},
  \bibinfo{author}{S.~Tsuno}, \bibinfo{author}{P.-Y. Bard},
  \bibinfo{author}{J.~Kristek}, \bibinfo{author}{M.~K{\"a}ser},
  \bibinfo{author}{M.~Stupazzini}, \bibinfo{author}{M.~Kristekova},
\newblock \bibinfo{title}{{Quantitative comparison of four numerical
  predictions of 3D ground motion in the Grenoble Valley, France}},
\newblock \bibinfo{journal}{Bulletin of the Seismological Society of America}
  \bibinfo{volume}{100} (\bibinfo{year}{2010}) \bibinfo{pages}{1427--1455}.
\bibitem[{Graves et~al.(2011)Graves, Jordan, Callaghan, Deelman, Field, Juve,
  Kesselman, Maechling, Mehta, Milner, Okaya, Small, and
  Vahi}]{Graves_etal2011}
\bibinfo{author}{R.~Graves}, \bibinfo{author}{T.~H. Jordan},
  \bibinfo{author}{S.~Callaghan}, \bibinfo{author}{E.~Deelman},
  \bibinfo{author}{E.~Field}, \bibinfo{author}{G.~Juve},
  \bibinfo{author}{C.~Kesselman}, \bibinfo{author}{P.~Maechling},
  \bibinfo{author}{G.~Mehta}, \bibinfo{author}{K.~Milner},
  \bibinfo{author}{D.~Okaya}, \bibinfo{author}{P.~Small},
  \bibinfo{author}{K.~Vahi},
\newblock \bibinfo{title}{Cybershake: A physics-based seismic hazard model for
  southern california},
\newblock \bibinfo{journal}{Pure Appl. Geophys.} \bibinfo{volume}{168}
  (\bibinfo{year}{2011}) \bibinfo{pages}{367--381}.
\bibitem[{Rayleigh(1885)}]{Rayleigh1885}
\bibinfo{author}{L.~Rayleigh},
\newblock \bibinfo{title}{On waves propagated along the plane surface of an
  elastic solid},
\newblock \bibinfo{journal}{Proceedings of the London Mathematical Society}
  \bibinfo{volume}{s1-17} (\bibinfo{year}{1885}) \bibinfo{pages}{4--11}.
\bibitem[{Cui et~al.(2013)Cui, Poyraz, Olsen, Zhou, Withers, Callaghan, Larkin,
  Guest, Choi, Chourasia et~al.}]{cui2013sc}
\bibinfo{author}{Y.~Cui}, \bibinfo{author}{E.~Poyraz}, \bibinfo{author}{K.~B.
  Olsen}, \bibinfo{author}{J.~Zhou}, \bibinfo{author}{K.~Withers},
  \bibinfo{author}{S.~Callaghan}, \bibinfo{author}{J.~Larkin},
  \bibinfo{author}{C.~Guest}, \bibinfo{author}{D.~Choi},
  \bibinfo{author}{A.~Chourasia}, et~al.,
\newblock \bibinfo{title}{Physics-based seismic hazard analysis on petascale
  heterogeneous supercomputers},
\newblock in: \bibinfo{booktitle}{High Performance Computing, Networking,
  Storage and Analysis (SC), 2013 International Conference for},
  \bibinfo{organization}{IEEE}, \bibinfo{year}{2013}, pp.
  \bibinfo{pages}{1--12}.
\bibitem[{Rodgers et~al.(2018)Rodgers, Pitarka, Petersson, Sj{\"o}green, and
  McCallen}]{rodgers2018hayward}
\bibinfo{author}{A.~J. Rodgers}, \bibinfo{author}{A.~Pitarka},
  \bibinfo{author}{N.~A. Petersson}, \bibinfo{author}{B.~Sj{\"o}green},
  \bibinfo{author}{D.~B. McCallen},
\newblock \bibinfo{title}{Broadband (0--4 hz) ground motions for a magnitude
  7.0 {H}ayward fault earthquake with three-dimensional structure and
  topography},
\newblock \bibinfo{journal}{Geophysical Research Letters} \bibinfo{volume}{45}
  (\bibinfo{year}{2018}) \bibinfo{pages}{739--747}.
\bibitem[{Reinarz et~al.(2020)Reinarz, Charrier, Bader, Bovard, Dumbser,
  K.~Duru, Gabriel, Gallard, K{\"o}ppel, Krenz, Rannabauer, Rezzolla, Samfass,
  Tavelli, and Weinzierl}]{ExaHyPE2019}
\bibinfo{author}{A.~Reinarz}, \bibinfo{author}{D.~E. Charrier},
  \bibinfo{author}{M.~Bader}, \bibinfo{author}{L.~Bovard},
  \bibinfo{author}{M.~Dumbser}, \bibinfo{author}{F.~F. K.~Duru},
  \bibinfo{author}{A.-A. Gabriel}, \bibinfo{author}{J.-M. Gallard},
  \bibinfo{author}{S.~K{\"o}ppel}, \bibinfo{author}{L.~Krenz},
  \bibinfo{author}{L.~Rannabauer}, \bibinfo{author}{L.~Rezzolla},
  \bibinfo{author}{P.~Samfass}, \bibinfo{author}{M.~Tavelli},
  \bibinfo{author}{T.~Weinzierl},
\newblock \bibinfo{title}{Exahype: An engine for parallel dynamically adaptive
  simulations of wave problems},
\newblock \bibinfo{journal}{Comput. Phys. Comm.}  (\bibinfo{year}{2020})
  \bibinfo{pages}{107251}.
\bibitem[{Dumbser et~al.(2013)Dumbser, Zanotti, Loub\`{e}re, and
  Diot}]{Dumbser:2013}
\bibinfo{author}{M.~Dumbser}, \bibinfo{author}{O.~Zanotti},
  \bibinfo{author}{R.~Loub\`{e}re}, \bibinfo{author}{S.~Diot},
\newblock \bibinfo{title}{{A Posteriori Subcell Limiting of the Discontinuous
  {Galerkin} Finite Element Method for Hyperbolic Conservation Laws}},
\newblock \bibinfo{journal}{J. Comput. Phys.} \bibinfo{volume}{278}
  (\bibinfo{year}{2013}) \bibinfo{pages}{47--75}.
\bibitem[{Weinzierl(2019)}]{Weinzierl2019}
\bibinfo{author}{T.~Weinzierl},
\newblock \bibinfo{title}{Peano -- the peano software--parallel,
  automaton-based, dynamically adaptive grid traversals},
\newblock \bibinfo{journal}{ACM Trans. Math. Softw.} \bibinfo{volume}{45}
  (\bibinfo{year}{2019}). \DOIprefix\doi{https://doi.org/10.1145/3319797}.
\bibitem[{Weinzierl and Mehl(2011)}]{WeinzierlMehl2011}
\bibinfo{author}{T.~Weinzierl}, \bibinfo{author}{M.~Mehl},
\newblock \bibinfo{title}{Peano -- a traversal and storage scheme for
  octree-like adaptive cartesian multi-scale grids},
\newblock \bibinfo{journal}{SIAM J. Sci. Comput.} \bibinfo{volume}{33}
  (\bibinfo{year}{2011}) \bibinfo{pages}{2732--2760}.
\bibitem[{de~la Puente et~al.(2007)de~la Puente, K{\"a}ser, Dumbser, and
  Igel}]{Puente2007}
\bibinfo{author}{J.~de~la Puente}, \bibinfo{author}{M.~K{\"a}ser},
  \bibinfo{author}{M.~Dumbser}, \bibinfo{author}{H.~Igel},
\newblock \bibinfo{title}{{An arbitrary high-order discontinuous Galerkin
  method for elastic waves on unstructured meshes - IV. Anisotropy}},
\newblock \bibinfo{journal}{Geophysical Journal International}
  \bibinfo{volume}{169} (\bibinfo{year}{2007}) \bibinfo{pages}{1210--1228}.
  \DOIprefix\doi{10.1111/j.1365-246X.2007.03381.x}.
\bibitem[{de~la Puente et~al.(2008)de~la Puente, Dumbser, K{\"a}ser, and
  Igel}]{Puente2008}
\bibinfo{author}{J.~de~la Puente}, \bibinfo{author}{M.~Dumbser},
  \bibinfo{author}{M.~K{\"a}ser}, \bibinfo{author}{H.~Igel},
\newblock \bibinfo{title}{{D}iscontinuous {G}alerkin methods for wave
  propagation in poroelastic media},
\newblock \bibinfo{journal}{Geophysics} \bibinfo{volume}{73}
  (\bibinfo{year}{2008}) \bibinfo{pages}{T77--T97}.
  \DOIprefix\doi{10.1190/1.2965027}.
\bibitem[{de~la Puente et~al.(2009)de~la Puente, Ampuero, and
  K{\"a}ser}]{delaPuenteAmpueroKaser2009}
\bibinfo{author}{J.~de~la Puente}, \bibinfo{author}{J.-P. Ampuero},
  \bibinfo{author}{M.~K{\"a}ser},
\newblock \bibinfo{title}{Dynamic rupture modeling on unstructured meshes using
  a discontinuous galerkin method},
\newblock \bibinfo{journal}{J. Geophys. Res.} \bibinfo{volume}{114}
  (\bibinfo{year}{2009}) \bibinfo{pages}{B10302}.
\bibitem[{Pelties et~al.(2012)Pelties, de~la Puente, Ampuero, Brietzke, and
  K{{\"a}}ser}]{PeltiesdelaPuenteAmpueroBrietzkeKaser2012}
\bibinfo{author}{C.~Pelties}, \bibinfo{author}{J.~de~la Puente},
  \bibinfo{author}{J.-P. Ampuero}, \bibinfo{author}{G.~B. Brietzke},
  \bibinfo{author}{M.~K{{\"a}}ser},
\newblock \bibinfo{title}{Three-dimensional dynamic rupture simulation with a
  high-order discontinuous galerkin method on unstructured tetrahedral meshes},
\newblock \bibinfo{journal}{J. Geophys. Res.} \bibinfo{volume}{117}
  (\bibinfo{year}{2012}) \bibinfo{pages}{B02309}.
\bibitem[{Wolf et~al.(2020)Wolf, Gabriel, and Bader}]{Wolf2020}
\bibinfo{author}{S.~Wolf}, \bibinfo{author}{A.~Gabriel},
  \bibinfo{author}{M.~Bader},
\newblock \bibinfo{title}{Optimization and local time stepping of an ader-dg
  scheme for fully anisotropic wave propagation in complex geometries}
  \bibinfo{volume}{12139} (\bibinfo{year}{2020}).
\bibitem[{Grazia et~al.(2013)Grazia, Mengaldo, Moxey, Vincent, and
  Sherwin}]{DeGraziaMengaldoMoxeyVincentSherwin2013}
\bibinfo{author}{D.~D. Grazia}, \bibinfo{author}{G.~Mengaldo},
  \bibinfo{author}{D.~Moxey}, \bibinfo{author}{P.~E. Vincent},
  \bibinfo{author}{S.~Sherwin},
\newblock \bibinfo{title}{Connections between the discontinuous galerkin method
  and high-order flux reconstruction schemes},
\newblock \bibinfo{journal}{Int. J. Numer. Meth. Fluids} \bibinfo{volume}{00}
  (\bibinfo{year}{2013}) \bibinfo{pages}{1--18}.
\bibitem[{Huynh(2007)}]{Huynh2007}
\bibinfo{author}{H.~T. Huynh}, \bibinfo{title}{A flux reconstruction approach
  to high-order schemes including discontinuous galerkin methods},
  \bibinfo{year}{2007}.
\bibitem[{Godunov(1959)}]{Godunov1959}
\bibinfo{author}{S.~K. Godunov},
\newblock \bibinfo{title}{A difference method for numerical calculation of
  discontinuous solutions of the equations of hydrodynamics},
\newblock \bibinfo{journal}{Mat. Sb. (N.S.).} \bibinfo{volume}{181}
  (\bibinfo{year}{1959}) \bibinfo{pages}{271--306}.
\bibitem[{Rusanov(1961)}]{Rusanov1961}
\bibinfo{author}{V.~V. Rusanov},
\newblock \bibinfo{title}{Calculation of interaction of non-stationary shock
  waves with obstacles},
\newblock \bibinfo{journal}{J. Comput. Math. Phys. USSR} \bibinfo{volume}{1}
  (\bibinfo{year}{1961}) \bibinfo{pages}{267--279}.
\bibitem[{Qiu(2008)}]{Qiu2008}
\bibinfo{author}{J.~Qiu},
\newblock \bibinfo{title}{Development and comparison of numerical fluxes for
  lwdg methods},
\newblock \bibinfo{journal}{Numer. Math. Theor. Meth. Appl.}
  \bibinfo{volume}{1} (\bibinfo{year}{2008}) \bibinfo{pages}{435--459}.
\bibitem[{Kirby and Karniadakis(2005)}]{KirbyKarniadakis2005}
\bibinfo{author}{R.~M. Kirby}, \bibinfo{author}{G.~E. Karniadakis},
\newblock \bibinfo{title}{Selecting the numerical flux in discontinuous
  galerkin methods for diffusion problems},
\newblock \bibinfo{journal}{J. Sci. Comput.} \bibinfo{volume}{22}
  (\bibinfo{year}{2005}) \bibinfo{pages}{385--411}.
\bibitem[{Kopriva et~al.(2017)Kopriva, Nordstr{\"o}m, and
  Gassner}]{KoprivaNordstromGassner2016}
\bibinfo{author}{D.~A. Kopriva}, \bibinfo{author}{J.~Nordstr{\"o}m},
  \bibinfo{author}{G.~J. Gassner},
\newblock \bibinfo{title}{Error boundedness of discontinuous galerkin spectral
  element approximations of hyperbolic problems},
\newblock \bibinfo{journal}{J. Sci. Comput.} \bibinfo{volume}{72}
  (\bibinfo{year}{2017}) \bibinfo{pages}{314--330}.
\bibitem[{Duru et~al.(2019)Duru, Rannabauer, Gabriel, and
  Igel}]{DuruGabrielIgel2017}
\bibinfo{author}{K.~Duru}, \bibinfo{author}{L.~Rannabauer},
  \bibinfo{author}{A.-A. Gabriel}, \bibinfo{author}{H.~Igel}, \bibinfo{title}{A
  new discontinuous galerkin method for elastic waves with physically motivated
  numerical fluxes}, \bibinfo{howpublished}{https://arxiv.org/abs/1802.06380},
  \bibinfo{year}{2019}.
\bibitem[{K{\"a}ser and Dumbser(2006)}]{KaeserDumbser2006}
\bibinfo{author}{M.~K{\"a}ser}, \bibinfo{author}{M.~Dumbser},
\newblock \bibinfo{title}{An arbitrary high-order discontinuous {G}alerkin
  method for elastic waves on unstructured {meshes-I. The} two-dimensional
  isotropic case with external source terms},
\newblock \bibinfo{journal}{Geophysical Journal International}
  \bibinfo{volume}{166} (\bibinfo{year}{2006}) \bibinfo{pages}{855--877}.
\bibitem[{Pelties et~al.(2014)Pelties, Gabriel, and Ampuero}]{Pelties2014}
\bibinfo{author}{C.~Pelties}, \bibinfo{author}{A.-A. Gabriel},
  \bibinfo{author}{J.-P. Ampuero},
\newblock \bibinfo{title}{Verification of an {ADER-DG} method for complex
  dynamic rupture problems},
\newblock \bibinfo{journal}{Geoscientific Model Development}
  \bibinfo{volume}{7} (\bibinfo{year}{2014}) \bibinfo{pages}{847--866}.
  \DOIprefix\doi{10.5194/gmd-7-847-2014}.
\bibitem[{Chan and Warburton(2017)}]{ChangWarburton2017}
\bibinfo{author}{J.~Chan}, \bibinfo{author}{T.~Warburton},
\newblock \bibinfo{title}{On the penalty stabilization mechanism for upwind
  discontinuous galerkin formulations of first order hyperbolic systems},
\newblock \bibinfo{journal}{Computs. and Math Applications}
  \bibinfo{volume}{74} (\bibinfo{year}{2017}) \bibinfo{pages}{3099--3110}.
\bibitem[{Wilcox et~al.(2010)Wilcox, Stadler, Burstedde, and
  Ghattas}]{Wilcox2010}
\bibinfo{author}{L.~C. Wilcox}, \bibinfo{author}{G.~Stadler},
  \bibinfo{author}{C.~Burstedde}, \bibinfo{author}{O.~Ghattas},
\newblock \bibinfo{title}{A high-order discontinuous galerkin method for wave
  propagation through coupled elastic-acoustic media},
\newblock \bibinfo{journal}{J. Comput. Phys.} \bibinfo{volume}{229}
  (\bibinfo{year}{2010}) \bibinfo{pages}{9373--9396}.
\bibitem[{Scholz(1998)}]{Scholz1998}
\bibinfo{author}{C.~H. Scholz},
\newblock \bibinfo{title}{Earthquakes and friction laws},
\newblock \bibinfo{journal}{Nature} \bibinfo{volume}{391}
  (\bibinfo{year}{1998}) \bibinfo{pages}{37--42}.
\bibitem[{Rice(1983)}]{Rice1983}
\bibinfo{author}{J.~R. Rice},
\newblock \bibinfo{title}{Constitutive relations for fault slip and earthquake
  instabilities},
\newblock \bibinfo{journal}{J. Appl. Mech.} \bibinfo{volume}{50}
  (\bibinfo{year}{1983}) \bibinfo{pages}{443--475}.
\bibitem[{Rice and Ruina(1983)}]{JRiceetal_83}
\bibinfo{author}{J.~R. Rice}, \bibinfo{author}{A.~L. Ruina},
\newblock \bibinfo{title}{Stability of steady frictional slipping},
\newblock \bibinfo{journal}{J. Appl. Mech.} \bibinfo{volume}{50}
  (\bibinfo{year}{1983}) \bibinfo{pages}{343--349}.
\bibitem[{Duru and Dunham(2016)}]{DuruandDunham2016}
\bibinfo{author}{K.~Duru}, \bibinfo{author}{E.~M. Dunham},
\newblock \bibinfo{title}{Dynamic earthquake rupture simulations on nonplanar
  faults embedded in 3d geometrically complex, heterogeneous elastic solids},
\newblock \bibinfo{journal}{J. Comput. Phys.} \bibinfo{volume}{305}
  (\bibinfo{year}{2016}) \bibinfo{pages}{185--207}.
\bibitem[{Toro(1999)}]{Toro1999}
\bibinfo{author}{E.~F. Toro}, \bibinfo{title}{Riemann solvers and numerical
  methods for fluid dynamics}, \bibinfo{year}{1999}.
\bibitem[{Kristekov{\'a} et~al.(2006)Kristekov{\'a}, Kristek, Moczo, and
  Day}]{Kristekova_etal2006}
\bibinfo{author}{M.~Kristekov{\'a}}, \bibinfo{author}{J.~Kristek},
  \bibinfo{author}{P.~Moczo}, \bibinfo{author}{S.~M. Day},
\newblock \bibinfo{title}{Misfit criteria for quantitative comparison of
  seismograms},
\newblock \bibinfo{journal}{Bull. Seism. Soc. Am.} \bibinfo{volume}{96}
  (\bibinfo{year}{2006}) \bibinfo{pages}{1836--1850}.
\bibitem[{Kristekov{\'a} et~al.(2009)Kristekov{\'a}, Kristek, and
  Moczo}]{Kristekova_etal2009}
\bibinfo{author}{M.~Kristekov{\'a}}, \bibinfo{author}{J.~Kristek},
  \bibinfo{author}{P.~Moczo},
\newblock \bibinfo{title}{Time-frequency misfit and goodness-of-fit criteria
  for quantitative comparison of time signals},
\newblock \bibinfo{journal}{Geophys. J. Int.} \bibinfo{volume}{178}
  (\bibinfo{year}{2009}) \bibinfo{pages}{813--825}.
\bibitem[{Favretto-Cristini et~al.(2011)Favretto-Cristini, Komatitsch,
  Carcione, and Cavallini}]{Favretto-Cristini_etal2011}
\bibinfo{author}{N.~Favretto-Cristini}, \bibinfo{author}{D.~Komatitsch},
  \bibinfo{author}{J.~M. Carcione}, \bibinfo{author}{F.~Cavallini},
\newblock \bibinfo{title}{{Elastic surface waves in crystals. Part 1: Review of
  the physics}},
\newblock \bibinfo{journal}{Ultrasonics} \bibinfo{volume}{51}
  (\bibinfo{year}{2011}) \bibinfo{pages}{653--660}.
\bibitem[{Komatitsch et~al.(2011)Komatitsch, Carcione, Cavallini, and
  Favretto-Cristini}]{Komatitsch_etal2011}
\bibinfo{author}{D.~Komatitsch}, \bibinfo{author}{J.~M. Carcione},
  \bibinfo{author}{F.~Cavallini}, \bibinfo{author}{N.~Favretto-Cristini},
\newblock \bibinfo{title}{{Elastic surface waves in crystals -- Part 2:
  Cross-check of two full-wave numerical modeling methods}},
\newblock \bibinfo{journal}{Ultrasonics} \bibinfo{volume}{51}
  (\bibinfo{year}{2011}) \bibinfo{pages}{878--889}.
\bibitem[{Duru et~al.(2020)Duru, Rannabauer, Gabriel, Kreiss, and
  Bader}]{DuruRannabauerGabrielKreissBader2019}
\bibinfo{author}{K.~Duru}, \bibinfo{author}{L.~Rannabauer},
  \bibinfo{author}{A.-A. Gabriel}, \bibinfo{author}{G.~Kreiss},
  \bibinfo{author}{M.~Bader},
\newblock \bibinfo{title}{A stable discontinuous galerkin method for the
  perfectly matched layer for elastodynamics in first order form},
\newblock \bibinfo{journal}{Numerische Mathematik} \bibinfo{volume}{146}
  (\bibinfo{year}{2020}) \bibinfo{pages}{729--782}.
  \DOIprefix\doi{https://doi.org/10.1007/s00211-020-01160-w}.
\bibitem[{Duru et~al.(2019)Duru, Gabriel, and Kreiss}]{Duru2019}
\bibinfo{author}{K.~Duru}, \bibinfo{author}{A.-A. Gabriel},
  \bibinfo{author}{G.~Kreiss},
\newblock \bibinfo{title}{On energy stable discontinuous {G}alerkin spectral
  element approximations of the perfectly matched layer for the wave equation},
\newblock \bibinfo{journal}{Computer Methods in Applied Mechanics and
  Engineering} \bibinfo{volume}{350} (\bibinfo{year}{2019})
  \bibinfo{pages}{898--937}.
\bibitem[{Marsden and Hughes(1994)}]{MarsdenHughes1994}
\bibinfo{author}{J.~E. Marsden}, \bibinfo{author}{T.~J.~R. Hughes},
  \bibinfo{title}{Mathematical Foundations of Elasticity},
  \bibinfo{publisher}{Dover Publications Inc. New York}, \bibinfo{year}{1994}.
\bibitem[{Nordstr{\"o}m(2006)}]{Norstrom2006}
\bibinfo{author}{J.~Nordstr{\"o}m},
\newblock \bibinfo{title}{Conservative finite difference formulations, variable
  coefficients, energy estimates and artificial dissipation},
\newblock \bibinfo{journal}{J. Sci. Comput.} \bibinfo{volume}{29}
  (\bibinfo{year}{2006}) \bibinfo{pages}{375--404}.
\bibitem[{Kopriva(2006)}]{Kopriva2006}
\bibinfo{author}{D.~Kopriva},
\newblock \bibinfo{title}{Metric identities and the discontinuous spectral
  element method on curvilinear meshes},
\newblock \bibinfo{journal}{J. Sci. Comput.} \bibinfo{volume}{26}
  (\bibinfo{year}{2006}) \bibinfo{pages}{301--327}.
\bibitem[{Thomas and Lombard(1979)}]{ThomasLombad1979}
\bibinfo{author}{P.~Thomas}, \bibinfo{author}{C.~Lombard},
\newblock \bibinfo{title}{Geometric conservation law and its application to
  flow computations on moving grids},
\newblock \bibinfo{journal}{AIAA J.} \bibinfo{volume}{17}
  (\bibinfo{year}{1979}) \bibinfo{pages}{1030--1037}.
\bibitem[{Kopriva and Gassner(2014)}]{KoprivaGassner2014}
\bibinfo{author}{D.~A. Kopriva}, \bibinfo{author}{G.~J. Gassner},
\newblock \bibinfo{title}{An energy stable discontinuous galerkin spectral
  element discretization for variable coefficient advection problems},
\newblock \bibinfo{journal}{SIAM J. Sci. Comput.} \bibinfo{volume}{36}
  (\bibinfo{year}{2014}) \bibinfo{pages}{A2076--A2099}.
\bibitem[{Kopriva and Gassner(2015)}]{KoprivaGassner2015}
\bibinfo{author}{D.~Kopriva}, \bibinfo{author}{G.~Gassner},
\newblock \bibinfo{title}{Geometry effects in nodal discontinuous galerkin
  methods on curved elements that are provably stable},
\newblock \bibinfo{journal}{Appl. Math. Comput.}  (\bibinfo{year}{2015})
  \bibinfo{pages}{1--17}.
\bibitem[{Gustafsson et~al.(1995)Gustafsson, Kreiss, and
  Oliger}]{GustafssonKreissOliger1995}
\bibinfo{author}{B.~Gustafsson}, \bibinfo{author}{H.-O. Kreiss},
  \bibinfo{author}{J.~Oliger}, \bibinfo{title}{Time dependent problems and
  difference methods}, \bibinfo{publisher}{John Wiley and Sons, New York},
  \bibinfo{year}{1995}.
\bibitem[{Achenbach(1973)}]{Achenbach1973}
\bibinfo{author}{J.~Achenbach}, \bibinfo{title}{Wave propagation in elastic
  solids}, volume~\bibinfo{volume}{6}, \bibinfo{publisher}{Applied Mathematics
  and Mechanics. North-Holland}, \bibinfo{year}{1973}.
\bibitem[{Petersson et~al.(2016)Petersson, O'Reilly, Sj{\"o}green, and
  Bydlon}]{Petersson_etal2016}
\bibinfo{author}{N.~A. Petersson}, \bibinfo{author}{O.~O'Reilly},
  \bibinfo{author}{B.~Sj{\"o}green}, \bibinfo{author}{S.~Bydlon},
\newblock \bibinfo{title}{Discretizing singular point sources in hyperbolic
  wave propagation problems},
\newblock \bibinfo{journal}{J. Comput. Phys.} \bibinfo{volume}{321}
  (\bibinfo{year}{2016}) \bibinfo{pages}{532--555}.
\bibitem[{Duru et~al.(2020)Duru, Fung, and Williams}]{DuruFungWilliams2020}
\bibinfo{author}{K.~Duru}, \bibinfo{author}{F.~Fung},
  \bibinfo{author}{C.~Williams}, \bibinfo{title}{Upwind summation by parts
  finite difference methods for large scale elastic wave simulations in complex
  geometries}, \bibinfo{howpublished}{https://arxiv.org/abs/2011.02600},
  \bibinfo{year}{2020}.
\bibitem[{Krischer et~al.(2015)Krischer, Megies, Barsch, Beyreuther, Lecocq,
  Caudron, and Wassermann}]{ObsPy2015}
\bibinfo{author}{L.~Krischer}, \bibinfo{author}{T.~Megies},
  \bibinfo{author}{R.~Barsch}, \bibinfo{author}{M.~Beyreuther},
  \bibinfo{author}{T.~Lecocq}, \bibinfo{author}{C.~Caudron},
  \bibinfo{author}{J.~Wassermann},
\newblock \bibinfo{title}{Obspy: a bridge for seismology into the scientific
  python ecosystem},
\newblock \bibinfo{journal}{Computational Science \& Discovery}
  \bibinfo{volume}{8} (\bibinfo{year}{2015}) \bibinfo{pages}{014003}.
  \DOIprefix\doi{https://doi.org/10.1088/1749-4699/8/1/014003}.
\bibitem[{Het{\'e}nyi et~al.(2018)Het{\'e}nyi, Molinari, Clinton, Bokelmann,
  Bond{\'a}r, Crawford, Dessa, Doubre, Friederich, Fuchs et~al.}]{Hetenyi2018}
\bibinfo{author}{G.~Het{\'e}nyi}, \bibinfo{author}{I.~Molinari},
  \bibinfo{author}{J.~Clinton}, \bibinfo{author}{G.~Bokelmann},
  \bibinfo{author}{I.~Bond{\'a}r}, \bibinfo{author}{W.~C. Crawford},
  \bibinfo{author}{J.-X. Dessa}, \bibinfo{author}{C.~Doubre},
  \bibinfo{author}{W.~Friederich}, \bibinfo{author}{F.~Fuchs}, et~al.,
\newblock \bibinfo{title}{The {AlpArray} seismic network: a large-scale
  {E}uropean experiment to image the {A}lpine orogen},
\newblock \bibinfo{journal}{Surveys in geophysics} \bibinfo{volume}{39}
  (\bibinfo{year}{2018}) \bibinfo{pages}{1009--1033}.
\bibitem[{Oeser et~al.(2016)Oeser, Bunge, and Mohr}]{Oeser2006}
\bibinfo{author}{J.~Oeser}, \bibinfo{author}{H.-P. Bunge},
  \bibinfo{author}{M.~Mohr}, \bibinfo{title}{Cluster design in the earth
  sciences: Tethys}, \bibinfo{howpublished}{in High Performance Computing and
  Communications - Second International Conference, HPCC 2006, Munich, Germany,
  Lecture Notes in Computer Science}, \bibinfo{year}{2016}.
  \DOIprefix\doi{https://doi.org/10.1007/118473664}.

\end{thebibliography}

\end{document}